\documentclass[11pt]{article}
\usepackage{amssymb,amsfonts,amsmath,curves}

\usepackage[english]{babel}
\def \beq {\begin{eqnarray}}
\def \eeq {\end{eqnarray}}
\def \beqn {\begin{eqnarray*}}
\def \eeqn {\end{eqnarray*}}
\def \nn {\nonumber}

 \makeatletter
 \@addtoreset{equation}{section}
 \makeatother

 \makeatletter
 \@addtoreset{enunciato}{section}
 \makeatother

 \newcounter{enunciato}[section]

 \newtheorem{ittheorem}{Theorem}
 \newtheorem{itlemma}{Lemma}
 \newtheorem{itproposition}{Proposition}
 \newtheorem{itdefinition}{Definition}
 \newtheorem{itremark}{Remark}
 \newtheorem{itclaim}{Claim}
 \newtheorem{itfact}{Fact}
 \newtheorem{itcorollary}{Corollary}
 \newtheorem{itconjecture}{Conjecture}

 \newenvironment{theorem}{\addtocounter{enunciato}{1}
 \begin{ittheorem}}{\end{ittheorem}}

 \newenvironment{lemma}{\addtocounter{enunciato}{1}
 \begin{itlemma}}{\end{itlemma}}

 \newenvironment{proposition}{\addtocounter{enunciato}{1}
 \begin{itproposition}}{\end{itproposition}}

 \newenvironment{definition}{\addtocounter{enunciato}{1}
 \begin{itdefinition}}{\end{itdefinition}}

 \newenvironment{remark}{\addtocounter{enunciato}{1}
 \begin{itremark}}{\end{itremark}}

 \newenvironment{claim}{\addtocounter{enunciato}{1}
 \begin{itclaim}}{\end{itclaim}}

 \newenvironment{fact}{\addtocounter{enunciato}{1}
 \begin{itfact}}{\end{itfact}}

 \newenvironment{corollary}{\addtocounter{enunciato}{1}
 \begin{itcorollary}}{\end{itcorollary}}

 \newenvironment{conjecture}{\addtocounter{enunciato}{1}
 \begin{itconjecture}}{\end{itconjecture}}

 \newcommand{\be}[1]{\begin{equation}\label{#1}}
 \newcommand{\ee}{\end{equation}}

 \newcommand{\bl}[1]{\begin{lemma}\label{#1}}
 \newcommand{\el}{\end{lemma}}

 \newcommand{\br}[1]{\begin{remark}\label{#1}}
 \newcommand{\er}{\end{remark}}

 \newcommand{\bt}[1]{\begin{theorem}\label{#1}}
 \newcommand{\et}{\end{theorem}}

 \newcommand{\bd}[1]{\begin{definition}\label{#1}}
 \newcommand{\ed}{\end{definition}}

 \newcommand{\bcl}[1]{\begin{claim}\label{#1}}
 \newcommand{\ecl}{\end{claim}}

 \newcommand{\bfact}[1]{\begin{fact}\label{#1}}
 \newcommand{\efact}{\end{fact}}

 \newcommand{\bp}[1]{\begin{proposition}\label{#1}}
 \newcommand{\ep}{\end{proposition}}

 \newcommand{\bc}[1]{\begin{corollary}\label{#1}}
 \newcommand{\ec}{\end{corollary}}

 \newcommand{\bcon}[1]{\begin{conjecture}\label{#1}}
 \newcommand{\econ}{\end{conjecture}}

 \newcommand{\bpr}{\begin{proof}}
 \newcommand{\epr}{\end{proof}}

 \newcommand{\bi}{\begin{itemize}}
 \newcommand{\ei}{\end{itemize}}

 \newcommand{\ben}{\begin{enumerate}}
 \newcommand{\een}{\end{enumerate}}

 \newenvironment{proof}{\noindent {\em Proof}.\,\,}{\hspace*{\fill}$\halmos$\medskip}

 \newcommand{\halmos}{\rule{1ex}{1.4ex}}
 \def \qed {{\hspace*{\fill}$\halmos$\medskip}}

\def \R {\mathbb R}

\def \N {\mathbb N}
\def \E {\mathbb E}
\def \P {\mathbb P}

\def \cH {\mathcal H}
\def \cC {\mathcal C}
\def \cL {\mathcal L}

\renewcommand{\O}{\Omega}

\newcommand{\Ga}{\Gamma}
\newcommand{\tet}{\theta}
\newcommand{\la}{\lambda}
\newcommand{\Gi}{\mathcal G}

\begin{document}

\title{The renormalization transformation for\\
two-type branching models}
\author{
D.A.\ Dawson$^{\, 1}$ \\
A.\ Greven$^{\, 2}$ \\
F.\ den Hollander$^{\, 3,4}$ \\
Rongfeng Sun$^{\,4,5}$ \\
J.M.\ Swart$^{\, 6}$
}
\date{May 16th, 2007}
\maketitle
\footnotetext[1]{School of Mathematics and Statistics, Carleton
University, Ottawa K1S 5B6, Canada, {\em ddawson@math.carleton.ca}}
\footnotetext[2]{Mathematisches Institut, Universit\"at Erlangen-N\"urnberg,
Bismarckstra{\ss}e 1 1/2, D-91054 Erlangen, Germany, {\em greven@mi.uni-erlangen.de}}
\footnotetext[3]{Mathematical Institute, Leiden University, P.O.\ Box 9512,
2300 RA Leiden, the Netherlands, {\em denholla@math.leidenuniv.nl}}\,
\footnotetext[4]{EURANDOM, P.O.\ Box 513, 5600 MB Eindhoven, the Netherlands}\,
\footnotetext[5]{MA 7-5, Fakult\"at II -- Institut f\"ur Mathematik, TU Berlin, Stra\ss e des 17.\ Juni 136,
10623 Berlin, {\em sun@math.tu-berlin.de}}\,
\footnotetext[6]{\' UTIA, Pod vod\'arenskou v\v e\v z\' i 4, 18208 Praha 8,
Czech Republic, {\em swart@utia.cas.cz}}

\vspace{-20pt}

\begin{abstract}
This paper studies countable systems of linearly and hierarchically interacting diffusions
taking values in the positive quadrant. These systems arise
in population dynamics for two types of individuals migrating between and interacting
within colonies. Their large-scale space-time behavior can be studied by means of a
renormalization program. This program, which has been carried out
successfully in a number of other cases (mostly one-dimensional), is based on the
construction and the analysis of a nonlinear renormalization transformation, acting on
the diffusion function for the components of the system and connecting the evolution
of successive block averages on successive time scales. We identify a general class
of diffusion functions on the positive quadrant for which this renormalization
transformation is well-defined and, subject to a conjecture on its boundary behavior,
can be iterated. Within certain subclasses, we identify the fixed points for the transformation
and investigate their domains of attraction. These domains of attraction
constitute the universality classes of the system under space-time scaling.

\bigskip\noindent
\emph{Keywords:} Interacting diffusions, space-time renormalization,
two-type populations, independent branching, catalytic branching,
mutually catalytic branching, universality.
\bigskip

\noindent
\emph{AMS 2000 subject classification:} 60J60, 60J70, 60K35.
\bigskip

\end{abstract}

\newpage


\section{Introduction}
\label{S1}

\subsection{Model and background}
\label{S1.1}

We are interested in the following system of coupled stochastic differential
equations (SDE):
\be{SDEcoupled}
dX_{\eta, i}(t)
= \sum_{\xi\in\O_N} a_N(\xi,\eta)\,[X_{\xi,i}(t) - X_{\eta,i}(t)]\,dt
+ \sqrt{2g_i(\vec X_{\eta}(t))}\, dB_{\eta,i}(t), \quad \eta\in\O_N, i=1,2.
\ee
Here $a_N(\cdot,\cdot)$ is the transition rate kernel of a random walk on $\O_N$, the hierarchical group
(or lattice) of order $N$ (see (\ref{hier})),
$\{\vec X_\eta\}_{\eta\in\O_N}$ with $\vec X_\eta = (X_{\eta, 1}, X_{\eta, 2})$
is a family of diffusions taking values in $[0,\infty)^2$, $g=(g_1, g_2)$ is a pair
of diffusion functions on $[0,\infty)^2$, and $\{\vec B_\eta\}_{\eta\in\O_N}$ with
$\vec B_\eta = (B_{\eta,1}, B_{\eta,2})$ is a family of independent standard Brownian
motions on $\R^2$. As initial condition, we
take
\be{initcond}
\vec X_\eta(0) = \vec\theta
= (\theta_1, \theta_2)\in[0,\infty)^2 \qquad \forall\,\eta\in\O_N.
\ee

Equation (\ref{SDEcoupled}) arises as the continuum limit of discrete models in {\em
population dynamics}. In these models, individuals live in colonies labeled by the
hierarchical group $\O_N$. Each colony $\eta\in\O_N$ consists of two types of individuals,
whose total masses are represented by the vector $\vec X_\eta$. Individuals {\em migrate}
between colonies according to the migration kernel $a_N(\cdot,\cdot)$. At each colony,
each individual undergoes {\em branching} at a rate
that depends on the total masses of the two types of individuals present at that colony.
The system in (\ref{SDEcoupled}) arises in the so-called ``small-mass-fast-branching''
limit, where the number of individuals in each colony tends to infinity, the mass of
each individual tends to zero, and the effective branching rate grows proportionally to
the number of individuals in each colony. The drift term in (\ref{SDEcoupled}) arises
from the migration, which is the only source of interaction \emph{between} colonies. The
diffusion term in (\ref{SDEcoupled}) arises from the branching, where $g_i(x)/x_i$ is
the state-dependent branching rate of the $i$-th type, which incorporates the
interaction between individuals \emph{within} a colony. For more background, see e.g.\
Sawyer and Felsenstein \cite{SF83}, Dawson and Perkins \cite{DP98}, Chapters 9--10 in
Ethier and Kurtz \cite{EK86}, Cox, Dawson and Greven \cite{CDG04}, Dawson, Gorostiza
and Wakolbinger \cite{DGW05}.

The goal of the present paper is to study the universality classes of
the large-scale space-time behavior of (\ref{SDEcoupled}). It turns
out that, for the specific form of the migration kernel
$a_N(\cdot,\cdot)$ given by (\ref{rwkernel}) below and in the limit as
$N\to\infty$, (\ref{SDEcoupled}) is susceptible to a renormalization
analysis. The {\em renormalization program} for hierarchically
interacting diffusions was introduced by Dawson and Greven
\cite{DG93a}, \cite{DG93b} for diffusions taking values in $[0,1]$. It
has since been extended to several other state spaces (see Greven
\cite{G05} for an overview). We will give more detailed references in
Section \ref{S1.3}. First we outline the main ingredients of the
renormalization program.

\subsection{Renormalization program}
\label{S1.2}

The lattice in (\ref{SDEcoupled}) is the {\em hierarchical group} of order $N$,
which is defined as
\be{hier}
\Omega_N = \left\{\eta=(\eta_i)_{i\in\N}\in \{0,1,\dots,N-1\}^\N\colon\,
\sum_{i\in\N} \eta_i<\infty\right\}
\ee
with coordinatewise addition modulo $N$. Define a shift
$\phi:\O_N\to\O_N$ by $(\phi\eta)_i:=\eta_{i+1}$ $(i\in\N)$. On
$\Omega_N$, the {\em hierarchical distance} is defined as
\be{1.4}
d(\eta,\xi) = \min\{k\in\N_0=\N\cup\{0\}\colon\,\phi^k\eta=\phi^k\xi\},
\ee
which is an {\em ultrametric}, i.e., $d(\eta,\xi) \leq d(\eta,\zeta) \vee
d(\xi,\zeta)$ for all $\eta,\xi,\zeta\in\O_N$. We choose the random walk
transition rate kernel in such a way that $a_N(\xi,\eta)$ depends
only on the hierarchical distance between $\xi$ and $\eta$. In view of what follows,
we write $a_N$ in the form
\be{rwkernel}
a_N(\xi,\eta) = \sum_{k \geq d(\xi,\eta)} c_{k-1}\, N^{1-2k}, \qquad
\xi,\eta\in\O_N,\,\xi \neq \eta,
\ee
where $(c_n)_{n\in \N_0}$ is a sequence of positive constants. Formula
(\ref{rwkernel}) says that the random walk associated with $a_N(\cdot, \cdot)$ jumps
with rate $c_{k-1}/N^{k-1}$ from $\eta$ to an arbitrary site in the
{\em $k$-block} $\{\xi\in\O_n:\phi^k\xi=\phi^k\eta\}$ around $\eta$.

The key objects in the renormalization analysis are the {\em $k$-block averages}:
\be{1.5}
Y_{\eta,i}^{[k]}(t) = \frac{1}{N^k} \sum\limits_{{\xi\in\Omega_N}
\atop {\phi^k\xi=\eta}} X_{\xi, i} (t),
\qquad \eta\in\Omega_N,\,i=1,2,\,k\in\N_0.
\ee
Using (\ref{rwkernel}), we may rewrite (\ref{SDEcoupled}) as
\be{SDEcoupled2}
\begin{aligned}
\!\!\!\! dX_{\eta, i}(t) = \sum_{k\geq 1} \frac{c_{k-1}}{N^{k-1}}
\left[Y_{\phi^k\eta,i}^{[k]}(t) - X_{\eta, i}(t) \right]\!dt
+ \sqrt{2g_i(\vec X_\eta(t))}\,\, dB_{\eta,i}(t), \quad \eta\in\O_N,\,i=1,2,
\end{aligned}
\ee
where each component $\vec X_{\eta}$ feels a drift towards the successive
averages of $k$-blocks containing $\eta$. It can be seen that the evolution of the
$1$-block averages is described in law by the SDE
\be{SDE1block}
\begin{aligned}
dY^{[1]}_{\eta, i}(tN)
&= \sum_{k\geq 1} \frac{c_{k}}{N^{k-1}}
\left[Y^{[k+1]}_{\phi^k\eta,i}(tN) - Y^{[1]}_{\eta,i}(tN)\right]\,dt\\
&\qquad + \sqrt{\frac{2}{N}\sum_{{\xi\in\O_N}\atop{\phi\xi=\eta}}
g_i(\vec X_\xi(tN))}\,dB_{\eta,i}(t),\qquad \eta\in\O_N,\,i=1,2,
\end{aligned}
\ee
where $\vec B_\eta=(B_{\eta,1},B_{\eta,2})$ is a
family of independent standard two-dimensional Brownian motions. Note
that in the limit $N\to\infty$, we expect both the drift and the
diffusion term in (\ref{SDE1block}) to be of order one, which means
that $\vec Y^{[1]}_\eta$ evolves on the time scale $tN$.

Let us next see {\em heuristically} what happens if we let $N\to\infty$, the so-called
{\em hierarchical mean-field limit}. If we let $N\to\infty$ in (\ref{SDEcoupled2}),
then the only drift term that survives is
$$
c_0\big[Y^{[1]}_{\phi\eta,i}(t) - X_{\eta,i}(t)\big]dt.
$$
Furthermore, $\vec Y^{[1]}_{\phi\eta}(t) \to \vec X_{(\cdot)}(0) \equiv \vec \theta$ for all $t\geq 0$,
because $\vec Y^{[1]}_{\phi\eta}$ evolves on the time scale $tN$. Therefore the system
$\{\vec X_\eta(t)\}_{\eta\in\O_N}$ converges in law to an \emph{independent} system of
diffusions, each satisfying the autonomous SDE
\be{sde0}
dZ_i(t) = c_0 (\theta_i - Z_i)\ dt + \sqrt{2g_i(\vec Z(t))}\ dB_i(t),
\qquad i=1,2.
\ee
This kind of behavior is frequently referred to as ``McKean-Vlasov
limit'' and ``propagation of chaos''.

With the above fact in mind, we move one step up in the hierarchy. Since $\vec X_\xi(t)$ evolves on the time scale $t$, for each fixed $t$ the family
\be{Xfamily}
\{\vec X_\xi(tN)\}_{{\xi\in\O_N}\atop{\phi\xi=\eta}}
\ee
decouples and converges almost instantly to the equilibrium distribution of
(\ref{sde0}) with the drift towards $\vec\theta$ replaced by a drift towards the
first block average $\vec Y^{[1]}_\eta(tN)$. Thus, we expect that
\be{conbl}
\frac{1}{N}\sum_{{\xi\in\O_N}\atop{\phi\xi=\eta}}g_i(\vec X_\xi(tN))
\sim \int_{[0,\infty)^2}\Ga^{c_0,g}_{\vec Y^{[1]}_\eta(tN)}(d\vec x)g_i(\vec x)
\quad \mbox{ as } N\to\infty \mbox{ for fixed } t,
\quad \eta\in\O_N,\,i=1,2,
\ee
where $\Ga^{c_0,g}_{\vec\tet}$ denotes the {\em equilibrium distribution} of (\ref{sde0}).
Thus, if we set
\be{Fcg}
(F_{c_0}g)_i(\vec\theta)=\int_{[0,\infty)^2}\Ga^{c_0,g}_{\vec\tet}(d\vec x)g_i(\vec x),
\qquad i=1,2,\ \vec\tet\in[0,\infty)^2,
\ee
then by (\ref{conbl}), for large $N$, the SDE (\ref{SDE1block}) for
the $1$-block averages $\vec Y^{[1]}_\eta$ takes exactly the same form
as the SDE (\ref{SDEcoupled2}) for the single components, provided
that we rescale time by a factor $N$ and replace the single component
diffusion functions $g_i$ by $(F_{c_0}g)_i$ $(i=1,2$). Here, $F_{c_0}$
plays the role of a {\em renormalization transformation} acting on the
pair of diffusion functions $g=(g_1,g_2)$.

We can {\em iterate} the above procedure. The upshot of this is that,
as $N\to\infty$, the $k$-block averages $\vec Y^{[k]}_\eta$ evolve on the time scale
$tN^k$ according to the SDE
\be{sdek}
dZ^{[k]}_i(t) = c_k \left(\theta_i - Z^{[k]}_i(t)\right)\,dt
+ \sqrt{2(F^{[k]}g)_i(\vec Z^{[k]}(t))}\,\,dB_i(t), \qquad i=1,2,
\ee
with diffusion functions $F^{[k]}g=(F^{[k]}g_1,F^{[k]}g_2)$ given by
\be{Fkdef}
F^{[k]} g = F_{c_{k-1}}\circ\cdots \circ F_{c_0}\, g, \qquad k\in \N_0.
\ee
In fact, putting the successive iterates together and observing the sequence of block
averages
\be{blav}
\Big(\vec Y_{\phi^k\eta}^{[k]}(sN^k), \vec Y_{\phi^{k-1}\eta}^{[k-1]}(sN^k), \cdots, \vec Y_\eta^{[0]}(sN^k)\Big)
\ee
on the time scale $sN^k$, as $N\to\infty$, we expect this sequence to converge
in distribution to a \emph{backward Markov chain}
\be{MCback}
\big(\vec M(-k), \vec M(-k+1), \cdots, \vec M(0)\big),
\ee
the so-called {\it interaction chain}, where
\begin{itemize}
\item[(1)]
The starting position $\vec M(-k)$ is distributed as the weak solution of (\ref{sdek}) at time
$s$ with initial condition $\vec Z^{[k]}(0) = \vec\theta$;
\item[(2)]
for $0\leq j \leq k-1$, the transition probability kernel from $\vec M(-j-1)$ to $\vec M(-j)$
is given by
\be{Mtrans}
\P\big[\vec M(-j)\in d\vec y\,\big|\,\vec M(-j-1)=\vec x\big]=\Ga^{c_j, F^{[j]}g}_{\vec x}(d\vec y),
\ee
where $\Ga^{c_j, F^{[j]}g}_{\vec x}(\,\cdot\,)$ denotes the equilibrium distribution of (\ref{sdek}) with $k$
replaced by $j$.
\end{itemize}
The distribution of $\vec M(-k)$ depends on $s$ because $\vec Y_{\phi^k\eta}^{[k]}(sN^k)$ evolves on the time
scale $sN^k$, while the transition probability kernel from $\vec M(-j-1)$ to $\vec M(-j)$ for $0\leq j \leq k-1$ is
independent of $t$ because, conditioned on $\vec Y_{\phi^{j+1}\eta}^{[j+1]}$, $\vec Y^{[j]}_{\phi^j\eta}$ equilibrates almost
instantly on the time scale $sN^k$. Note that $(F^{[k]}g)_i(\vec\theta) =\E[g_i(\vec M(0))\, |\, \vec M(-k) =
\vec\theta]$, where $\E$ denotes expectation with respect to the interaction chain.

With these heuristics in mind, the {\em renormalization program} consists of the following two steps:
\begin{itemize}
\item[(I)]
\underline{Stochastic part:}
Show that for all scales $k\in\N$, in the hierarchical mean-field limit $N\to\infty$,
the block average in (\ref{1.5}) converges in law to the solution of the SDE in (\ref{sdek}),
and the sequence of block averages in (\ref{blav})  converges in law to the interaction
chain in (\ref{MCback}).
\item[(II)]
\underline{Analytic part:} Analyze the renormalization transformation $F_c$ and the iterates
$F^{[n]}$,  $n\in\N_0$.
\end{itemize}

Assuming that the stochastic part of the renormalization program can
be completed, the large-scale space-time behavior of
(\ref{SDEcoupled}) in the limit $N\to\infty$ is
characterized by the behavior of $F^{[n]}$ as $n\to\infty$, in
particular, by its {\em fixed shapes} and their {\em universality classes}.

Here, by a fixed shape we mean a pair of diffusion functions
$g=(g_1,g_2)$ such that $F_cg=\la g$ for some $c,\la>0$. We speak of a
{\em downgoing fixed shape, fixed point} or {\em upgoing fixed shape}
depending on whether $\la<1,=1$, or $>1$. Note that since the factor
$\la$ can always be absorbed in time-scaling, such fixed shapes
correspond to models that are mapped into themselves after a suitable
rescaling of space and time. Indeed, if we set $c_k=c\lambda^k$ $(k\geq 0)$,
then such a fixed shape satisfies $F^{[k]}g=\la^kg$ because the SDE associated
with $(c_k, F^{[k]}g)$ is simply a time change of the SDE associated with $(c,g)$,
which induces the same renormalization transformation. For the
interacting model in (\ref{SDEcoupled2}), this means that the
$k$-block averages evolve on the time scale $tN^k\la^k$ according to
the diffusion function $g$. We note that our definition of a fixed
shape deviates from the definition used in some earlier work, e.g.\ Fleischmann and Swart~\cite{FS06}. What is called a fixed shape there is, in our
terminology, a {\em joint fixed shape} for all $c>0$, i.e., a $g$ such
that for all $c>0$ there exists a $\la=\la(c)$ with $F_cg=\la g$.

By a universality class, we mean a set $\Gi$ of diffusion functions with the property
that, given $(c_k)_{k\in\N_0}$, for each $g\in\Gi$ there exist
scaling constants $(s_n)_{n\in\N}$ such that $s_nF^{[n]}g$ converges to the same limit (possibly up to a
multiplicative constant). Typically, the limit will be a fixed shape
or an {\em asymptotic fixed shape} (for the latter, see Fleischmann and Swart~\cite{FS06}).
Note that each joint fixed shape gives rise to a universality class, namely all models
within a given universality class exhibit the same large-scale  space-time behavior.

Apart from being relevant in the study of large-scale space-time behavior,
fixed shapes also give rise to continuum models, by taking the
so-called {\em hierarchical mean-field continuum limit}, which is a
spatial continuum limit of the hierarchical lattice $\O_N$ with
$N\to\infty$. These continuum models also exhibit universality on
small space-time scales, which is governed by the same renormalization
transformation $F_c$ and its iterates $F^{[n]}$, $n\in\N_0$. For more details, see
Cox, Dawson and Greven \cite{CDG04}, and Dawson, Greven and
Z\"ahle~\cite{DGZ}.

The large-scale space-time behavior of (\ref{SDEcoupled}) depends both
on the diffusion function $g$ and on the potential-theoretic
properties of the random walk with transition rate kernel
(\ref{rwkernel}). Based on earlier work, we expect nontrivial
universality classes to arise only when $\sum_{n\in\N_0} c_n^{-1}
=\infty$, which is the ``necessary and sufficient'' condition for the
random walk with transition rate kernel $a_N(\cdot, \cdot)$ on $\O_N$ to be
recurrent (except for a side condition that becomes irrelevant in the
limit $N\to\infty$; see Sawyer and Felsenstein~\cite{SF83}). For
linear systems such as (\ref{SDEcoupled}), the recurrence of the
random walk is usually associated with {\em clustering}; see e.g.\
Dawson and Greven~\cite{DG93b}, Cox and Greven~\cite{CG94},
Swart~\cite{S00}. In our context, clustering means that the solution
of (\ref{SDEcoupled}) converges in law to a mixture of distributions,
each of which is concentrated on the configuration $\vec X_\eta= \vec
x$, $\eta\in\O_N$, for some $\vec x\in [0,\infty)^2$ with $g_1(\vec
x)= g_2(\vec x)=0$. The choice of $(c_n)_{n\in \N_0}$ determines the
pattern of cluster formation, such as whether only small clusters appear, or only large
clusters appear, or clusters of all scales appear. The latter is known as  {\em diffusive clustering} (see e.g.\
Dawson and Greven~\cite{DG93b}).

With the above facts in mind, the {\em analytic part} of the renormalization
program can be more precisely formulated as follows.
\begin{itemize}
\item[1.] Find classes of diffusion functions on which the renormalization transformations
$F_c$ and their iterates $F^{[n]}$, $n\in\N_0$, are well-defined.
\item[2.] Determine all the (asymptotic) fixed shapes.
\item[3.] Determine the universality classes of diffusion functions that, for
given $(c_n)_{n\in\N_0}$ and after appropriate rescaling, converge to
these (asymptotic) fixed shapes, and determine the associated scaling constants.
\end{itemize}

\subsection{Literature}
\label{S1.3}

The {\em full} renormalization program has been successfully carried out for
hierarchically interacting diffusions taking values in:
\begin{itemize}
\item[(1)]
the compact interval $[0,1]$ (Dawson and Greven \cite{DG93a}, \cite{DG93b}, Baillon,
Cl\'ement, Greven and den Hollander \cite{BCGH95}), where the Wright-Fisher diffusion
is the unique fixed shape and is globally attracting with a scaling that is independent
of the diffusion function;
\item[(2)]
the halfline $[0,\infty)$ (Dawson and Greven \cite{DG96}, Baillon, Cl\'ement, Greven
and den Hollander \cite{BCGH97}), where the Feller branching diffusion is the unique
fixed point and is globally attracting with a scaling that depends on the asymptotic
behavior of the diffusion function at infinity.
\end{itemize}
For higher-dimensional diffusions, the {\em analytic part} has been carried out for:
\begin{itemize}
\item[(3)]
isotropic diffusions taking values in a compact convex subset of $\R^d$
(den Hollander and Swart \cite{HS98}, Swart~\cite{S04}), where the
diffusion function with constant curvature is the unique fixed shape and is globally
attracting with a scaling that is independent of the diffusion function;
\item[(4)]
a class of probability-measure-valued diffusions (Dawson, Greven and Vaillancourt
\cite{DGV95}, Dawson and March \cite{DM95}), where the Fleming-Viot process is the
unique fixed shape and is globally attracting with a scaling that is independent of
the diffusion function;
\item[(5)]
a class of catalytic Wright-Fisher diffusions taking values in
$[0,1]^2$ (Fleischmann and Swart \cite{FS06}), where the diffusion
function for the first component is an autonomous Wright-Fisher
diffusion and the diffusion function for the second component is an
autonomous Wright-Fisher diffusion function multiplied by a {\it
catalyzing function} depending only on the first component. The
renormalization transformation effectively acts on the catalyzing
function. There are four attracting shapes for the catalyzing
function, depending on whether the initial catalyzing function is zero
or strictly positive at the boundary points of $[0,1]$, and these
attracting shapes are globally attracting with a scaling that is
independent of the catalyzing function.
\end{itemize}
The {\em stochastic part} for higher-dimensional diffusions has only
been completed for interacting Fleming-Viot processes (Dawson, Greven
and Vaillancourt \cite{DGV95}) and for mutually catalytic branching
diffusions taking values in $[0,\infty)^2$ (Cox, Dawson and Greven
\cite{CDG04}).

All previous studies deal with diffusions that have certain
simplifying properties. In the one-dimensional cases (1) and (2), as
well as in the two-dimensional case (5), the equilibrium of (\ref{sde0}) is reversible. As
a result, many explicit calculations can be performed that are crucial
for the analysis. For certain diffusions with compact state space,
which includes the cases (1), (3) and (4), there is a common
underlying structure (called ``invariant harmonics'', see Swart~\cite{S00}) that
allows the determination of the unique fixed shape and its domain of
attraction. In all cases where the state space is compact, the scaling
needed for convergence to an attracting shape depends only on
$(c_n)_{n\in\N_0}$, not on the diffusion function $g$. This is different in
case (2), where the state space is not compact. In all cases
except case (5), the fixed shapes turn out to be joint fixed
shapes for all $c>0$.

The goal of the present paper is to carry out the {\em analytic part}
of the renormalization program for a \emph{general} class of branching
diffusions taking values in $[0,\infty)^2$. The multi-dimensionality
and the non-compactness of the state space pose significant
challenges. Due to the multidimensionality, the well-definedness of
the renormalization transformation is nontrivial. The structure of the
fixed points/shapes turns out to be rather rich. In fact, we will
prove that, under certain restrictions, the class of fixed points is a
{\em 4-parameter family of diffusions with independent branching,
catalytic branching and mutually catalytic branching as the extremal
fixed points}, and they are joint fixed points of $F_c$ for all
$c>0$. Moreover, we will prove that all diffusion functions that are
comparable to these fixed points in an appropriate sense fall in their
domains of attraction.

\subsection{Outline}

The rest of the paper is organized as follows. In Section \ref{S2} we
formulate our main results, which come with varying degrees of
restrictions on the diffusion functions. Section \ref{S3} contains the
proof of the ergodicity of the SDE (\ref{sde0}), and basic properties
of the renormalization transformation. Section \ref{S4} proves
the identification of fixed points/shapes. Sections \ref{S5}
and \ref{S6} identify the domains of attraction for the
fixed points. In Appendices \ref{S3.1} and \ref{A2} we collect some
technical results needed for the proofs.


\section{Main results}
\label{S2}

In Section \ref{S2.1}, we formulate a key class of diffusion functions
$\cC$, for which the SDE (\ref{sde0}) has a unique weak
solution. Section \ref{S2.2} contains a theorem on the ergodicity of
the SDE (\ref{sde0}), defines the renormalization transformation,
formulates a subclass $\cH_{0^+} \subset \cC$ on which the
renormalization transformation is well-defined and, subject to a
conjecture on the preservation of certain boundary properties, can be
iterated. Section \ref{S2.3} gives the definition of certain
generalized fixed points/shapes, and identifies some special fixed
points/shapes. Section \ref{S2.4} contains results on the
identification of fixed points/shapes in $\cH_{0^+}$ under additional
regularity assumptions. Section \ref{S2.5} contains our main result on
the domains of attraction to the fixed points under further
assumptions. Lastly, Section \ref{S2.6} provides a brief discussion of
these results and lists some future challenges.

\subsection{Key class and uniqueness for the autonomous SDE}
\label{S2.1}

The renormalization transformation $F_c$ is based on (\ref{sde0}), which is the SDE for
the vector $\vec X(t) = (X_1(t),X_2(t)) \in [0,\infty)^2$ written out as
\be{SDEaut}
\begin{aligned}
dX_1(t) &= c\,[\theta_1 - X_1(t)]\,dt +
\sqrt{2g_1(X_1(t),X_2(t))}\, dB_1(t),\\
dX_2(t) &= c\,[\theta_2 - X_2(t)]\,dt +
\sqrt{2g_2(X_1(t),X_2(t))}\, dB_2(t),
\end{aligned}
\ee
where $c>0$, $\vec\theta = (\theta_1, \theta_2) \in [0,\infty)^2$, and $\vec B(t)
= (B_1(t), B_2(t))$ are independent standard Brownian motions on $\R^2$. The
corresponding generator is
\be{generator}
(L_{\vec\theta}^{c,g} f)(\vec x) = c \sum_{i=1}^2 (\theta_i-x_i)
\frac{\partial}{\partial x_i} f(\vec x) + \sum_{i=1}^2 g_i(\vec x)
\frac{\partial^2}{\partial x_i^2} f(\vec x)\ , \qquad f\in C_c^2([0,\infty)^2).
\ee
Note that, due to the absence of mixed partial derivatives, $L_{\vec \theta}^{c,g}$ can be
interpreted as the generator of a two-type \emph{branching diffusion} with state-dependent branching
rates $g_i(\vec x)/x_i$ $(i=1,2)$.

Abbreviate
\be{Adefs}
A_1 = [0,\infty) \times \{0\}, \quad A_2 = \{0\} \times [0,\infty).
\ee
We will say that a function $f\colon\,[0,\infty)^2\to[0,\infty)$ has \emph{boundary property}
\be{bdprops}
\begin{array}{llll}
&(\partial_1)   &\mbox{if} &\lim\limits_{\vec x\to\vec y} \frac{f(\vec x)}{x_1} = \gamma(\vec y)
\,\,\,\forall\,\vec y\in A_1 \cup A_2 \mbox{ with } \gamma \mbox{ continuous and } >0
\mbox{ on } A_1 \cup A_2,\\[0.4cm]
&(\partial_2)   &\mbox{if} &\lim\limits_{\vec x\to\vec y} \frac{f(\vec x)}{x_2} = \gamma(\vec y)
\,\,\,\forall\,\vec y\in A_1 \cup A_2 \mbox{ with } \gamma \mbox{ continuous and } >0
\mbox{ on } A_1 \cup A_2,\\[0.4cm]
&(\partial_{12})   &\mbox{if} &\lim\limits_{\vec x\to\vec y} \frac{f(\vec x)}{x_1x_2} = \gamma(\vec y)
\,\,\,\forall\,\vec y\in A_1 \cup A_2 \mbox{ with } \gamma \mbox{ continuous and } >0
\mbox{ on } A_1 \cup A_2.
\end{array}
\ee
Throughout the paper, the pair $g=(g_1,g_2)$ will be assumed to be in the following class.

\bd{def:H}{\bf[Class $\cC$]}\\
Let $\cC$ be the class of functions $g(\vec x) = (g_1(\vec x),g_2(\vec x))$ satisfying:\\
$(i)$ For $i=1,2$, $g_i$ is continuous on $[0,\infty)^2$ and $>0$ on $(0,\infty)^2$.\\
$(ii)$ For $i=1,2$, $g_i$ satisfies boundary property $(\partial_i)$ or $(\partial_{12})$.
\ed
Note that for $(g_1, g_2)\in \cC$ we can write $g_i(\vec x) = x_i \gamma_i(\vec x)$ or
$g_i(\vec x) = x_1 x_2 \gamma_i(\vec x)$ for some positive continuous function $\gamma_i$
on $[0,\infty)^2$, depending on whether $g_i$ satisfies boundary property $(\partial_i)$ or
$(\partial_{12})$. Note also that $g_1$ and $g_2$ vanish on $A_2$, respectively, $A_1$,  which is necessary to
guarantee that the diffusion stays within $[0,\infty)^2$. Thus, if we denote the
\emph{effective boundary} of $g$ by
\be{partgdef}
\partial g = \{\vec x\in[0,\infty)^2\colon\, g_1(\vec x) = g_2(\vec x) = 0\},
\ee
then $\partial g$ can be either of the following:
\be{bdpos}
A_1 \cap A_2, \quad A_1, \quad A_2, \quad A_1 \cup A_2.
\ee
These boundary constraints allow for the system (\ref{SDEaut}) to be treated as a
perturbation of either of the following diffusions:
\begin{itemize}
\item[(1)] {\it Independent branching}:
$(g_1, g_2)= (b_1 x_1, b_2 x_2)$, $b_1, b_2 >0, \partial g = A_1 \cap A_2$.
\item[(2)] {\it Catalytic branching}:
either $(g_1, g_2)= (b_1 x_1, c_2 x_1x_2)$, $b_1, c_2>0,
\partial g = A_2$;
or $(g_1, g_2) = (c_1 x_1x_2, b_2 x_2)$, $c_1, b_2 >0,
\partial g = A_1$.
\item[(3)] {\it Mutually catalytic branching}:
$(g_1, g_2) = (c_1 x_1 x_2, c_2 x_1 x_2)$, $c_1, c_2>0,
\partial g = A_1 \cup A_2$.
\end{itemize}
Such a perturbation is behind the following result of Athreya, Barlow, Bass and Perkins~\cite{ABBP02},
and Bass and Perkins \cite{BPprep}, which provides the starting point of our analysis. The latter paper
developed out of Dawson and Perkins \cite{DP05}, where H\"older continuity is assumed
rather than continuity, but the result there is not restricted to two dimensions as in~\cite{BPprep}.

\bt{thm:martproblem} {\bf [Well-posedness of martingale problem]}
{\rm (\cite{ABBP02}, \cite{BPprep})}\\
For all $c>0$, $g\in\cC$, $\vec\theta\in [0,\infty)^2$ and $\vec x\in [0,\infty)^2$, with the possible exception
of the case when $\vec x=(0,0)$, $\vec\theta\in (0,\infty)^2$, and either $g_1$ or $g_2$ satisfies boundary property
$(\partial_{12})$, the martingale problem associated with the generator in {\rm (\ref{generator})} has a unique solution with starting position $\vec x$.
\et

As a consequence of Theorem \ref{thm:martproblem}, the SDE {\rm (\ref{SDEaut})} has a unique weak solution
for all $\vec\theta \in [0,\infty)^2$ and $\vec x\in [0,\infty)^2$,
with the possible exception of the case when $\vec x=(0,0)$, $\vec\theta\in(0,\infty)^2$, and either $g_1$ or $g_2$
satisfies boundary property $(\partial_{12})$. For each fixed $\vec\theta\in [0,\infty)^2$, the SDE (\ref{SDEaut}) defines
a Feller process satisfying the strong Markov property (see e.g.\ Theorem 4.4.2 in Ethier and Kurtz \cite{EK86} and
Corollary 11.1.5 in Stroock and Varadhan \cite{SV79}).

\medskip\noindent
{\bf Remark 1:} When $\vec\theta \in (0,\infty)^2$, $g\in\cal C$, $g_1$ and $g_2$
satisfy $(\partial_1)$, resp.\ $(\partial_2)$, the well-posedness of the martingale problem
was established in Athreya, Barlow, Bass and Perkins~\cite{ABBP02} for all initial conditions $\vec x\in [0,\infty)^2$.
When $\vec\theta\in (0,\infty)^2$, $g\in \cal C$, and $g_1, g_2$ both satisfy $(\partial_{12})$,
the well-posedness is established in Bass and Perkins~\cite{BPprep} for all intial condition $\vec x\in [0,\infty)^2\backslash\{(0,0)\}$.
Both \cite{ABBP02} and \cite{BPprep} use local perturbation arguments and the results are not restricted to
linear drift as considered here. Since the perturbation arguments are local, this implies that well-posedness
also holds for mixed boundaries, i.e.,\ $g_1$ satisfies $(\partial_1)$ and $g_2$ satisfies
$(\partial_{12})$, or vice versa. When either $g_1$ or $g_2$ satisfies $(\partial_{12})$, Lemma~35 of
Dawson and Perkins~\cite{DP05} shows that, for all $\vec x\in [0,\infty)^2\backslash{(0,0)}$, with probability 1 the unique
weak solution of (\ref{SDEaut}) with initial condition $\vec x$ never hits $(0,0)$, and hence we can restrict the state
space to $[0,\infty)^2\backslash \{(0,0)\}$. When $\vec\theta \in \partial [0,\infty)^2$, the local analysis
of \cite{ABBP02} and \cite{BPprep} still applies until the diffusion first hits the absorbing boundary, at which time
the diffusion becomes one-dimensional, a situation for which the well-posedness of the martingale problem is standard.

\medskip\noindent
{\bf Remark 2:} The proof given in \cite{ABBP02} requires the drift to be strictly positive in each
component on $\partial[0,\infty)^2$. However, as pointed out in Bass and Perkins \cite{BP04}, it is
sufficient that the inward normal component of the drift is strictly positive on $\partial [0,\infty)^2$,
which holds in our setting when $\vec\theta\in (0,\infty)^2$.

\medskip\noindent
{\bf Remark 3:} It would be considerably more difficult to deduce from Theorem \ref{thm:martproblem}
the well-posedness of the martingale problem for the system (\ref{SDEcoupled}), for
which one would need to restrict the state space. To deduce the Feller property,
one would need to restrict the state space even further and impose growth
conditions on the diffusion function $g$, typically $g_1(\vec x)+ g_2(\vec x)
= O(x_1^2+x_2^2)$ (see, e.g., Shiga and Shimizu~\cite{SS80}, Cox, Dawson and
Greven~\cite{CDG04}). We will not resolve these issues here, since they
belong to the {\it stochastic part} of the renormalization program, which remains
open.

\subsection{Equilibrium distribution and renormalization transformation}
\label{S2.2}

Our first result shows that (\ref{SDEaut}) has a unique equilibrium for
the class $\cC$. The proof will be given in Section \ref{S3.2}. Henceforth
$\cL$ denotes law.

\bt{thm:SDEeq} {\bf [Equilibrium distribution]}\\
For all $g\in\cC$, $\vec\theta\in [0,\infty)^2$ and $c>0$, {\rm (\ref{SDEaut})}
has a unique equilibrium distribution $\Gamma^{c,g}_{\vec\theta}$, which is continuous
in $\vec\theta$ with respect to weak convergence of probability measures, and
\be{eqconv}
\cL(\vec X(t))_{\displaystyle \quad \Longrightarrow \quad \atop t\to\infty}
\Gamma^{c,g}_{\vec\theta} \qquad \forall\,\vec X(0)\in [0,\infty)^2.
\ee
\et
The convergence in (\ref{eqconv}) is crucial for the {\em stochastic part} of the renormalization
program (not considered here), while the uniqueness of the equilibrium is crucial for the definition
of the {\em renormalization transformation}, which we now define.

\bd{def:transform} {\bf [Renormalization transformation]}\\
The renormalization transformation $F_c$, acting on $g\in\cC$, is defined as
\be{transform}
(F_c g)_i(\vec\theta) = \int_{[0,\infty)^2} g_i(\vec x)\,
\Gamma^{c,g}_{\vec\theta}(d\vec x),   \qquad \vec\theta\in [0,\infty)^2,\,c>0,\,i=1,2.
\ee
\ed

\noindent
Henceforth we will denote expectation with respect to $\Gamma^{c,g}_{\vec\theta}$ by
$\E^{c,g}_{\vec\theta}$.

Without restrictions on the growth of $g$ at infinity, it is possible that
$F_cg$ is infinite. We therefore need to consider a tempered subclass of $\cal C$.

\bd{def:Ha} {\bf [Class $\cH_{0^+}$]}\\
$(i)$ For $a \geq 0$, let $\cH_a\subset\cC$ be the class of all $g\in\cC$
satisfying
\be{Ha}
g_1(x_1, x_2) + g_2(x_1, x_2) \leq C(1+x_1)(1+x_2) + a (x_1^2 + x_2^2),
\qquad (x_1,x_2)\in [0,\infty)^2,
\ee
for some $0<C=C(g)<\infty$. \\
$(ii)$ Let
\be{H0}
\cH_{0^+} = \bigcap_{a>0} \cH_a.
\ee
\ed

\noindent
Note that $\cH_{0^+}$ is much larger than $\cH_0$. In particular, $\cH_{0^+}$
includes diffusion functions that along the axes grow faster than linear but slower
than quadratic.

Our second result shows that $F_c$ is well-defined on the class $\cH_a$ when $0\leq a<c$,
preserves the effective boundary, and preserves the growth bound in (\ref{Ha}) though
with a different coefficient. The proof will be given in Section \ref{S3.3}.

\bt{thm:FCdomain} {\bf [Finiteness, continuity, preservation of $\partial g$ and growth bound]}\\
For $c>0$ and $0\leq a<c$, if $g\in\cH_a$, then $F_cg$ is finite and continuous on
$[0,\infty)^2$, $\partial F_cg =\partial g$, and $F_cg$ satisfies {\rm (\ref{Ha})}
with $a$ replaced by $\frac{c}{c-a}a$.
\et

To proceed with our analysis, we need the following:
\bcon{conj:FCbdprop} {\bf [Preservation of boundary properties]}\\
Let $g\in \cH_{0^+}$. \\
$(i)$ For $i=1,2$, if $g_i$ satisfies $(\partial_i)$, then so does $(F_cg)_i$ for all $c>0$.\\
$(ii)$ For $i=1,2$, if $g_i$ satisfies $(\partial_{12})$, then so does $(F_cg)_i$ for all $c>0$.
\econ

\noindent
In Section \ref{S3.4} we will explain why this conjecture is plausible. Combining
Theorem \ref{thm:FCdomain} with Conjecture \ref{conj:FCbdprop}, we get:

\bc{cor:H0pres} {\bf [Preservation of class $\cH_{0^+}$]}\\
For all $c>0$, the class $\cH_{0^+}$ is preserved under $F_c$, i.e., $F_c g\in\cH_{0^+}$
for all $g\in\cH_{0^+}$.
\ec

\noindent
The latter is a key property, because it allows us to \emph{iterate}
$F_c$ on $\cH_{0+}$ and investigate the orbit $F^{[n]}g =
F_{c_{n-1}}\circ\cdots\circ F_{c_0}g$, $n\in\N_0$. We will not need
Conjecture~\ref{conj:FCbdprop} or Corollary~\ref{cor:H0pres} until we
study the iterates $F^{[n]}$ in Section \ref{S2.5}.

The subquadratic growth bound imposed by $\cH_{0^+}$ cannot be relaxed: we will see in
Corollary \ref{cor:iterfixshape} below that $F_c$ cannot be iterated indefinitely on
$\cH_a$ for any $a>0$.

\subsection{Definition and examples of fixed points and fixed shapes}
\label{S2.3}

We next give the definition of fixed points and fixed shapes of $F_c$. Generalizing our definition given in the
introduction, we allow for the case where $F_cg=\la g$ with $\la$ not a constant but a diagonal matrix. These {\em
generalized fixed shapes} do not give rise to universality classes as defined in Section~\ref{S1.2}, but they may be
relevant for studying finer properties of the orbit $(F^{[n]}g)_{n\in\N_0}$.

\bd{def:fixedshapes} {\bf [Generalized fixed shapes and points]}\\
The pair $g=(g_1, g_2) \in \cH_a$ with $a\in [0,c)$ is called a generalized fixed shape of $F_c$ if
\be{fixsh}
F_c (g_1, g_2) = (\lambda_1 g_1, \lambda_2 g_2) \quad \mbox{ for some }
\quad \lambda_1,\lambda_2>0.
\ee
If $\lambda_1=\lambda_2$, then $g$ is called a fixed shape, and if
$\lambda_1=\lambda_2=1$, then $g$ is called a fixed point of $F_c$.
\ed

Our third result identifies a family of fixed points and (generalized) fixed shapes of $F_c$. The proof
is nontrivial because of integrability issues, and will be given in Section \ref{S3.3}.

\bt{thm:fixshap} {\bf [Examples of fixed points and fixed shapes]} \\
{\rm (i)} The pair
\be{fp}
(g_1,g_2) = (b_1x_1 + c_1x_1x_2,\,b_2x_2 + c_2x_1x_2)
\ee
is a fixed point of $F_c$ in $\cH_{0^+}$ for all $c>0$ and all
$b_1,b_2,c_1,c_2\geq 0$ with $(b_1+c_1)(b_2+c_2)>0$.\\
{\rm (ii)} The pair
\be{fs}
(g_1,g_2) = (a_1x_1^2 + b_1x_1 + c_1x_1x_2,\,a_2x_2^2 + b_2x_2 + c_2x_1x_2)
\ee
is a generalized fixed shape of $F_c$ in $\cH_{a_1 \vee a_2}$ for all $c>0$, $0<a_1,a_2<c$ and
$b_1,b_2,c_1,c_2\geq 0$. The corresponding scaling constants are
\be{fseig}
\lambda_1 = \frac{c}{c-a_1},\,\lambda_2=\frac{c}{c-a_2}.
\ee
\et

\noindent
Diffusion functions of the form in (\ref{fp}) are \emph{mixtures} of independent branching,
catalytic branching and mutually catalytic branching (recall Section \ref{S2.1}), all of
which are in the class $\cH_{0+}$. We will see in Theorem \ref{thm:FPS} below that, under
additional regularity conditions, such mixtures are the {\it only} fixed points of
$F_c$. Diffusion functions of the form in (\ref{fs}) are mixtures of these fixed points and
the Anderson branching diffusion $(g_1, g_2) = (a_1 x_1^2, a_2 x_2^2)$. The latter do
not fall in the class $\cH_{0+}$.

The following corollary of Theorem \ref{thm:fixshap}
shows that $F_cg$ cannot be defined for all $g\in \cH_a$ with $a\geq
c$, and $F_c$ cannot be iterated indefinitely on $\cH_a$ for any
$a>0$. The proof will be given in Section \ref{S3.3}.

\bc{cor:iterfixshape}{\bf [Divergence of iterated fixed shapes]} \\
Let $g_i(\vec x) = \alpha_i x_i^2 +\beta_i x_i +\gamma_i x_1 x_2$ with $\alpha_i>0$ and $\beta_i, \gamma_i\geq 0$,
$i=1,2$. Let $(c_n)_{n\in\N_0}$ be the positive sequence that defines $F^{[n]}$, $($see~$(\ref{Fkdef}))$. Let $n_0 =\min\{n\in\N :
(\alpha_1\vee \alpha_2)\sum_{i=0}^{n-1} c_i^{-1}  \geq 1\}$. Then
\be{iterfs}
\Big((F^{[n]}g)_1, (F^{[n]}g)_2\Big)
= \Big(\frac{1}{1-\alpha_1 \sum_{i=0}^{n-1} c_i^{-1}}\, g_1,\ \frac{1}{1-\alpha_2\sum_{i=0}^{n-1}c_i^{-1}}\, g_2\Big),
\qquad 0\leq n < n_0,
\ee
while $(F^{[n_0]}g)_1+(F^{[n_0]}g)_2 \equiv \infty$ on $(0,\infty)^2$.
\ec

\subsection{Identification of fixed points and fixed shapes}
\label{S2.4}

Our fourth result rules out generalized fixed shapes in $\cH_{0^+}$ with an upgoing component. The proof will be given
in Section \ref{S4.3}.

\bt{thm:upfs}{\bf [No fixed shapes in $\cH_{0^+}$ with an upgoing component]} \\
For $c>0$, there is no $g\in \cH_{0^+}$ such that either $(F_cg)_1= \lambda_1 g_1$
with $\lambda_1>1$ or $(F_cg)_2= \lambda_2 g_2$ with $\lambda_2>1$.
\et

Our fifth result does the same for generalized fixed shapes with a downgoing component, but only
under mild additional regularity conditions. The proof will be given in Section \ref{S4.3}. Below, in line with general
topological notation, $\liminf_{\vec x\to(\infty,\infty)}$ denotes the infimum of all limits along sequences
tending to $(\infty,\infty)$.

\bt{thm:downfs}{\bf [Sufficient conditions for no downgoing fixed shapes in $\cH_{0^+}$]} \\
Let $c>0$.\\
$(i)$
There is no $g\in\cH_{0+}$ such that $F_c(g_1, g_2) = (\lambda_1 g_1, \lambda_2 g_2)$ with
$0< \lambda_1, \lambda_2 <1$ and
\beq
\label{fscond4}
\liminf_{\vec x\to (\infty, \infty)}
\left(\frac{g_1(\vec x)}{x_1^2} + \frac{g_2(\vec x)}{x_2^2}\right) =0.
\eeq

\noindent
$(ii)$ There is no $g\in\cH_{0+}$ such that $(F_cg)_1 = \lambda_1 g_1$ for some $0<\lambda_1<1$
and $g$ satisfies any of the following conditions:
\be{fscond1}
\!\!\!\!\!\!\!\!\!\!\!\!\!\!\!\!\!\!\!\!\!\!\!\!\!\!\!\!\!\!\!\!\!\!\!\!\!\!\!\!\!\!\!\!\!\!\!\!\!\!\!\!\!\!\!\!\!\!\!\!
\!\!\!\!\!\!\!\!\!\!\!\!\!\!\!\!\!\!\!
\bullet \qquad \qquad \qquad \quad \qquad \qquad \qquad g_1>0 \mbox{ on } A_1\setminus\{(0,0)\},
\ee
\be{fscond2}
\!\!\!\!\!\!\!\!\!\!\!\!\!\!\!\!\!\!\!\!\!\!\!\!\!\!\!\!\!\!\!\!\!\!\!\!\!\!\!\!\!\!\!\!\!\!\!\!\!\!\!\!\!\!\!\!\!\!\!\!
\!\!\!\!\!\!\!\!\!\!\!\!\!\!\!\!\!\!\!\!\!\!\!
\bullet \qquad\qquad\qquad \qquad \qquad \quad \qquad\ \  \liminf_{\vec x\to (\infty, \infty)}\ \frac{g_1(\vec x)}{x_1 x_2} >0.
\ee
A similar result holds with the indices $1$ and $2$ interchanged.
\et
{\bf Remark:} Conditions (\ref{fscond4}) and (\ref{fscond2}) are complementary. Note that one particular case not
covered by conditions (\ref{fscond4})--(\ref{fscond2}) is when $g_1$ vanishes on both axes, $g_1(\vec x)=o(x_1x_2)$
as $\vec x\to (\infty,\infty)$, and $g_2(\vec x) = x_1 x_2$. In that case we cannot rule out the possibility of $g_1$ being a downgoing
fixed shape.
\medskip

In Theorem \ref{thm:fixshap} we identified a 4-parameter family of fixed points. To show that
these are the only fixed points, we need to impose strong additional regularity conditions.

Abbreviate
\be{Rdefs}
R_\infty = \{(0,\infty),(\infty,0), (\infty,\infty)\}
\ee
and
\be{hdefs}
h_{(\infty,0)}(\vec x) = x_1, \qquad h_{(0,\infty)}(\vec x) = x_2,
\qquad h_{(\infty,\infty)}(\vec x) = x_1x_2.
\ee

\bd{def:H0r} {\bf [Class $\cH_0^r$]}\\
Let $\cH_0^r$ be the set of $g\in\cH_0$ satisfying
\begin{eqnarray}
&(i)& \inf_{\vec x\in [s,\infty)^2}\ g_i(\vec x) > 0
\quad \forall\,s>0,\,i=1,2,
\label{gcon1} \\
&(ii)& \lim_{\vec x\to \vec z} \frac{g_i(\vec x)}{h_{\vec z}(\vec x)}
= \lambda_{i,\vec z} \in [0,\infty)
\quad \forall\,\vec z \in R_\infty,\,i=1,2.
\label{gcon3}
\end{eqnarray}
\ed

\noindent
Note that $\cH^r_0 \subset \cH_0 \subset \cH_{0^+}$.
Also note that, because $g_1$ vanishes on $A_2$ and $g_2$ on $A_1$, necessarily
$\lambda_{1,(0,\infty)}=\lambda_{2,(\infty,0)}=0$.

Our sixth result is the following. The proof will be given in Section \ref{S4.1}.

\bt{thm:FPS} {\bf [Identification of fixed points in $\cH_0^r$]}\\
Let $c>0$ and $g=(g_1, g_2) \in \cH_0^r$. If $F_c(g_1, g_2) = (g_1, g_2)$, then
\be{mixfix}
\begin{aligned}
g_1(\vec x) = \lambda_{1, (\infty,0)} x_1 + \lambda_{1, (\infty,\infty)} x_1 x_2, \\
g_2(\vec x) = \lambda_{2, (0,\infty)} x_2 + \lambda_{2, (\infty,\infty)} x_1 x_2,
\end{aligned}
\ee
where $\lambda_{i, \vec z}$, $\vec z\in R_\infty$, are defined in $(\ref{gcon3})$.
\et

\subsection{Domain of attraction of fixed points}
\label{S2.5}

Our seventh and final result is on the domain of attraction of the iterated maps $F^{[n]} = F_{c_{n-1}}\circ\cdots\circ
F_{c_0}$, $n\in\N_0$, for a fixed positive sequence $(c_n)_{n\in\N_0}$. We show that, provided $\inf_{n\in\N_0} c_n >0$ and $\sum_{n\in\N_0} c_n^{-1} =\infty$, all diffusion functions that are comparable
to a mixture of the fixed points fall into its domain of attraction. In Section \ref{S5}, we will give the proof
for the special case $c_n\equiv c$, while in Section \ref{S6}, we prove the result for varying $c_n$.

\bt{thm:DA} {\bf [Domain of attraction of fixed points]}\\
Let $(c_n)_{n\in \N_0}$ be a sequence such that $\inf_{n\in \N_0} c_n >0$ and $\sum_{n\in\N_0} c_n^{-1} = \infty$. Let $g\in\cH_0^r$ be such that
\begin{equation}
\label{DAgcon1}
g_i(\vec x) \geq \alpha_i x_i + \beta_i x_1 x_2, \qquad \alpha_i, \beta_i\geq 0,\ \alpha_i+\beta_i>0,\ i=1,2.
\end{equation}
Then
\begin{equation}
\label{DAFnlim}
\lim_{n\to\infty} (F^{[n]}g)_i(\vec\theta) = \sum_{\vec z\in R_\infty}
\lambda_{i,\vec z}\,h_{\vec z}(\vec\theta) \qquad \forall\,\vec\theta\in [0,\infty)^2,\,i=1,2,
\end{equation}
where $h_{\vec z}$, $\lambda_{i, \vec z}$, $\vec z\in R_\infty$, are defined in
$(\ref{hdefs})$ and $(\ref{gcon3})$.
\et

\noindent
What this says is that under the iterates $F^{[n]}$, \emph{any} $g$ that is properly \emph{minorized}
and has the \emph{same behavior at infinity} as a mixture of the fixed points, converges to that mixture pointwise as $n\to\infty$.
\medskip

\noindent
{\bf Remark 1:} Note that Theorem \ref{thm:DA} implicitly assumes Conjecture \ref{conj:FCbdprop}. To be formally correct, in Theorem \ref{thm:DA} we should replace $\cH_0^r$ by the largest subclass of $\cH_0^r$ that is preserved by $F_c$ for all $c>0$.
\medskip

\noindent
{\bf Remark 2:} The condition $\inf_{n\in\N_0}c_n>0$ means that we
partially exclude the regime of {\em large clusters} (see e.g.\ Dawson
and Greven~\cite{DG93b}). We do not believe this assumption to be essential. As long as $\sum_{n\in\N_0} c_n^{-1}
=\infty$, i.e., the associated random walk on $\Omega_N$ with transition rate kernel $a_N(\cdot, \cdot)$ is recurrent, we expect there to be universality and the convergence in (\ref{DAFnlim}) to hold.

\subsection{Discussion and future challenges}
\label{S2.6}

The results in Sections \ref{S2.2}--\ref{S2.5} constitute a \emph{partial completion} of
the analytic part of the renormalization program outlined in Section \ref{S1.2}.
We have formulated $\cH_{0+}$ as the class on which the renormalization transformation is properly defined and, apart from Conjecture \ref{conj:FCbdprop}, can be iterated. We have proved absence of upgoing fixed shapes in
this class, and absence of downgoing fixed shapes under mild regularity conditions, given
by (\ref{fscond4})--(\ref{fscond2}). Furthermore, we have identified our 4-parameter family of fixed points in (\ref{fp}) as the only fixed points in a subclass $\cH_0^r$ of the smaller class $\cH_0$, given by the strong regularity conditions (\ref{gcon1})--(\ref{gcon3}). Finally, we have found the domain of attraction of these
fixed points in $\cH_0^r$ supplemented with the lower bound (\ref{DAgcon1}), i.e., diffusion
functions that are comparable to a mixture of the fixed shapes. There are several open problems remaining, the chief
among which are:

\begin{itemize}
\item[$(1)$] Verify Conjecture \ref{conj:FCbdprop}, i.e., establish that the renormalization transformation can be iterated on $\cH_{0^+}$.

\item[$(2)$] Remove assumptions (\ref{fscond4})--(\ref{fscond2}) in the proof of the absence of downgoing fixed shapes in $\cH_{0^+}$.

\item[$(3)$] Show that the fixed points in $(\ref{fp})$ are the only fixed points in $\cH_{0^+}$. In particular, remove assumption (\ref{gcon3}) and the bound $g_1(\vec x)+g_2(\vec x) \leq C(1+x_1)(1+x_2)$ in $\cH_0^r\subset \cH_0$.

\item[$(4)$] Strengthen (2) and (3) by determining whether it is actually true that the fixed shapes in (\ref{fs}) are the only fixed shapes in $\cC$.

\item[$(5)$] Study the orbit of $(F^{[n]}g)_{n\in\N_0}$ when the behavior of $g$ at infinity is \emph{different} from that of the fixed points. In that case we still expect convergence, but only after $F^{[n]}g$ is
scaled with $n$ in some appropriate manner. For diffusions on the halfline $[0,\infty)$, this study
was successfully completed in Baillon, Cl\'ement, Greven and den Hollander \cite{BCGH97}, which
raises some hope that it can be carried through on the quadrant as well.

\end{itemize}

The questions we treated in this paper and the open problems we
just mentioned have close connections to probabilistic potential
theory of diffusions and Markov chains taking values in the
quadrant. Our proofs strongly lean on the observation that the fixed
points we build are mixtures of {\em extremal universal harmonic
functions} of the interaction chains described in Section~\ref{S1.2}.
The problem of finding all fixed points then requires identifying the
{\em universal Martin boundary} of these Markov chains. The reader
interested in this point of view can find the necessary concepts in
Pinsky~\cite{P95}. Harmonic functions have played an important role in
earlier studies of the analytic part of the renormalization
program. In particular, the convergence proofs in the cases (1), (3)
and (4) mentioned in Section~\ref{S1.3} all depend on a special
property of these models, called ``invariant harmonics'' (see
Swart~\cite{S00}). Case (2) uses moment equations combined with
comparison arguments, while case (5) uses a representation in terms of
a superprocess. Due to multi-dimensionality and non-compactness, these
tools either do not apply or are insufficient for our
model. However, our present methods have their limitations as
well. In particular, in their present state they can only be used to
prove convergence to {\em joint fixed points of $F_c$ for all $c>0$},
as opposed to fixed shapes, or cases where there might be different
fixed points of $F_c$ for different values of $c$. Moreover, we can treat only functions that are
perturbations of these fixed points, albeit in a rather large class.

Another interesting question is to study multi-type branching models
with more than two types. The class of random catalytic networks introduced in Dawson and
Perkins~\cite{DP05} and generalized in Kliem~\cite{K07} provide a rich class of fixed points of the
renormalization transformation. However our results here do not extend trivially to
higher dimensions, because we need the well-posedness of the martingale
problem (Theorem~\ref{thm:martproblem}), which is more delicate in higher
dimensions. Also, our proof of the formula (\ref{eq:mixmoment}) for the mixed moment
$X_1X_2$ does not extend to mixed moments of higher order.


\section{Proof of Theorems \ref{thm:SDEeq}, \ref{thm:FCdomain},
\ref{thm:fixshap} and Corollary \ref{cor:iterfixshape}} \label{S3}

In Section \ref{S3.2} we give the proof of Theorem \ref{thm:SDEeq}, in Section \ref{S3.3}
that of Theorems \ref{thm:FCdomain}, \ref{thm:fixshap} and Corollary \ref{cor:iterfixshape}.
In Section \ref{S3.4} we discuss Conjecture \ref{conj:FCbdprop}. Along the way we need a
proposition on moment equations for the equilibrium distribution $\Gamma_{\vec\theta}^{c,g}$,
which will turn out to be fundamental in our analysis. This proposition is formulated and
proved in Appendix \ref{S3.1}.

\subsection{Proof of Theorem \ref{thm:SDEeq}}
\label{S3.2}

We break down the proof of Theorem \ref{thm:SDEeq} into four parts: {\em existence},
{\em uniqueness}, {\em weak continuity} and {\em convergence}. For uniqueness and
convergence, we need to distinguish between $\vec\theta \in (0,\infty)^2$ and
$\vec\theta\in \partial [0,\infty)^2$.

\noindent
{\bf Existence:}

\bpr
If we denote the distribution of $\vec X(t)$ by $\mu_t$, with $\mu_0= \delta_{\vec x}$ for
some arbitrary $\vec x\in[0,\infty)^2$, then it suffices to show that $\{\nu_t\colon\,\nu_t
= \frac{1}{t} \int_0^t \mu_s ds \}_{t\geq 0}$ forms a tight family of distributions
on $[0,\infty)^2$. Indeed, we can then find a sequence $(t_n)$ tending to infinity
such that $\nu_{t_n}$ converges weakly to a limiting distribution $\nu$. Consequently,
for any $f\in C^2_c([0,\infty)^2)$,
\begin{eqnarray}
\int (L^{c,g}_{\vec\theta}f)(\vec x)\ \nu(d\vec x)  \nn
&=& \lim_{n\to\infty} \int (L^{c,g}_{\vec \theta}f)(\vec x)\ \nu_{t_n}(d\vec x) \\
&=& \lim_{n\to\infty} \frac{1}{t_n} \int_0^{t_n} \int (L_{\vec\theta}^{c,g} f)(\vec x) \mu_s(d\vec x)\ ds  \nn \\
&=& \lim_{n\to\infty} \frac{1}{t_n}
\E_{\mu_0}\Big[ \int_0^{t_n} (L_{\vec\theta}^{c,g} f)(\vec X(s)) ds \Big]  \nn \\
&=& \lim_{n\to\infty} \frac{1}{t_n} \E_{\mu_0}[f(\vec X(t_n)) - f(\vec X(0)] = 0,
\end{eqnarray}
where the first line uses that $\nu_{t_n}$ converges weakly to $\nu$, the
second line uses the definition of $\nu_{t_n}$, the third lines uses the definition of $\mu_s$ and Fubini,
and the fourth line uses that $f(\vec X(t))-f(\vec X(0)) - \int_0^t (L_{\vec\theta}^{c,g}f)( \vec X(s))ds$ is
a martingale and $f$ is bounded. Since $\int (L_{\vec\theta}^{c,g}f)(\vec x)\nu(d\vec x)=0$
for all $f\in C^2_c([0,\infty)^2)$, which form an algebra of functions that is dense
in the space of continuous functions on $[0,\infty)^2$ vanishing at $\infty$, it
follows from Theorem 4.9.17 in Ethier and Kurtz \cite{EK86} that $\nu$ is an
equilibrium distribution for (\ref{SDEaut}).

Tightness of the family $\{\nu_t\}_{t\geq 0}$ follows from the following lemma.

\bl{tightlem}{\bf[Tightness estimate]}\\
Let $(\vec X(t))_{t\geq 0}$ be the unique solution of the martingale problem for
$L_{\vec\theta}^{c,g}$ with initial condition $\vec X(0)=\vec x$. Then
\begin{equation} \label{tigh}
\E\left[X_i(t)-\theta_i\right] \leq (x_i-\theta_i)e^{-ct}, \qquad i=1,2,\, t\geq 0.
\end{equation}
\end{lemma}

\bpr
For any $\rho_1,\rho_2>0$, the function $f(t, \vec x)=\sum_{i=1}^2\rho_i
(x_i-\theta_i)e^{ct}$ satisfies
\be{Lparrel}
\Big(L^{c,g}_{\vec\theta}+\frac{\partial}{\partial t}\Big)f(t, \vec x)
= c {\sum_{i=1}^2}\rho_i(\theta_i-x_i)e^{ct}
+ c{\sum_{i=1}^2}\rho_i(x_i-\theta_i)e^{ct}=0,
\ee
and therefore the process $\sum_{i=1}^2\rho_i(X_i(t)-\theta_i)e^{ct}$ is a local
martingale. Introduce stopping times
\be{taun}
\tau_n =\inf\left\{t\geq 0:{\sum_{i=1}^2}\rho_iX_i(t)\geq n\right\},
\qquad n\in\N.
\ee
Then
\begin{eqnarray}
\sum_{i=1}^2 \rho_i(x_i-\theta_i)
&=& \sum_{i=1}^2\rho_i
\E\left[\Big(X_i(t\wedge\tau_n)-\theta_i\Big)e^{c(t\wedge\tau_n)}\right] \\
&=& \sum_{i=1}^2\rho_i \E\left[(X_i(t)-\theta_i)e^{ct}1_{\{\tau_n>t\}}\right]
+ \sum_{i=1}^2\rho_i \E\left[(X_i(\tau_n)-\theta_i)
e^{c\tau_n}1_{\{\tau_n\leq t\}}\right]. \nn
\end{eqnarray}
For $n\geq \sum_{i=1}^2\rho_i\theta_i$, the second term in the right-hand side
is nonnegative, so letting $n\to\infty$ we find that
\be{sumrel}
\sum_{i=1}^2 \rho_i \E\left[X_i(t)-\theta_i\right]\,
e^{ct}\leq \sum_{i=1}^2\rho_i (x_i-\theta_i).
\ee
Since $\rho_1,\rho_2>0$ are arbitrary, we arrive at (\ref{tigh}).
\epr

This completes the proof of the existence.
\epr

\noindent
{\bf Uniqueness:}

\bpr
We distinguish between $\vec\theta$ in the interior resp. on the boundary of $[0,\infty)^2$.

\medskip\noindent
\underline{$\vec\theta \in (0,\infty)^2$}: By Theorem \ref{thm:martproblem}, the
unique weak solution $(\vec X(t))_{t\geq 0}$ of (\ref{SDEaut}) is a strong Markov process. By
Remark 1 following Theorem~\ref{thm:martproblem}, we restrict the state space to be
$[0,\infty)^2\backslash\{(0,0)\}$ for the cases where weak uniqueness is not known when
$\vec X(0)=(0,0)$. If
$(\vec X(t))_{t\geq 0}$ has two distinct equilibrium distributions, then we can find two extremal
equilibrium distributions $\mu$ and $\nu$ that are singular with respect to each
other (see e.g.\ Theorem 6.9 in Varadhan \cite{V01}). This implies that there exist
$\vec x, \vec y\in[0,\infty)^2$ such that the transition kernels $p_t(\vec x, \cdot)$ and $p_t(\vec y, \cdot)$
are mutually singular for all $t>0$. However, if $\vec x, \vec y \in (0,\infty)^2$, then we can first apply Theorem
\ref{thm:supportthm} to transport the diffusions started at $\vec x$, resp.~$\vec y$, to a
common small neighborhood with positive probability, and subsequently apply Corollary \ref{cor:domaindiffdensity} to see that $p_t(\vec x, \cdot)$
and $p_t(\vec y, \cdot)$ cannot be singular for all $t>0$. On the other hand, when either $\vec x$ or $\vec y\in \partial [0,\infty)^2$, it suffices
to note that the drift in (\ref{SDEaut}) forces the diffusion to enter $(0,\infty)^2$ instantly, which we justify shortly. Then, again by
Theorem \ref{thm:supportthm}, the diffusion can be kept in $(0,\infty)^2$ up to any fixed time with positive probability, which
reduces it to the case $\vec x, \vec y \in (0,\infty)^2$.

We now show that, for $\vec X(0)=\vec x\in \partial [0,\infty)^2$, $(\vec X(t))_{t\geq 0}$ enters $(0,\infty)^2$ instantly. Consider first
the case $\vec X(0) \in \{0\}\times (0,\infty)$. Let $\vec X(0)=(0,y)$ with $y>0$, and let $\tau_\epsilon =
\inf\{t>0\colon\,|X_2(t)-X_2(0)|\geq y/2 \textrm{ or } X_1(t) \geq \epsilon\}$.
Then $X_1(t\wedge \tau_\epsilon) - \int_0^{t\wedge\tau_\epsilon}c(\theta_1-X_1(s))ds$
is a martingale, and
\begin{equation}\label{eq:leaveboundary1}
\E[X_1(t\wedge\tau_\epsilon)]
= \E \left[\int_0^{t\wedge\tau_\epsilon} c(\theta_1-X_1(s))ds\right].
\end{equation}
Letting $t\to\infty$, we find that for $\epsilon$ small,
\begin{equation}\label{eq:leaveboundary2}
\epsilon \geq \E[X_1(\tau_\epsilon)]
= \E \left[\int_0^{\tau_\epsilon} c(\theta_1-X_1(s))ds\right]
\geq \frac{c\,\theta_1}{2} \E[\tau_\epsilon].
\end{equation}
Therefore $\E[\tau_\epsilon] \to 0$ as $\epsilon \downarrow 0$, which is possible only
if $(\vec X(t))_{t\geq 0}$ enters $(0,\infty)^2$ instantly. The case $\vec X(0) \in (0,\infty)
\times \{0\}$ is analogous. For $\vec X(0)=(0,0)$, a similar argument shows
that $(\vec X(t))_{t\geq 0}$ enters $[0,\infty)^2 \backslash \{(0,0)\}$ instantly, which reduces it to
the previous cases.

\medskip\noindent
\underline{$\vec\theta \in \partial [0,\infty)^2$}: If $\theta_1=0$, then
$\E_{\vec\theta}^{c,g}[X_1] = \theta_1 =0$ for any equilibrium distribution
$\Gamma_{\vec\theta}^{c,g}$ by Proposition \ref{prop:momenteqn}. In particular,
$\Gamma_{\vec\theta}^{c,g}$ is concentrated on $\{0\}\times[0,\infty)$. Furthermore,
$(X_1(t))_{t\geq 0}$ is a local supermartingale, and hence $\{0\}\times [0,\infty)$ is an
absorbing set. The equilibria for $(\vec X(t))_{t\geq 0}$ are therefore exactly the equilibria
for $(\vec X(t))_{t\geq 0}$ restricted to the axis $\{0\}\times[0,\infty)$, which is a
one-dimensional diffusion. The proof of the existence and the uniqueness of the
equilibrium distribution for this one-dimensional diffusion can be deduced either
from explicit calculations as in Baillon, Cl\'ement, Greven and den Hollander~\cite{BCGH97}, or
from the same argument as above for the two-dimensional diffusion with $\vec\theta\in(0,\infty)^2$.
The situation is similar if $\theta_2=0$.
\epr

\noindent
{\bf Weak continuity:}

\bpr
We will show that $\Gamma_{\vec\theta}^{c,g}$ is weakly continuous in $\vec\theta$.
Let $(\vec\theta_n)$ be a sequence such that $\vec\theta_n\to \vec\theta$ in
$[0,\infty)^2$. It suffices to show that $\{\Gamma_{\vec\theta_n}^{c,g}\}_{n\in\N}$
is tight, and that any weak limit point of $\Gamma_{\vec\theta_n}^{c,g}$ is an
equilibrium distribution for the SDE (\ref{SDEaut}), which must be the unique
$\Gamma_{\vec\theta}^{c,g}$. Tightness of $\{\Gamma_{\vec\theta_n}^{c,g}\}_{n\in\N}$
follows from (\ref{eq:1stmoment}). Suppose that $\Gamma_{\vec\theta_n}^{c,g}$ converges
weakly to a distribution $\nu$.  Then for any $f\in C^2_c([0,\infty)^2)$,
\begin{eqnarray}
\int_{[0,\infty)^2} (L_{\vec\theta}^{c,g} f)(x) \nu(dx)
\!\!\!\!&=&\!\!\!\! \int_{[0,\infty)^2} (L_{\vec\theta_n}^{c,g} f)(x) \Gamma_{\vec\theta_n}^{c,g}(dx)
+ \int_{[0,\infty)^2} \Big[(L_{\vec\theta}^{c,g}- L_{\vec\theta_n}^{c,g}) f\Big](x)
\Gamma_{\vec\theta_n}^{c,g}(dx) \nonumber \\
&& \qquad + \int_{[0,\infty)^2}
(L_{\vec\theta}^{c,g} f)(x) \big[\nu(dx)-\Gamma_{\vec\theta_n}^{c,g}(dx)\big], \label{Lrels}
\end{eqnarray}
where the first term is zero because $\Gamma_{\vec\theta_n}^{c,g}$ is an equilibrium
distribution for the SDE in (\ref{SDEaut}) with parameter $\vec\theta_n$, the
second term tends to 0 as $n\to\infty$ because $f\in C^2_c([0,\infty)^2)$ and
$\| L_{\vec\theta}^{c,g} f(x) - L_{\vec\theta_n}^{c,g}f (x)\|_{\infty}\to 0$ as
$\theta_n\to\theta$, and the third term tends to 0 as $n\to\infty$ by the weak
convergence of $\Gamma_{\vec\theta_n}^{c,g}$ to $\nu$. Therefore $\int_{[0,\infty)^2}
(L_{\vec\theta}^{c,g} f)(x) \nu(dx) = 0 $ for all $f\in C^2_c([0,\infty)^2)$. By
Theorem 4.9.17 in Ethier and Kurtz \cite{EK86}, it follows that $\nu$ must be an
equilibrium distribution for (\ref{SDEaut}), and hence $\nu = \Gamma_{\vec\theta}^{c,g}$.
\epr

\noindent
{\bf Convergence:}

\bpr
We again distinguish between $\vec\theta$ in the interior resp.\ on the boundary of $[0,\infty)^2$.

\medskip\noindent
$\underline{\vec\theta \in (0,\infty)^2}$: Firstly, note that by Theorem \ref{thm:supportthm} and the fact
that $(\vec X(t))_{t\geq 0}$ started from $\partial [0,\infty)^2$ enters $(0,\infty)^2$ instantly (see the paragraph containing
(\ref{eq:leaveboundary1}--\ref{eq:leaveboundary2}) above), the equilibrium distribution $\Gamma_{\vec\theta}^{c,g}$
must assign positive measure to every open subset of $(0,\infty)^2$.

Secondly, we show that for almost all $\vec x \in [0,\infty)^2$ with respect to $\Gamma_{\vec\theta}^{c,g}$,
${\cal L}(\vec X(t) | \vec X(0)=\vec x)$ converges weakly to $\Gamma_{\vec\theta}^{c,g}$ as $t\to\infty$.
We achieve this by showing that, for almost all $(\vec x, \vec y) \in [0,\infty)^2\times [0,\infty)^2$ with respect to the
product measure $\Gamma_{\vec \theta}^{c,g}\times \Gamma_{\vec\theta}^{c,g}$, we can couple two solutions
$(\vec X(t))_{t\geq 0}$ and $(\vec Y(t))_{t\geq 0}$ of (\ref{SDEaut}) starting from $\vec x$, resp.~$\vec y$, such that $\lim_{t\to\infty} \P(\vec X(t) \neq \vec Y(t))=0$.
This goes as follows.

Let $\epsilon, \delta>0$ be chosen as in Corollary \ref{cor:domaindiffdensity}, where $b(\vec x) = c(\vec\theta -\vec x)$ and
$a(\vec x)=\left({g_1(\vec x)\atop 0} {0\atop g_2(\vec x)}\right)$ on $[0,\infty)^2$ (the definition of $(a,b)$ in the rest of
the plane $\R^2$ is irrelevant, for instance one may define it by reflection), $D=\{\vec x\in [0,\infty)^2 : \Vert \vec x-(1,1)\Vert< \frac{1}{2}\}$ and
$\vec x^* = (1,1)$. Note that $a(\cdot)$ is nondegenerate on $D$ for $g\in\cal C$.
If $(\vec X(t))_{t\geq 0}, (\vec Y(t))_{t\geq0}$ are two independent copies of the strong Markov process defined by (\ref{SDEaut}), then the joint process
$(\vec X(t), \vec Y(t))_{t\geq 0}$ is strong Markov and, by the same argument as for a single diffusion $(\vec X(t))_{t\geq 0}$, the joint process has a
unique equilibrium given by the product measure $\Gamma_{\vec \theta}^{c,g}\times \Gamma_{\vec\theta}^{c,g}$, which implies that the stationary process
$(\vec X(t), \vec Y(t))_{t\geq 0}$ with ${\cal L}(\vec X(0), \vec Y(0)) = \Gamma_{\vec \theta}^{c,g}\times \Gamma_{\vec\theta}^{c,g}$ is ergodic (see e.g.\
Theorem 6.9 in Varadhan~\cite{V01} and the remarks thereafter). Since $\Gamma_{\vec \theta}^{c,g}\times \Gamma_{\vec\theta}^{c,g}$ assigns positive measure to
$B_\epsilon(\vec x^*)\times B_\epsilon(\vec x^*)$, by the ergodic theorem almost surely $(\vec X(t), \vec Y(t))_{t\geq 0}$ visits the set
$B_\epsilon(\vec x^*)\times B_\epsilon(\vec x^*)$ after any finite time $T$. In particular, for almost all $(\vec x, \vec y)$ with respect to
$\Gamma_{\vec \theta}^{c,g}\times \Gamma_{\vec\theta}^{c,g}$,
almost surely the Markov process $(\vec X(t), \vec Y(t))_{t\geq 0}$ starting from $(\vec x, \vec y)$ visits $B_\epsilon(\vec x^*)\times B_\epsilon(\vec x^*)$ after
any finite time $T$. For such a pair $(\vec x, \vec y)$, we construct the coupled process as follows. Start the independent processes $(\vec X(t))_{t\geq0}$ and
$(\vec Y(t))_{t\geq0}$
with initial conditions $\vec x$, resp.~$\vec y$. Then $\tau = \inf\{t\geq 0 : (\vec X(t), \vec Y(t)) \in B_\epsilon(\vec x^*)\times B_\epsilon(\vec x^*)\} <\infty$
almost surely. By Corollary \ref{cor:domaindiffdensity}, the conditional transition probability kernels
$\mu_{\vec X}=\P(\vec X(\tau +\delta) \in \cdot\, |(\vec X(\tau), \vec Y(\tau)))$ and $\mu_{\vec Y}=\P(\vec Y(\tau +\delta) \in \cdot\, |(\vec X(\tau), \vec Y(\tau)))$
have a common part $\mu_{\vec X, \vec Y}$ with measure at least $\frac{1}{2}$. From $\mu_{\vec X}\times \mu_{\vec Y}$, we can take out
$\mu_{\vec X, \vec Y}\times \mu_{\vec X, \vec Y}$, which has measure at least $\frac{1}{4}$, and couple $(\vec X(\tau+\delta+t))_{t\geq 0}$ and $(\vec Y(\tau+\delta+t))_{t\geq 0}$
so that they coincide for all $t\geq 0$ and evolve as the strong Markov process defined by (\ref{SDEaut}) with initial measure $\mu_{\vec X, \vec Y}$.
With respect to the remaining measure $\mu_{\vec X}\times \mu_{\vec Y} - \mu_{\vec X, \vec Y}\times \mu_{\vec X, \vec Y}$, we let $(\vec X(\tau+\delta+t))_{t\geq 0}$ and
$(\vec Y(\tau+\delta+t))_{t\geq 0}$ continue to evolve independently.
Since $\mu_{\vec X}\times \mu_{\vec Y} - \mu_{\vec X, \vec Y}\times \mu_{\vec X, \vec Y}$ is absolutely
continuous with respect to $\mu_{\vec X}\times \mu_{\vec Y}$, a.s.\ $(\vec X(\tau+\delta+t), (\vec Y(\tau+\delta+t))_{t\geq 0}$ will visit
$B_\epsilon(\vec x^*)\times B_\epsilon(\vec x^*)$ again. We can therefore iterate the above coupling procedure. Each iteration reduces the probability that
$\vec X$ and $\vec Y$ have not been successfully coupled by a factor $\frac{1}{4}$. Continue the iteration indefinitely to get the desired coupling
between $\vec X$ and $\vec Y$. We comment that, unlike in the context of Harris chains (see e.g.\ Section 5.6 of
Durrett~\cite{Du96}) where one would need $\P(\vec X(\delta)\in \cdot | \vec X(0)=\vec x)$ to be dominated from below by a
positive measure uniformly for $\vec x\in B_\epsilon(\vec x^*)$, to get a successful coupling it suffices that
$\P(\vec X(\delta) \in \cdot |\vec X(0)=\vec x)$ and $\P(\vec X(\delta)\in \cdot |\vec X(0)=\vec y)$ overlap with
probability at least $\alpha$ for some $\alpha>0$ uniformly for all $\vec x,\vec y\in B_\epsilon(\vec x^*)$.

Next we show that, for Lebesgue almost every $\vec x\in [0,\infty)^2$, ${\cal L}(\vec X(t) |\vec X(0)=\vec x) \Longrightarrow \Gamma_{\vec\theta}^{c,g}$
as $t\to\infty$. Let $A=\{ \vec x\in [0,\infty)^2 : {\cal L}\big(\vec X(t)|\vec X(0)=\vec x\big) \not\Longrightarrow \Gamma_{\vec\theta}^{c,g}\}$.
By Theorem \ref{thm:martproblem} and the remark following it, the process defined by (\ref{SDEaut}) is Feller continuous, and therefore $A$
is Borel-measurable. If $A$ has positive Lebesgue measure, then we can find a simply connected bounded open domain $D\subset (0,\infty)^2$ with smooth
boundary such that $A\cap D$ has positive Lebesgue measure. We have shown above that $\Gamma_{\vec\theta}^{c,g}(A)=0$, and hence $\Gamma_{\vec\theta}^{c,g}(A\cap D)=0$.
If $(\vec X(t))_{t\geq 0}$ is the stationary solution of (\ref{SDEaut}) with marginal distribution $\Gamma_{\vec\theta}^{c,g}$, then
$\E[\int_0^T 1_{\vec X(t) \in A\cap D} dt] =0$ for all $T>0$. On the other hand, by Theorem \ref{thm:occupation}, we have for every $\vec x\in D$ that
$\E[\int_0^{\tau_D} 1_{\vec X(t) \in A\cap D}dt\, |\, \vec X(0)=\vec x] >0$. Since $\Gamma_{\vec\theta}^{c,g}$ assigns positive probability to $D$, we have
$$
\int_D \E\left[\int_0^{\tau_D} 1_{\vec X(t) \in A\cap D}dt\, \Big|\, \vec X(0)=\vec x\right] \Gamma_{\vec\theta}^{c,g}(d\vec x) >0.
$$
By the monotone convergence theorem, we can choose $T$ sufficiently large such that
$$
\int_D \E\left[\int_0^{\tau_D\wedge T} 1_{\vec X(t) \in A\cap D}dt\, \Big|\, \vec X(0)=\vec x\right] \Gamma_{\vec\theta}^{c,g}(d\vec x) >0,
$$
the left-hand side of which is in turn dominated by $\E[\int_0^T 1_{\vec X(t) \in A\cap D} dt] =0$, which is a contradiction. Therefore $A$ has Lebesgue
measure 0.

Lastly, we show that ${\cal L}(\vec X(t) | \vec X(0)=\vec x) \Longrightarrow \Gamma_{\vec\theta}^{c,g}$ for all $\vec x\in [0,\infty)^2$. Indeed, for $\vec x\in (0,\infty)^2$, let $\epsilon>0$ be such that $B_{\epsilon}(\vec x) \subset (0,\infty)^2$. By Corollary \ref{cor:domaindiffdensity} applied to $D=B_\epsilon(\vec x)$,
the transition kernel $\mu_t^{B_\epsilon(\vec x)}(\vec x, \cdot)$ with killing at the boundary of $B_\epsilon(\vec x)$ is absolutely continuous with respect to
Lesbesgue measure. Since, for Lebesgue almost every $\vec y \in B_\epsilon(\vec x)$, ${\cal L}(\vec X(t+s) | \vec X(t)=\vec y) \Longrightarrow \Gamma_{\vec\theta}^{c,g}$
as $s\to\infty$ and $\mu_t^{B_\epsilon(\vec x)}(\vec x, B_\epsilon(\vec x))\uparrow 1$ as $t\downarrow 0$ (see (\ref{eq:escapeprob})), we have
${\cal L}(\vec X(t)) | \vec X(0)=\vec x) \Longrightarrow \Gamma_{\vec\theta}^{c,g}$. The case $\vec x\in \partial [0,\infty)^2$ follows from our previous
observation that $\vec X(t)$ starting from $\vec x$ enters $(0,\infty)^2$ instantly (see (\ref{eq:leaveboundary1})--(\ref{eq:leaveboundary2})).

\medskip\noindent
\underline{$\vec\theta \in \partial [0,\infty)^2$}: Without loss of generality we may
assume that $\theta_1=0$. If $X_1(0)=0$, then $X_1(t)=0$ for all $t\geq 0$ and
$(X_1(t), X_2(t)) = (0, X_2(t))$ is effectively a one-dimensional diffusion with
diffusion function $g_2(0, x_2)$. By the same argument as before, albeit much simpler,
this one-dimensional diffusion is ergodic, and the convergence in (\ref{eqconv})
holds. If $X_1(0)\neq 0$, then it suffices to show that $X_1(t)\to 0$ a.s.\ and
$\cL(X_2(t)) \Longrightarrow \Gamma_{\vec\theta}^{c,g}$ as $t\to\infty$, where $\Gamma_{\vec\theta}^{c,g}$
is taken as a measure on $[0,\infty)$.

Note that $X_1(t)$ is a local supermartingale and $X_1(t) \wedge 1$ is a bounded
supermartingale, so that $X_1(t)\wedge 1 \to Y$ a.s. as $t\to\infty$ for some
non-negative random variable $Y$. By the bounded convergence theorem and (\ref{tigh}),
\be{}
\E[Y] = \lim_{t\to\infty} \E[X_1(t)\wedge 1] \leq \lim_{t\to\infty} X_1(0) e^{-ct} = 0 .
\ee
Therefore $Y\equiv 0$ and $X_1(t)\to 0$ a.s.\ as $t\to\infty$.

To show that $\cL(X_2(t)) \Longrightarrow \Gamma_{\vec\theta}^{c,g}$ as $t\to\infty$, it
suffices to show that $\E[ \phi(X_2(t))] \to \E_{\vec\theta}^{c,g}[\phi(X_2)]$ as $t\to\infty$ for any
$\phi \in C^2_c[0,\infty)$. Abbreviate
\be{aludef}
\alpha = \E_{\vec\theta}^{c,g}[\phi(X_2)] \quad \mbox{ and } \quad
u(t, \vec x) = \E[ \phi(X_2(t))\mid \vec X(0) = \vec x].
\ee
For $\vec X(0)\in [0,\infty)^2$ with $X_1(0)=0$, $(\vec X(t))_{t\geq0}$ is effectively a
one-dimensional diffusion that is ergodic, and hence $u(t, \vec x) \to \alpha$
as $t\to\infty$ for each $\vec x \in \{0\}\times [0,\infty)$. We claim that in
fact $u(t, \vec x)\to\alpha$ uniformly on compact intervals of the form
$\{0\}\times [0,K]$. To see why, note that if $Y(t)$ and $Z(t)$ are solutions
of the one-dimensional SDE
\be{Xeq1}
dX(t) = c (\theta_2 - X(t))\ dt + \sqrt{2g_2(0, X(t))}\, dB_t
\ee
with initial condition $Y(0)=y < Z(0)=z$, then $Z(t)$ stochastically dominates $Y(t)$
for all $t\geq 0$, i.e., if $F_{t,y}(v) = \P (Y(t) < v|Y(0)=y)$, then $F_{t,y}(v) \geq F_{t,z}(v)$
for all $t, v\geq 0$. Let $F_{\infty}(v) = \Gamma_{\vec\theta}^{c,g}(-\infty,v)$. Then, for
any $x_2\geq 0$, $F_{t,x_2}(v)\to F_{\infty}(v)$ as $t\to\infty$ for all but countably many $v\in[0,\infty)$.
For any $x_2\in [0,K]$, $K>0$, we can write
\be{urels}
u(t, (0,x_2)) =\!\! \int_0^{\infty}\!\!\!\! \phi(v) dF_{t,x_2}(v)
= -\!\!\int_0^\infty\!\!\!\! \phi'(v) F_{t,x_2}(v) dv
=\!\! \int_0^\infty\!\!\!\! (\phi'_{-}(v) - \phi'_{+}(v)) F_{t,x_2}(v)dv,
\ee
where $\phi'_{+}(v) = \phi'(v) \vee 0$ and $\phi'_{-}(v) = - (\phi'(v) \wedge 0)$. Since
\be{intoinft}
\int_0^\infty \phi'_{-}(v) F_{t,K}(v) dv
\leq \int_0^\infty \phi'_{-}(v) F_{t,x_2}(v) dv
\leq \int_0^\infty \phi'_{-}(v) F_{t,0}(v) dv,
\ee
where both ends of the inequality tend to $\int_0^\infty \phi'_{-}(v) F_{\infty}(v)dv$
by the bounded convergence theorem, $\int_0^\infty \phi'_{-}(v) F_{t,x_2}(v) dv$
converges uniformly to $\int_0^\infty \phi'_{-}(v) F_{\infty}(v)dv$ for $x_2\in [0,K]$ as $t\to\infty$.
A similar statement holds for $\int_0^\infty \phi'_{+}(v) F_{t,x_2}(v) dv$. Therefore
$u(t, \vec x)$ converges uniformly to $\alpha$ on $\{0\}\times [0,K]$.

Let $\vec X(0) \in [0,\infty)^2$ be arbitrary. By (\ref{tigh}), $(X_2(t))_{t\geq 0}$
is tight, and hence for any $\epsilon>0$ we can choose $K$ large enough so that
$\P(X_2(t)>K)\leq \epsilon$ for all $t\geq 0$. Since $u(t, \vec x)\to \alpha$
uniformly on $\{0\}\times [0,K]$, we can choose $t_1$ large enough so that
$\sup_{x_2\in[0,K]} |u(t_1, (0,x_2))-\alpha|\leq \epsilon/2$. Since
$\big\{(\vec X(t))_{t\geq 0}\big\}_{\vec X(0)\in[0,\infty)^2}$ defines a
Feller process (see the remark below Theorem \ref{thm:martproblem}),
$u(t_1, \vec x)$ is continuous in $\vec x\in[0,\infty)^2$. We can therefore choose
$\delta>0$ sufficiently small so that $\sup_{\vec x\in [0,\delta]\times [0,K]}
|u(t_1, \vec x)-\alpha|\leq\epsilon$. Since $X_1(t)\to 0$ a.s., we can choose
$t_2$ large enough so that $\P(X_1(t)>\delta) \leq \epsilon$ for all $t\geq t_2$.
Then, by the Markov property, for any $t\geq t_1+t_2$ we have
\be{uasym}
\begin{aligned}
u(t, \vec X(0)) &= \E[u(t_1, \vec X(t-t_1))]\\
&= \E\Big[u(t_1, \vec X(t-t_1)) 1_{\vec X(t-t_1)\in[0,\delta]\times[0,K]}\Big]\\
&\qquad + \E\Big[u(t_1, \vec X(t-t_1)) 1_{\vec X(t-t_1)\notin[0,\delta]\times[0,K]}\Big].
\end{aligned}
\ee
Since $\P\big(\vec X(t-t_1)\notin[0,\delta]\times[0,K]\big) \leq 2\epsilon$ and
$\|u\|_\infty\leq \|\phi\|_\infty$, $\alpha \leq \|\phi\|_\infty$, we easily verify from (\ref{uasym})
that
\be{ualmin}
| u(t, \vec X(0))-\alpha | \leq \epsilon
+ 4 \epsilon \|\phi\|_\infty  \quad \mbox{ for all }\ t\geq t_1+t_2.
\ee
Since $\epsilon>0$ is arbitrary, $u(t,\vec X(0))\to \alpha$ as $t\to\infty$, and
hence $\cL(\vec X(t)) \Longrightarrow \Gamma_{\vec\theta}^{c,g}$.
\epr

\subsection{Proof of Theorems \ref{thm:FCdomain}, \ref{thm:fixshap} and
Corollary \ref{cor:iterfixshape}}
\label{S3.3}

{\bf Proof of Theorem \ref{thm:FCdomain}.} Let $g=(g_1, g_2) \in \cH_a$ for some
$0\leq a<c$. Then, by (\ref{Ha}), there exists a $0<C=C(g)<\infty$ such that
\be{gsumbds}
g_1(\vec x)+ g_2(\vec x) \leq C(1+x_1)(1+x_2) + a(x_1^2+x_2^2),
\qquad (x_1,x_2) \in [0,\infty)^2.
\ee
The finiteness of $F_c g$ follows from Proposition \ref{prop:momenteqn}(ii). If
$\vec \theta_n\to \vec\theta$ for some $\vec\theta\in [0,\infty)^2$, then, by Proposition
\ref{prop:momenteqn}(iii), $g_1, g_2$ are uniformly integrable with respect to
$\{\Gamma^{c,g}_{\vec\theta_n}\}_{n\in\N}$. Combining this with the fact, shown in Theorem
\ref{thm:SDEeq} and proved in Section \ref{S3.2}, that $\Gamma_{\vec\theta_n}^{c,g}$
converges weakly to $\Gamma_{\vec\theta}^{c,g}$ as $\vec\theta_n\to\vec\theta$, we have
$\E_{\vec\theta_n}^{c,g} [ g_i(\vec X)] \to \E_{\vec \theta}^{c,g} [ g_i(\vec X)]$,
i.e., $(F_c g)_i(\vec\theta_n) \to (F_c g)_i(\vec \theta)$ for $i=1, 2$ (recall
(\ref{transform})).

By the moment equations (\ref{eq:1stmoment}--\ref{eq:mixmoment}), we have \be{Fgids}
\begin{aligned}
(F_cg)_1(\vec \theta) + (F_cg)_2(\vec \theta)
& =  \E_{\vec\theta}^{c,g} [g_1(\vec X) + g_2(\vec X)] \\
& \leq  \E_{\vec\theta}^{c,g} [C(1+X_1)(1+X_2) + a(X_1^2+X_2^2)] \\
& = C(1+\theta_1)(1+\theta_2) + a(\theta_1^2+\theta_2^2) + \frac{a}{c}\big((F_cg)_1(\vec\theta) +
(F_cg)_2(\vec\theta)\big).
\end{aligned}
\ee
Therefore
\be{sumFcsbds}
(F_c g)_1(\vec\theta) + (F_c g)_2(\vec\theta)
\leq \frac{c}{c-a} \Big( C(1+\theta_1)(1+\theta_2) + a(\theta_1^2+\theta_2^2) \Big).
\ee
Consequently, if $g\in \cH_{0^+}$, then $F_cg$ satisfies (\ref{sumFcsbds}) for all $a>0$,
and so it satisfies the subquadratic growth bound imposed by the class $\cH_{0^+}$.

To show $\partial F_cg = \partial g$, note that $F_cg\geq 0$ is obvious. If $\vec\theta\in
(0,\infty)^2$, then the equilibrium distribution $\Gamma_{\vec\theta}^{c,g}$ has positive
mass in $(0,\infty)^2$, and so $(F_cg)(\vec\theta)>0$ follows from the fact that $g>0$
on $(0,\infty)^2$. If $\theta_1=0$, then, by (\ref{eq:1stmoment}), $\Gamma_{\vec\theta}^{c,g}$
is concentrated on the vertical axis $A_2$. Since $g_1$ vanishes on $A_2$, it follows that
$(F_cg)_1(\vec\theta) = 0$. Moreover, $(F_cg)_2(\vec\theta)=0$ if and only if $g_2$ vanishes
on $A_2$ (recall (\ref{partgdef}--\ref{bdpos})). A similar result holds for $\theta_2=0$.
\qed

\medskip\noindent
{\bf Proof of Theorem \ref{thm:fixshap}.}
Theorem \ref{thm:fixshap}\,(i) follows immediately from (\ref{eq:1stmoment}--\ref{eq:mixmoment}).
To prove Theorem \ref{thm:fixshap}\,(ii), note that, by (\ref{eq:1stmoment}--\ref{eq:2ndmoment}),
\beq
(F_c g)_1 (\vec\theta)
&=& \E_{\vec\theta}^{c,g}[a_1 X_1^2 + b_1 X_1 + c_1 X_1 X_2 ]
= a_1 \E_{\vec\theta}^{c,g}[X_1^2] + b_1\theta_1 + c_1\theta_1\theta_2
\nonumber \\
&=& a_1 \theta_1^2 + \frac{a_1}{c} (F_c g)_1 + b_1\theta_1 + c_1\theta_1\theta_2 = g_1(\vec \theta) + \frac{a_1}{c} (F_c g)_1(\vec\theta).
\eeq
Solving for $(F_cg)_1(\vec\theta)$, we get $(F_c g)_1(\vec\theta) = \frac{c}{c-a_1}
g_1(\vec\theta)$. Similarly, we have $(F_c g)_2 = \frac{c}{c-a_2} g_2$ for
$g_2 = a_2 x_2^2 + b_2 x_2 + c_2 x_1x_2$. The assumption $(b_1+c_1)(b_2+c_2)>0$ is meant to rule out
the uninteresting case $g_1=0$ or $g_2=0$.
\qed

\medskip\noindent
{\bf Proof of Corollary \ref{cor:iterfixshape}.}
Equation (\ref{iterfs}) follows from Theorem \ref{thm:fixshap}\,(ii) by induction.
Note that if $\alpha_i \sum_{k=0}^{n_0-1} c_k^{-1} \geq 1$ for either $i=1$ or 2, then the coefficient of
$x_i^2$ in $(F^{[n_0-1]}g)_i(\vec x)$ is $\alpha_i /[1-\alpha_i \sum_{k=0}^{n_0-2}c_k^{-1}] \geq c_{n_0-1}$. To show
$(F^{[n_0]}g)_1 +(F^{[n_0]}g)_2 = \infty$ on $(0,\infty)^2$, it therefore
suffices to show $(F_cg)_1+(F_cg)_2 \equiv \infty$ on $(0,\infty)^2$ for $g$
of the form $g_i(\vec x) = \alpha_i x_i^2 +\beta_i x_i +\gamma_i x_1x_2$ with $\alpha_1\vee \alpha_2\geq c$.
Without loss of generality, assume $\alpha_1 \geq c$. The proof of Proposition \ref{prop:momenteqn}\,(ii) shows
that the moment equations (\ref{eq:1stmoment}--\ref{eq:2ndmoment}) are valid as
long as $(F_cg)_1(\vec \theta) +(F_cg)_2(\vec \theta) = \E_{\vec \theta}^{c,g}[g_1+g_2]
<\infty$. Assume $(F_cg)_1(\vec\theta)+(F_cg)_2(\vec\theta) < \infty$ for some
$\vec\theta\in (0,\infty)^2$. Then
\be{2ndmomest}
\E_{\vec\theta}^{c,g} [X_1^2] = \theta_1^2
+ \frac{1}{c} (F_cg)_1(\vec\theta) = \theta_1^2
+ \frac{\alpha_1}{c}\E_{\vec\theta}^{c,g}[X_1^2]
+ \frac{\beta_1}{c}\theta_1 +\frac{\gamma_1}{c}\theta_1\theta_2,
\ee
which is not possible for $\alpha_1\geq c$. Therefore we must have $(F_cg)_1(\vec\theta)
+(F_cg)_2(\vec\theta) = \infty$ for all $\vec\theta\in (0,\infty)^2$.
\qed

\subsection{Discussion of Conjecture \ref{conj:FCbdprop}}
\label{S3.4}

In this section we explain why Conjecture \ref{conj:FCbdprop} is plausible.
We focus on the case where $g_1,g_2$ both satisfy boundary property $(\partial_{12})$
in (\ref{bdprops}), i.e., $g_1(\vec x) = x_1x_2\gamma_1(\vec x)$ and $g_2(\vec x)
= x_1x_2\gamma_2(\vec x)$ with $\gamma_1,\gamma_2>0$ continuous on $[0,\infty)^2$.

Consider the tilted equilibrium
\be{eqtilt}
\hat\Gamma^{c,g}_{\vec\theta}(d\vec x)
= \frac{x_1x_2}{\theta_1\theta_2}\Gamma^{c,g}_{\vec\theta}(d\vec x),
\qquad \vec\theta \in (0,\infty)^2,
\ee
where (\ref{eq:mixmoment}) implies the proper normalization. The conjecture amounts
to showing that, as $\vec\theta\to\vec\theta^*\in\partial[0,\infty)^2$, this
tilted equilibrium converges weakly to some probability distribution on $[0,\infty)^2$,
say $\hat\Gamma_{\vec \theta^*}^{c,g}(d\vec x)$, that is weakly continuous
in $\vec\theta^*$ and, in addition, $\gamma_i(\vec x)$ is uniformly integrable with
respect to $\hat\Gamma_{\vec \theta}^{c,g}(d\vec x)$ for $\vec\theta$ in a small neighborhood
of $\vec \theta^*$. Indeed, this observation is immediate from the identity
\be{inttilt}
\int_{[0,\infty)^2} \gamma_i(\vec x)\,\hat\Gamma^{c,g}_{\vec\theta}(d\vec x)
= \frac{1}{\theta_1\theta_2}\,(F_cg)_i(\vec\theta), \qquad i=1,2.
\ee

Now, recalling the generator in (\ref{generator}), we note that $\hat\Gamma^{c,g}_{\vec\theta}
(d\vec x)$ is the equilibrium associated with the time-changed diffusion given by the generator
\begin{eqnarray}
\!\!\!\!\!&&\!\!\!\!\! (\hat L_{\vec\theta}^{c,g} f)(\vec x)
= \frac{c(\theta_1-x_1)}{x_1x_2} \frac{\partial}{\partial x_1} f(\vec x)
+ \frac{c(\theta_2-x_2)}{x_1x_2} \frac{\partial}{\partial x_2} f(\vec x)
+ \gamma_1(\vec x) \frac{\partial^2}{\partial x_1^2} f(\vec x)
+ \gamma_2(\vec x) \frac{\partial^2}{\partial x_2^2} f(\vec x),\nonumber \\
&&\!\!\!\!f\in C^2_c([0,\infty)^2),\quad (\vec\theta -\vec x)\cdot \nabla f(\vec x) =0 \,\mbox{ on }
\,\partial[0,\infty)^2. \label{gentimechd}
\end{eqnarray}

Let $\vec\theta\to\vec\theta^*=(\alpha,0)$ for some $\alpha>0$. Then, at least heuristically,
we get a limiting generator
\be{limgen}
\begin{aligned}
&(\hat L_{(\alpha,0)}^{c,g} f)(\vec x)
= \frac{c(\alpha-x_1)}{x_1x_2} \frac{\partial}{\partial x_1} f(\vec x)
- \frac{c}{x_1} \frac{\partial}{\partial x_2} f(\vec x)
+ \gamma_1(\vec x) \frac{\partial^2}{\partial x_1^2} f(\vec x)
+ \gamma_2(\vec x) \frac{\partial^2}{\partial x_2^2} f(\vec x),\\
&f\in C^2_c([0,\infty)^2),\, (\vec\theta^*-\vec x) \cdot \nabla
f(\vec x) =0  \,\mbox{ on }
\,\partial[0,\infty)^2\backslash\{\vec\theta^*\}, \,
\frac{\partial}{\partial x_1} f(\vec\theta^*) = \frac{\partial}{\partial x_2} f(\vec\theta^*) = 0.
\end{aligned}
\ee
Here, the diffusion part has no singularity at the boundary, but the drift
part does. As the process approaches the vertical axis $A_2$ it feels a growing drift downwards
and to the right, while as it approaches the horizontal axis $A_1$ it feels a growing drift
horizontally towards $(\alpha,0)$ and a constant drift downwards. Therefore,
again heuristically, this generator describes a process that is {\em obliquely reflected
in the direction of $(\alpha,0)$ upon hitting $A_2$, and upon hitting $A_1$
jumps to $(\alpha,0)$ instantly and then moves back into the interior
by reflection}. Like the original diffusion with generator (\ref{generator}),
this process ought to exist, be weakly unique, and have an ergodic equilibrium
$\hat \Gamma_{\vec\theta^*}^{c,g}$ that is weakly continuous in $\vec\theta^*\in [0,\infty)^2$.


\section{Proof of Theorems \ref{thm:upfs}, \ref{thm:downfs} and \ref{thm:FPS}}
\label{S4}

Section \ref{S4.1} contains the proof of Theorem~\ref{thm:FPS}, which is an immediate consequence of
Proposition~\ref{prop:chainharm}. Section \ref{S4.2} contains some preliminary lemmas needed for the proof
of Proposition \ref{prop:chainharm}. Section \ref{S4.3} provides the proof of Proposition \ref{prop:chainharm}
and of Theorems \ref{thm:upfs} and \ref{thm:downfs}.

\subsection{Proof of Theorem \ref{thm:FPS}}
\label{S4.1}

The proof of Theorem \ref{thm:FPS} is based on an asymptotic analysis of the
\emph{homogeneous} Markov chain $\vec M^{c,g}=(\vec M^{c,g}(n))_{n\in\N_0}$ with
transition probability kernel given by $p(\vec\theta, d\vec y) = \Gamma_{\vec\theta}^{c,g}
(d\vec y)$, the unique equilibrium distribution of (\ref{SDEaut}). For $F_c g = g$,
$\vec M^{c,g}$ is in fact the {\it interaction chain} in (\ref{MCback}). Throughout the rest of the
section, unless specified otherwise, we will denote the Markov chain $\vec M^{c,g}$ by $\vec X$. For $F_c g = g$
and $g\in \cH_{0^+}^r$, both $g_1$ and $g_2$ are harmonic functions of $\vec X$, i.e., both $(g_1(\vec
X(n)))_{n\in\N_0}$ and $(g_2(\vec X(n)))_{n\in\N_0}$ are martingales. Theorem \ref{thm:FPS} then
follows immediately from the following proposition.

\bp{prop:chainharm}{\bf [Harmonic functions of $\vec X=\vec M^{c,g}$]}\\
If $g\in\cH_{0^+}$ and satisfies {\rm (\ref{gcon1})} in the definition of $\cH_0^r$,
then every nonnegative harmonic function $f$ of $\vec X = \vec M^{c,g}$, i.e., every $f$
such that
\be{fmartrel}
\E\Big[f(\vec X(n))\,\Big|\, \vec X(0)=\vec \theta\Big] = f(\vec \theta)
\qquad\forall\,\vec\theta\in [0,\infty)^2,\,n\in\N_0,
\ee
which furthermore satisfies the constraints
\begin{eqnarray}
&(i)&  f(\vec x) \leq C(1+x_1)(1+x_2) \quad \mbox{ for some } 0<C=C(f)<\infty,
\label{constrain1} \\[0.3cm]
&(ii)& \lim_{\vec x\to\vec z} f(\vec x) = 0 \qquad \forall\, \vec z \in \partial g,
\label{constrain2} \\
&(iii)& \lim_{\vec x\to\vec z} \frac{f(\vec x)}{h_{\vec z}(\vec x)}
= \lambda_{f,\vec z} \in [0,\infty) \qquad \forall\, \vec z \in R_\infty,
\label{constrain3}
\end{eqnarray}
is of the form
\be{fform}
f(\vec x) = \sum_{\vec z\in R_\infty} \lambda_{f,\vec z} h_{\vec z}(\vec x)
= \lambda_{f,(\infty,0)}\, x_1 + \lambda_{f,(0,\infty)}\, x_2
+ \lambda_{f,(\infty,\infty)}\, x_1 x_2,
\ee
with $h_{\vec z}$, $\vec z\in R_\infty$, given by {\rm (\ref{Rdefs})--(\ref{hdefs})}.
\ep

The proof of Proposition~\ref{prop:chainharm} will be given in
Section~\ref{S4.3}. The strategy is to first $h$-transform
$\vec X$ (see Definition \ref{htransform} below) to a new process $\vec X^h=(\vec X^h(n))_{n\in\N_0}$ using
\be{}
h(\vec x) = (1+x_1)(1+x_2),
\ee
i.e., $\vec X^h$ is defined as the homogeneous Markov chain with
transition probability kernel
$$
p(\vec\theta,d\vec y)= h(\vec y)\Gamma_{\vec\theta}^{c,g} (d\vec y)/h(\vec \theta),
$$
which is well-defined since $h(\vec x)$ is a harmonic function
of $\vec M^{c,g}$. The function $f$ is harmonic for $\vec M^{c,g}$ if and only if
$f/h$ is harmonic for $\vec X^h$. The constraint in (\ref{constrain1}) guarantees
that $f/h$ is bounded, the constraints in (\ref{constrain2}--\ref{constrain3})
guarantee that $f/h$ is continuous up to the boundary
\be{Rdef}
R = \partial g \cup R_\infty,
\ee
while the constraint in (\ref{gcon1}) guarantees that $\lim_{n\to\infty} \vec X^h(n)
\in R$ a.s. It is then standard to show that $f/h$ is uniquely determined by its
values at $R$, which will imply (\ref{fform}).

The proofs of Theorems \ref{thm:upfs} and \ref{thm:downfs} are also based on an
asymptotic analysis of the Markov chain $\vec M^{c,g}$, even though when $g$ is
not a fixed point of $F_c$, it no longer corresponds to the interaction chain in (\ref{MCback}).

\subsection{Preliminary lemmas}
\label{S4.2}

The key results in this section are Proposition \ref{prop:hlimitchain} and
Corollary \ref{cor:trapprob}.

Let $\vec X=\vec M^{c,g}$ be as stated before Proposition \ref{prop:chainharm}. First we list some moment
equations for $\vec X(n)$, $n\in \N_0$, which follow immediately from Proposition \ref{prop:momenteqn}.

\bl{lem:itermomeqn} {\bf [Moment equations for $\vec X= \vec M^{c,g}$]}\\
Let $c>0$, and $g\in \cH_a$ for some $0\leq a<c$. Fix $\vec X(0)=\vec\theta \in [0,\infty)^2$. Then for all $n\in\N_0$,
\begin{eqnarray}
&& \E[X_i(n)] = \theta_i,
\qquad \qquad \qquad \qquad \qquad i= 1, 2,\label{eq:1stitermom}\\
&& \E[X_1(n) X_2(n)] = \theta_1 \theta_2. \label{eq:mixitermom}
\end{eqnarray}
If $((F_cg)_1, (F_cg)_2) = (\lambda_1 g_1, \lambda_2 g_2)$ for some
$\lambda_1, \lambda_2>0$, then
\begin{eqnarray}
&& \E[g_i(\vec X(n))] = \lambda_i^n g_i(\vec \theta),
\qquad\qquad\qquad \ i=1,2, \label{eq:iterfp} \\
&& \E[X^2_i(n)] = \theta_i^2
+ \frac{1}{c} \sum_{j=1}^n \lambda_i^j g_i(\vec\theta), \quad \qquad i=1,2.
\label{eq:2nditermom}
\end{eqnarray}
\el

In the proof of Theorem \ref{thm:FPS}, we will need Doob's $h$-transform of a Markov
chain, which we recall here. For more information on the
$h$-transform, see e.g.\ Section 4.1 of Pinsky~\cite{P95}.

\bd{htransform}{\bf [$h$-transform]} \\
Let $X=(X(n))_{n\in\N_0}$ be a Markov chain with state space $E$ and $n$-step transition
probability kernel $p_n(x,dy)$. If $h$ is a nonnegative (not identically zero) harmonic
function of $X$, i.e., $(h(X(n)))_{n\in\N_0}$ is a nonnegative martingale, then
the $h$-transform of $X$, denoted by $X^h$, is defined as the Markov chain on the
space $\{x \in E\colon\, h(x) > 0\}$ with $n$-step transition probability kernel
$p^h_n(x,dy) = p_n(x,dy) h(y)/h(x)$.
\ed

\noindent
The next two lemmas are immediate consequences of Definition \ref{htransform}.

\bl{lem:htransform}{\bf [Harmonic functions of $X^h$]} \\
Let $X$, $h$ and $X^h$ be as in Definition {\rm \ref{htransform}}. If $f$ is a harmonic
function of $X$, then $f/h$ restricted to $\{x\in E\colon\, h(x)>0\}$ is a harmonic function of $X^h$.
The converse is true if $h(x)>0$ for all $x\in E$.
\el

\bl{lem:densityhtransform}{\bf [Absolute continuity of $X^h$ w.r.t. $X$ at bounded stopping times]}\\
Let $X$, $h$ and $X^h$ be as in Definition {\rm \ref{htransform}}. If $X(0)=X^h(0)=x\in E$ where
$h(x)>0$, and $\tau$ is a bounded stopping time, then the law of $X^h(\tau)$ is absolutely continuous
with respect to the law of $X(\tau)$ with density $\frac{h(\cdot)}{h(x)}$.
\el

\noindent The next proposition is the key to establishing Proposition
\ref{prop:chainharm}. Such a result is referred to as almost sure {\em extinction versus unbounded growth}, see e.g.\ Fleischmann and Swart~\cite{FS06}.

\bp{prop:hlimitchain}{\bf [Almost sure limit of $h$-transform of $\vec X=\vec M^{c,g}$]} \\
Let $c>0$, and let $g\in \cH_{0^+}$ satisfy condition {\rm (\ref{gcon1})}. Let $h(\vec x) = (1+x_1)(1+x_2)$
and let $\vec X^h$ be the $h$-transform of $\vec X$. Then, for any $\vec X^h(0)\in [0,\infty)^2$, almost surely,
$\lim_{n\to\infty} \vec X^h(n) =\vec X^h(\infty)$ exists and $\vec X^h(\infty) \in R$ $($see $(\ref{Rdef}))$.
\ep
Before giving the proof of Proposition~\ref{prop:hlimitchain}, which we defer to the end of this subsection,
we first state and prove a corollary and another prerequisite lemma.

\bc{cor:trapprob}{\bf [Trapping probabilities]} \\
Let $c$, $g$, $h$, $\vec X^h$ and $\vec X^h(\infty)$ be as in Proposition {\rm \ref{prop:hlimitchain}}.\\
$(i)$
\be{inftyinfty}
\P[\vec X^h(\infty) = (\infty, \infty)]
= \frac{X^h_1(0)X^h_2(0)}{(1+X^h_1(0))(1+X^h_2(0))}.
\ee
$(ii)$ If $(0,\infty)\times \{0\} \notin \partial g$, then
\be{inftyzero}
\P[\vec X^h(\infty) = (\infty, 0)] = \frac{X^h_1(0)}{(1+X^h_1(0))(1+X^h_2(0))}.
\ee
$(iii)$ If $\{0\}\times (0,\infty) \notin \partial g$, then
\be{zeroinfty}
\P[\vec X^h(\infty) = (0,\infty)] = \frac{X^h_2(0)}{(1+X^h_1(0))(1+X^h_2(0))}.
\ee
\ec

\bpr
By Lemmas \ref{lem:itermomeqn} and \ref{lem:htransform},
\be{fsext}
f_1(\vec x) = \frac{x_1x_2}{(1+x_1)(1+x_2)},
\quad  f_2(\vec x) = \frac{x_1}{(1+x_1)(1+x_2)},
\quad  f_3(\vec x) = \frac{x_2}{(1+x_1)(1+x_2)},
\ee
are bounded harmonic functions of $\vec X^h$, and therefore $(f_i(\vec X^h(n)))_{n\in\N_0}$,
$i=1,2,3$, are bounded martingales. Since, by Proposition \ref{prop:hlimitchain},
$\vec X^h(n) \to \vec X^h(\infty) \in R$ a.s.\ as $n\to\infty$, we have
\be{fiexp}
f_i(\vec X^h(0)) = \E_{\vec\theta}^{c,g}[ f_i(\vec X^h(\infty))], \qquad i=1,2,3.
\ee
Now (\ref{inftyinfty}--\ref{zeroinfty}) follow from the following observations: (1)
$f_1((\infty,\infty))=1$ and $f_1= 0$ on $R\backslash\{(\infty,\infty)\}$; (2) if
$(0,\infty)\times \{0\} \notin \partial g$, then $f_2((\infty,0)) = 1$ and $f_2=0$
on $R\backslash \{(\infty,0)\}$; (3) if $\{0\}\times (0,\infty) \notin \partial g$,
then $f_3((0,\infty)) = 1$ and $f_3=0$ on $R\backslash\{(0,\infty)\}$.
\epr

The proof of Proposition \ref{prop:hlimitchain} in turn relies on the next lemma, which gives a
lower bound for $\hat\Gamma_{\vec\theta, h}^{c,g}(d\vec x) =\Gamma_{\vec\theta}^{c,g}(d\vec x)
h(\vec x)/h(\vec \theta)$, the transition kernel of $\vec X^h$ with $h(\vec x)=(1+x_1)(1+x_2)$, that is uniform in both $g$ and $\vec\theta$. The uniformity in $g$ is not needed for the proof of Proposition \ref{prop:hlimitchain},
but will be crucial for the proof of Theorem \ref{thm:DA} in Section \ref{S5}.

\bl{lem:kerest}{\bf [Uniform lower bound on $\hat\Gamma_{\vec\theta,h}^{c,g}(d \vec x)$]} \\
Let ${\cal A} \subset \cH_{0^+}$.\\
$(i)$ For any $\vec\theta \in [0,\infty)^2$, if
\begin{equation}
\label{kbcond2}
\exists\ \epsilon'>0 \mbox{ such that }
\inf_{g\in {\cal A} \atop \vec x \in B_{\epsilon'}(\vec \theta)} g_i(\vec x) >0
\quad \mbox{ for } i=1 \mbox{ or } i=2
\end{equation}
with $B_{\epsilon'} (\vec \theta) = \{\vec x \in [0,\infty)^2\colon\, \Vert \vec x -\vec
\theta\Vert \leq \epsilon'\}$, then
\begin{equation}
\label{kb1}
\exists\ \epsilon>0 \mbox{ such that } \inf_{g\in {\cal A} \atop
\vec x\in B_{\epsilon}(\vec \theta)} \hat\Gamma_{\vec x,h}^{c,g}\left([0,\infty)^2\backslash
B_{\epsilon}(\vec \theta)\right) > 0.
\end{equation}
$(ii)$ For any $\alpha > 0$, if
\begin{equation}
\label{kbcond1}
\exists\ \epsilon', N'>0 \mbox{ such that }
\inf_{g\in {\cal A} \atop \vec x\in [N',\infty)\times[\alpha-\epsilon',\alpha+\epsilon']}
g_2(\vec x) >0
\end{equation}
and
\be{kbcond3}
\begin{aligned}
&\forall\, a>0, \ \exists\, C_a \in [0, \infty) \mbox{ such that, uniformly for all $\vec x\in [0,\infty)^2$ and $g\in\cal A$,  }\\
&\qquad \qquad g_1(\vec x) + g_2(\vec x) \leq C_a(1+x_1)(1+x_2) + a(x_1^2+x_2^2),
\end{aligned}
\ee
then
\begin{equation}
\label{kb2}
\exists\, \epsilon,\,N>0 \mbox{ such that }
\inf_{g\in {\cal A} \atop \vec x\in [N,\infty) \times [\alpha-\epsilon,\alpha+\epsilon]}
\hat\Gamma_{\vec x,h}^{c,g}\left([0,\infty)^2 \backslash [N,\infty) \times
[\alpha-\epsilon,\alpha+\epsilon] \right) > 0.
\end{equation}
A statement similar to $(\ref{kb2})$ holds for vertical strips of the form $[\alpha -\epsilon,
\alpha+\epsilon]\times [N,\infty)$ if, in $(\ref{kbcond1})$, $g_2$ is replaced by $g_1$ and
$[N',\infty)\times [\alpha-\epsilon',\alpha+\epsilon']$ is replaced by
$[\alpha-\epsilon',\alpha+\epsilon'] \times [N',\infty)$.
\el

\bpr
We first prove (\ref{kb1}) and (\ref{kb2}) with $\hat\Gamma^{c,g}_{\vec x,h}$ replaced by
$\Gamma^{c,g}_{\vec x}$. The main tool is the following moment equation valid for $g\in \cH_{0^+},
\vec \theta\in [0,\infty)^2$ and $i=1,2$:
\begin{equation}
\label{expmomeqn}
\E_{\vec \theta}^{c,g}\left[\frac{1}{(1+X_i)^2}\right] = \frac{1}{1+\theta_i} \E_{\vec \theta}^{c,g}\left[\frac{1}{1+X_i}\right]
+ \frac{2}{c(1+\theta_i)}\E_{\vec x}^{c,g}\left[\frac{g_i(\vec X)}{(1+X_i)^3}\right],
\end{equation}
where $\vec X=(\vec X(t))_{t\geq 0}$ in this proof denotes the stationary solution of the SDE (\ref{SDEaut}). By
stationarity, $\cL(\vec X(s)) = \Gamma_{\vec \theta}^{c,g}$ for all $s\geq 0$. Hence
\begin{equation}
\label{Mdef}
M_i(t) = \frac{1}{1+X_i(t)} - \frac{1}{1+X_i(0)} - \int_0^t L_{\vec\theta}^{c,g}\!\!\left(\frac{1}{1+x_i}\right)\Bigg|_{\vec
  x= \vec X(s)}ds, \quad i=1,2,
\end{equation}
are local martingales, where
\be{Ldesfext}
L_{\vec \theta}^{c,g} = c(\theta_1-x_1)\frac{\partial} {\partial x_1}
+ c(\theta_2-x_2)\frac{\partial}{\partial x_2} + g_1(\vec x)\frac{\partial^2}{\partial x_1^2}
+ g_2(\vec x)\frac{\partial}{\partial x_2^2}.
\ee
Since $\E_{\vec \theta}^{c,g}[X_i(s)] = \theta_i$ and $\E_{\vec \theta}^{c,g}[g_i(\vec X(s))]
= (F_cg)_i(\vec \theta)<\infty$ by Proposition \ref{prop:momenteqn},
we have
\begin{eqnarray}
\E_{\vec\theta}^{c,g}\Big[\sup_{0\leq s\leq t}
|M_i(s)|\Big] &\leq& 2 + \E_{\vec\theta}^{c,g}\Big[\int_0^t \left(c|\theta_i-X_i(s)|
+ 2g_i(\vec X(s))\right) ds \Big]
\nonumber \\
&\leq & 2+ 2 t \Big(c\theta_i + (F_cg)_i(\vec \theta)\Big) < \infty.
\end{eqnarray}
Therefore $M_i=(M_i(t))_{t\geq 0}$, $i=1,2$, are in fact martingales, and $\E_{\vec\theta}^{c,g}[M_i(t)]=0$.
By the stationarity of $\vec X$, we have
\begin{equation}
\label{Lzero}
\E_{\vec \theta}^{c,g}\!\!\left[L_{\vec \theta}^{c,g}\!\!\left(\frac{1}{1+x_i}\right)\!\Bigg|_{\vec x=\vec X(s)}\right]
= \E_{\vec \theta}^{c,g}\!\!\left[-c\cdot\frac{1+\theta_i-1-X_i}{(1+X_i)^2} + \frac{2g_i(\vec X)}{(1+X_i)^3}\right]=0, \quad i=1, 2.
\end{equation}
Rearranging terms, we obtain (\ref{expmomeqn}).

\medskip\noindent
\underline{(\ref{kb1})}: Suppose that (\ref{kb1}) with $\hat \Gamma_{\vec x,h}^{c,g}$ replaced by
$\Gamma_{\vec x}^{c,g}$ is false. Then
\begin{equation}
\inf_{g\in {\cal A}\atop \vec x\in B_{\epsilon}(\vec \theta)}
\Gamma_{\vec x}^{c,g}\left([0,\infty)^2\backslash B_{\epsilon}(\vec \theta)\right)=0
\qquad \forall\epsilon>0.
\end{equation}
By (\ref{kbcond2}), we may assume without loss of generality that $\inf_{g\in {\cal A}, \vec x
\in B_{\epsilon_0}(\vec \theta)} g_1(\vec x) =\delta >0$ for some $\epsilon_0>0$. In particular,
$\inf_{g\in {\cal A}, \vec x \in B_{\epsilon}(\vec \theta)} g_1(\vec x) \geq \delta $ for all
$\epsilon\in [0,\epsilon_0]$. Fix $\epsilon\in  [0,\epsilon_0]$. Let $\vec x^{(n)}\in
B_\epsilon(\vec\theta)$ and $g^{(n)}\in {\cal A}$ be chosen such that $\Gamma_{\vec x^{(n)}}^{c,g^{(n)}}
([0,\infty)^2\backslash B_\epsilon(\vec \theta))= o(1)$ as $n\to\infty$. In (\ref{expmomeqn}) with
$i=1$, substitute $\vec x^{(n)}$ and $g^{(n)}$ for $\vec \theta$ and $g$.
Then
\be{lhsrhs}
\begin{aligned}
\mbox{l.h.s.} &\leq \frac{1}{(1+\theta_1-\epsilon)^2} + o(1), \\
\mbox{r.h.s.} &\geq \frac{1}{(1+\theta_1)(1+\theta_1+\epsilon)} + \frac{2}{c(1+\theta_1)} \times (1-o(1))\times \frac{\delta}{(1+\theta_1+\epsilon)^3},
\end{aligned}
\ee
where we applied Jensen's inequality to obtain $\frac{1}{(1+\theta_1)^2}$ in the estimate for the r.h.s.
For $\epsilon>0$ sufficiently small and $n$ sufficiently large, the above two equations are
incompatible, and therefore (\ref{kb1}) with $\hat \Gamma_{\vec x,h}^{c,g}$ replaced by
$\Gamma_{\vec x}^{c,g}$ holds. Since $h(\vec x) = (1+x_1)(1+x_2)\geq 1$ on $[0,\infty)^2$
and is bounded on $B_\epsilon(\vec\theta)$, it is easy to see by the definition
of $\hat \Gamma^{c,g}_{\vec x, h}$ that (\ref{kb1}) also holds .

\medskip\noindent
\underline{(\ref{kb2})}: The proof that (\ref{kb2}) holds with $\hat \Gamma_{\vec x,h}^{c,g}$
replaced by $\Gamma_{\vec x}^{c,g}$ is the same as above and we leave the details to the reader.
To get (\ref{kb2}), we argue as follows.

Choose $\epsilon \in (0,\alpha)$ and $N_0>0$ such that
\be{betarels}
\beta_{\alpha,\epsilon, N_0} = \inf_{g\in {\cal A} \atop \vec x\in [N_0,\infty)\times
[\alpha-\epsilon,\alpha+\epsilon]} \Gamma_{\vec x}^{c,g}\left([0,\infty)^2 \backslash [N_0,\infty)\times
[\alpha-\epsilon,\alpha+\epsilon] \right)>0.
\ee
By Proposition \ref{prop:momenteqn}, we have
\begin{equation}
\E_{\vec x}^{c,g}[(X_1-x_1)^2] = \frac{1}{c}(F_cg)_1(\vec x), \qquad g\in \cH_{0^+},
\ \vec x\in [0,\infty)^2.
\end{equation}
Therefore
\begin{equation}
\Gamma_{\vec x}^{c,g}\{\vec y\in [0,\infty)^2\colon\, y_1<x_1/2\} \leq \Gamma_{\vec x}^{c,g}
\{\vec y\in [0,\infty)^2\colon\, |y_1-x_1|\geq x_1/2\} \leq \frac{4(F_cg)_1(\vec x)}{c x_1^2}.
\end{equation}
We claim that
\be{4Fcg}
\lim_{x_1\to\infty} \sup_{g\in {\cal A} \atop x_2\in[\alpha -\epsilon, \alpha+\epsilon]} \frac{4(F_cg)_1(\vec x)}{cx_1^2} = 0.
\ee
Assume (\ref{4Fcg}) for the moment. Since $\beta_{\alpha,\epsilon, N}$
is nondecreasing in $N$, we can choose $N>N_0$ sufficiently large such that
\begin{equation}
\label{infrel}
\inf_{g\in {\cal A} \atop \vec x \in [N,\infty)\times [\alpha-\epsilon,\alpha+\epsilon]}
\Gamma_{\vec x}^{c,g}\Big\{\vec y\in [0,\infty)^2
\backslash [N,\infty)\times [\alpha-\epsilon,\alpha+\epsilon]\colon\, y_1\geq
\frac{x_1}{2}\Big\} \geq \frac{\beta_{\alpha,\epsilon,N_0}}{2}.
\end{equation}
Then
\be{estest}
\begin{aligned}
&\inf_{g\in {\cal A} \atop \vec x\in [N,\infty)\times [\alpha-\epsilon,\alpha+\epsilon]}
\hat\Gamma_{\vec x,h}^{c,g}\left([0,\infty)^2 \backslash [N,\infty)\times
[\alpha-\epsilon,\alpha+\epsilon] \right)  \\
&\qquad = \inf_{g\in{\cal A} \atop \vec x\in [N,\infty)\times [\alpha-\epsilon,\alpha+\epsilon]}
\int_{[0,\infty)^2\backslash [N,\infty)\times [\alpha-\epsilon,\alpha+\epsilon]}
\frac{h(\vec y)}{h(\vec x)} \Gamma_{\vec x}^{c,g}(d\vec y) \\
&\qquad \geq \inf_{g\in{\cal A} \atop \vec x\in [N,\infty)\times [\alpha-\epsilon,\alpha+\epsilon]}
\int_{\left\{y_1\geq x_1/2, \atop [0,\infty)^2\backslash [N,\infty)\times [\alpha-\epsilon,\alpha+\epsilon]\right\}}
\frac{(1+y_1)(1+y_2)}{(1+x_1)(1+x_2)} \Gamma_{\vec x}^{c,g}(d\vec y) \\
&\qquad \geq \inf_{g\in{\cal A} \atop \vec x\in [N,\infty)\times [\alpha-\epsilon,\alpha+\epsilon]}
\frac{1+x_1/2}{(1+x_1)(1+\alpha+\epsilon)} \int_{\left\{y_1\geq x_1/2, \atop [0,\infty)^2\backslash
[N,\infty)\times [\alpha-\epsilon,\alpha+\epsilon]\right\}} \Gamma_{\vec x}^{c,g}(d\vec y) \\
&\qquad \geq \frac{\beta_{\alpha,\epsilon,N_0}}{4(1+\alpha+\epsilon)} > 0,
\end{aligned}
\ee
which establishes (\ref{kb2}).

To verify (\ref{4Fcg}), note that, by condition (\ref{kbcond3}) and Proposition \ref{prop:momenteqn},
\begin{equation}
\label{g12est}
\begin{aligned}
\E_{\vec x}^{c,g}[g_1+g_2] &\leq \E_{\vec x}^{c,g}\Big[C_a(1+X_1)(1+X_2)+a(X_1^2+X_2^2)\Big] \\
& = C_a(1+x_1)(1+x_2) + a(x_1^2+x_2^2) + \frac{a}{c}\E_{\vec x}^{c,g}[g_1+g_2]
\qquad \forall\,g\in {\cal A}.
\end{aligned}
\end{equation}
Solving for $\E_{\vec x}^{c,g}[g_1+g_2]$, we get
$$
\E_{\vec x}^{c,g}[g_1+g_2] = (F_cg)_1(\vec x) + (F_cg)_2(\vec x) \leq \frac{c}{c-a}
\Big(C_a(1+x_1)(1+x_2)+a(x_1^2+x_2^2)\Big) \qquad \forall\,g\in{\cal A}.
$$
Therefore
\begin{equation}
\label{x1lim}
\limsup_{x_1\to\infty} \sup_{g\in {\cal A} \atop x_2\in [\alpha-\epsilon, \alpha+\epsilon]}
\frac{4(F_cg)_1(\vec x)}{cx_1^2} \leq \frac{4ca}{c-a}.
\end{equation}
Since $a>0$ can be made arbitrarily small, (\ref{4Fcg}) follows.
\epr

\noindent
{\bf Proof of Proposition \ref{prop:hlimitchain}:}
By Lemma \ref{lem:itermomeqn}, $h_1(\vec x) = 1+x_1$, $h_2(\vec x) = 1+x_2$ and $h(\vec x)
= (1+x_1)(1+x_2)$ are harmonic for $\vec X$. Hence, by Lemma \ref{lem:htransform},
$h_1(\vec x)/h(\vec x)=1/(1+x_2)$ and $h_2(\vec x)/h(\vec x)=1/(1+x_1)$ are harmonic for
$\vec X^h$. Therefore $(1/(1+X^h_1(n)))_{n\in\N_0}$ and $(1/(1+X^h_2(n)))_{n\in\N_0}$ are
nonnegative martingales and, by the martingale convergence theorem, $\vec X^h(n)\to\vec X^h(\infty)
\in [0,\infty]^2$ a.s.\ as $n\to\infty$. We need to show that
\begin{itemize}
\item[(i)] $\P \Big[\vec X^h(\infty)\in [0,\infty)^2,
\vec X^h(\infty)\notin \partial g\Big]=0$.
\item[(ii)] $\P\Big[X^h_1(\infty) = \infty, X^h_2(\infty)\in (0,\infty)\Big]
= \P\Big[X^h_2(\infty) = \infty, X^h_1(\infty)\in (0,\infty)\Big] = 0$.
\end{itemize}

If (i) is false, then there exists a $\vec\theta \in[0,\infty)^2\backslash\partial g$ such that,
for all $B\subset[0,\infty)^2$ with $\vec\theta\in {\rm int}(B)$, $\P[\vec X^h(n) \in B
\mbox{ for all } n \mbox{ large enough}]>0$. In particular, we must have
\begin{equation}
\label{infexitprob}
\inf_{\vec x\in B} \frac{1}{h(\vec x)}\int_{[0,\infty)^2\backslash B}
h(\vec y)\Gamma_{\vec x}^{c,g}(d \vec y) = 0 \qquad
\forall\, B\subset [0,\infty)^2 \mbox{ with } \vec\theta\in\mbox{int}(B).
\end{equation}
Otherwise, there is a uniform probability of escaping from $B$ at each step, and $\vec X^h$ cannot
be confined in $B$ forever with positive probability.

If (ii) is false, then (considering without loss of generality the first part of (ii)) there exists
an $\alpha\in (0,\infty)$ such that
\begin{equation}
\P\Big[\vec X^h(n) \in [N,\infty)\times [\alpha-\epsilon, \alpha+\epsilon]
\mbox{ for all $n$ large enough} \Big] > 0 \qquad \forall\, \epsilon\in (0,\alpha),\,N>0.
\end{equation}
In particular, we must have
\begin{equation}
\label{infexitprob2}
\inf_{\vec x \in [N,\infty)\times [\alpha-\epsilon,\alpha+\epsilon]}
\frac{1}{h(\vec x)}\int_{[0,\infty)^2\backslash [N,\infty)\times [\alpha-\epsilon,\alpha+\epsilon]}
h(\vec y)\Gamma_{\vec x}^{c,g}(d \vec y) = 0 \qquad \forall\,\epsilon>0,\, N>0.
\end{equation}
But both (\ref{infexitprob}) and (\ref{infexitprob2}) contradict Lemma \ref{lem:kerest} applied to
${\cal A} = \{ g\}$, where conditions (\ref{kbcond1}--\ref{kbcond3}) in Lemma \ref{lem:kerest}
are easily verified by our assumption that $g\in \cH_{0^+}$ and that $g$ satisfies (\ref{gcon1}).
Therefore we must have $\lim_{n\to\infty} \vec X^h(n) =\vec X^h(\infty) \in R$ a.s.
\qed

\subsection{Proof of Proposition \ref{prop:chainharm} and Theorems \ref{thm:upfs}
and \ref{thm:downfs}}
\label{S4.3}

{\bf Proof of Proposition \ref{prop:chainharm}:} Let $f$ be a nonnegative harmonic
function of $\vec X=\vec M^{c,g}$ satisfying the constraints in (\ref{constrain1}--\ref{constrain3}).
Since $x_1$, $x_2$ and $x_1x_2$ are harmonic for $\vec X$, so is $f_0(\vec x) =
f(\vec x) - \lambda_{f,(0,\infty)}x_2-\lambda_{f,(\infty,0)}x_1-\lambda_{f,(\infty,\infty)}x_1x_2$.
Let $\vec X^h$ denote the $h$-transform of $\vec X$ with $h(\vec x) = (1+x_1)(1+x_2)$.
Then, by Lemma \ref{lem:htransform}, $f_0/h$ is harmonic for $\vec X^h$, and so
\be{f0hratio}
\frac{f_0(\vec \theta)}{h(\vec\theta)}
= \E\Big[\frac{f_0(\vec X^h(n))}{h(\vec X^h(n))} \,\Big|\,\vec X^h(0) = \vec\theta\Big]
\qquad \forall \ n\in\N,\,\vec\theta \in [0,\infty)^2.
\ee
Constraint (\ref{constrain1}) implies that $f_0/h$ is bounded, constraint (\ref{constrain3})
implies that $\lim_{\vec x\to \vec z} f_0(\vec x)/h(\vec x)$ $= 0$ for all $\vec z\in R_\infty$,
while constraints (\ref{constrain2}--\ref{constrain3}) imply that $\lim_{\vec x\to \vec z}
f_0(\vec x)/h(\vec x) = 0$ for all $\vec z\in \partial g$. Since, by Proposition \ref{prop:hlimitchain},
$\lim_{n\to\infty} \vec X^h(n)= \vec X^h(\infty) \in R\, (=\partial g \cup R_\infty)$ a.s., letting
$n\to\infty$ in (\ref{f0hratio}) and applying the bounded convergence theorem, we obtain
$f_0/h \equiv 0$ and $f_0\equiv 0$. Therefore $f(\vec x) = \lambda_{f,(0,\infty)}x_2
+\lambda_{f,(\infty,0)} x_1+ \lambda_{f,(\infty,\infty)}x_1x_2$.
\qed

\noindent
{\bf Proof of Theorem \ref{thm:upfs}:} Suppose the claim is false. Then, without loss of
generality, we may assume that $((F_cg)_1, (F_cg)_2) = (\lambda_1 g_1, \lambda_2 g_2)$
for some $g\in \cH_{0^+}$, $\lambda_1>1$, $\lambda_1 \geq \lambda_2 >0$. By Definition
\ref{def:Ha}, for any $a>0$ there exists a $0<C_a<\infty$ such that $g_1(\vec x)+g_2(\vec x)
\leq C_a(1+x_1)(1+x_2)+a(x_1^2+x_2^2)$. Fix $\vec X(0)=\vec\theta\in [0,\infty)^2$, then by Lemma \ref{lem:itermomeqn},
we have \be{Eg1exp}
\begin{aligned}
\lambda_1^n g_1(\vec\theta) &\quad =\quad \E[g_1(\vec X(n))] \\
&\quad \leq\quad \E\Big[C_a(1+X_1(n))(1+X_2(n))+a(X_1^2(n)+X_2^2(n))\Big] \\
&\quad \leq\quad C_a(1+\theta_1)(1+\theta_2) + a(\theta_1^2+\theta_2^2)
+ \frac{a}{c}\sum_{j=1}^n
\Big(\lambda_1^j g_1(\vec\theta)+\lambda_2^j g_2(\vec\theta)\Big).
\end{aligned}
\ee
Since $\lambda_1>1$ and $\lambda_1\geq \lambda_2 >0$, dividing both sides of the
above inequality by $\lambda_1^n$ and letting $n\to\infty$, we get
\be{g1up}
g_1(\vec\theta) \leq \frac{a\lambda_1}{c(\lambda_1 -1)}\left[g_1(\vec\theta)
+ 1_{\lambda_1=\lambda_2}\,g_2(\vec\theta)\right].
\ee
Since $a>0$ can be made arbitrarily small, (\ref{g1up}) implies that $g_1(\vec\theta)
\leq 0$, which is a contradiction.
\qed

\bigskip\noindent
{\bf Proof of Theorem \ref{thm:downfs}:}
(i) Assume that, for some $g\in \cH_{0^+}$ with $\liminf_{\vec x\to (\infty,\infty)}
[g_1(\vec x)/x_1^2+g_2(\vec x))/x_2^2]=0$, $F_c(g_1, g_2) = (\lambda_1 g_1, \lambda_2 g_2)$
for some $0<\lambda_1,\lambda_2<1$. Fix $\vec X(0)=\vec\theta \in [0,\infty)^2$. By Lemma \ref{lem:itermomeqn}, we have
\be{sumiter1}
\E[(X_i(n)-\theta_i)^2] = \frac{1}{c}\sum_{k=1}^n \lambda_i^k
g_i(\vec\theta) < \frac{\lambda_i}{c(1-\lambda_i)} g_i(\vec\theta), \qquad \forall\,n\in \N.
\ee
Next, choose $\vec\theta$ such that $\frac{\lambda_i g_i(\vec\theta)}{c(1-\lambda_i)} \leq
\frac{\theta_i^2}{16}$ for $i=1,2$, which is possible by the above assumptions. Then, by the Chebychev inequality,
\be{sumiter2}
\P\left(\vec X(n) \in \left[\frac{\theta_1}{2},\frac{3\theta_1}{2}\right]
\times\left[\frac{\theta_2}{2}, \frac{3\theta_2}{2}\right] \
\Big|\ \vec X(0)=\vec\theta\right) \geq \frac{1}{2} \qquad \forall\, n\in\N,
\ee
and hence
\be{sumiter3}
\E\left[g_i(\vec X(n))\right]
\geq \frac{1}{2}\ \inf_{\scriptscriptstyle \vec x\in \left[\frac{\theta_1}{2},
\frac{3\theta_1}{2}\right]\times\left[\frac{\theta_2}{2}, \frac{3\theta_2}{2}\right]}
g_i(\vec x) >0 \qquad \forall\, n\in\N,
\ee
which contradicts the assumption that $\E[g_i(\vec X(n))] = \lambda_i^n
g_i(\vec\theta)\to 0$ as $n\to\infty$.
\medskip

\noindent
(ii) We consider the conditions (\ref{fscond1}) and (\ref{fscond2}) separately.

\medskip\noindent
\underline{(\ref{fscond1})}:
Assume that $(F_cg)_1 = \lambda_1 g_1$ with $\lambda_1<1$ and $g_1(x_1,0)>0$ for all $x_1>0$
for some $g\in\cH_{0+}$. For $\vec\theta = (\theta_1,0)$ with $\theta_1\geq 0$,
$\vec\Gamma_{\vec\theta}^{c,g}(d\vec x)$ is supported on the horizontal axis $A_1$
and is in fact the equilibrium distribution of the one-dimensional diffusion
\beq
dX_1(t) = c(\theta_1 -X_1)dt + \sqrt{2g_1(X_1, 0)} dB_1(t).
\eeq
Therefore the mapping $g_1(x_1, 0) \mapsto (F_cg)_1(x_1, 0)$ is the renormalization transformation
for diffusions on the halfline which, by Lemma 2 and Theorem 2 in Baillon, Cl\'ement, Greven and
den Hollander~\cite{BCGH97}, cannot have a fixed shape with scaling constant $\lambda_1\neq 1$.

\medskip\noindent
\underline{(\ref{fscond2})}:
Assume that $(F_cg)_1 = \lambda_1 g_1$ with $\lambda_1 \in (0,1)$ for some $g\in \cH_{0^+}$
such that $\liminf_{\vec x\to (\infty, \infty)}$ $g_1(\vec x)/x_1x_2 = \epsilon >0$. Then the
$h$-transformed Markov chain $\vec X^h$ with $h(\vec x)= (1+x_1)(1+x_2)$ satisfies
$\E[(g_1/h)(\vec X^h(n))]=\lambda_1^n(g_1/h)(\vec X^h(0))$. If $X^h_1(0),X^h_2(0) >0$, then, by
Corollary \ref{cor:trapprob},
\be{Prelnz}
\P[\vec X^h(\infty) = (\infty,\infty)] = \frac{X^h_1(0)X^h_2(0)}{(1+X^h_1(0))(1+X^h_2(0))}>0
\ee
and
\begin{eqnarray}
0 &=& \lim_{n\to\infty}\lambda_1^n\frac{g_1(\vec X^h(0))}{h(\vec X^h(0))}
= \lim_{n\to\infty} \E\left[\frac{g_1(\vec X^h(n))}{h(\vec X^h(n))}\right] \\
&\geq& \frac{X^h_1(0)X^h_2(0)}{(1+X^h_1(0))(1+X^h_2(0))} \liminf_{\vec x\to(\infty,\infty)}
\frac{g_1(\vec x)}{h(\vec x)}>0,
\end{eqnarray}
which is a contradiction. \qed

\section{Proof of Theorem \ref{thm:DA} with constant $c_n$}
\label{S5}

\bpr
Assume $c_n\equiv c>0$, in which case $F^{[n]} = F_c^n$. The proof is based on an analysis of the {\it interaction chain} introduced in Section
\ref{S1.2}. Let $g$ satisfy the conditions in Theorem \ref{thm:DA}. Let
$\vec X=(\vec X(-n))_{n\in\N_0}$ be the (inhomogeneous) \emph{backward Markov chain} on
$[0,\infty)^2$ with transition probability kernel
\begin{equation}
\P\left(\vec X(-n) \in d\vec x \,\Big|\, \vec X(-n-1)=\vec \theta\right)
= \Gamma_{\vec\theta}^{c, F_c^ng}(d\vec x).
\end{equation}
Denote the transition probability kernel from time $-m$ to time $-n>-m$ by
$K^{-m, -n}(\vec x, d\vec u)$. By Proposition \ref{prop:momenteqn}, the functions
$1$, $x_1$, $x_2$ and $x_1x_2$ are harmonic for $\vec X$. Let $\vec X^h
=(\vec X^h(-n))_{n\in\N_0}$ denote the $h$-transform of $\vec X$ with $h(\vec x)=(1+x_1)(1+x_2)$.
Then $1$, $\frac{x_1}{1+x_1}$, $\frac{x_2}{1+x_2}$ and $\frac{x_1 x_2}{h(\vec x)}$ are harmonic
for $\vec X^h$. Now change variables and let
\be{}
\vec Y(-n) = \phi(\vec X^h(-n)),
\ee
with $\phi\colon\, [0,\infty)^2\to [0,1)^2$ given by
\be{phidefcomp}
\phi(x_1, x_2)=\left(\frac{x_1}{1+x_1},\frac{x_2}{1+x_2}\right).
\ee
Then $\vec Y=(\vec Y(-n))_{n\in\N_0}$ is a backward Markov chain on $[0,1)^2$ with $1$, $y_1$,
$y_2$ and $y_1 y_2$ harmonic. Denote its transition probability kernel from time $-m$ to time
$-n>-m$ by $\hat K^{-m,-n}(\vec y, d\vec v)$. Then $\hat K^{-m,-n}$ and $K^{-m,-n}$ are related
via
$$
\int_{[0,\infty)^2}\! f(\vec x)\, K^{-m, -n}(\vec\theta, d\vec x)
= h(\vec\theta)\int_{[0,1)^2}\! \Big(\frac{f}{h}\circ\phi^{-1}\Big)(\vec y)\,
\hat K^{-m,-n}(\phi(\vec\theta), d\vec y) \qquad \forall\,f\ {\rm measurable}.
$$
In particular,
\begin{eqnarray}
\label{fcjgeqn}
(F_c^jg)_i(\vec\theta) &=& \int_{[0,\infty)^2}
g_i(\vec x)\, K^{-j,0}(\vec\theta, d\vec x) \\
&=& h(\vec\theta) \int_{[0,1)^2} \Big(\frac{(F_c^N g)_i}{h}\circ\phi^{-1}\Big)(\vec y)
\,\hat K^{-j,-N}(\phi(\vec\theta),
d\vec y), \qquad 0\leq N \leq j,\ i=1,2,
\nonumber
\end{eqnarray}
since $(F_c^jg)_i(\vec\theta) = \E[(F_c^Ng)_i(\vec X(-N))|\vec X(-j)=\vec\theta]$
for all $0\leq N \leq j$. For $j \in \N$, if we let
\be{}
\vec Y^{(j)}=(\vec Y^{(j)}(-n))_{n\in\N_0}
\ee
denote the Markov chain $\vec Y$ started at time $-j$ with $\vec Y^{(j)}(-j) = \phi(\vec\theta)$, and for all $-n<-j$
set $\vec Y^{(j)}(-n) = \phi(\vec\theta)$, then we can rewrite (\ref{fcjgeqn}) as
\begin{equation}
\label{fcjgeqn2}
(F_c^jg)_i(\vec\theta) = h(\vec\theta) \E\left[\left(\frac{g_i}{h}\circ\phi^{-1}\right)
\left(\vec Y^{(j)}(0)\right)\right]
= h(\vec\theta) \E\left[\left(\frac{(F_c^N g)_i}{h}\circ\phi^{-1}\right)
\left(\vec Y^{(j)}(-N)\right)\right].
\end{equation}

To establish (\ref{DAFnlim}), and hence Theorem \ref{thm:DA} for $c_n\equiv c$, we need the following lemma,
the proof of which is postponed.

\bl{lem:supplim}
For any fixed $N\in\N_0$, all weak limit points of $\{\vec Y^{(j)}(-N)\}_{j\in\N}$ as $j\to\infty$ are supported
on $\phi(R_\infty) \cup ([0,1)\times\{0\}) \cup (\{0\}\times [0,1))$.
\el

We first complete the proof subject to Lemma \ref{lem:supplim}. Without loss of generality,
take $i=1$. Note that, since $g\in \cH_0^r$, we have $g_1(\vec x)+g_2(\vec x)\leq C(1+x_1)(1+x_2)$ for some $C>0$.
Consequently, by the moment equations (\ref{eq:1stmoment})--(\ref{eq:mixmoment}), the family of functions
\be{seq1}
\left\{\left(\frac{(F_c^kg)_1}{h}\circ \phi^{-1}\right)(\vec y)\right\}_{k\in\N_0,\vec y\in [0,1)^2}
\ee
is uniformly bounded. Now fix $\vec\theta \in [0,\infty)^2$. If $\{j'_m\}_{m\in\N}$ is any subsequence along which $\lim_{m\to \infty}
(F_c^{j'_m}g)_1(\vec\theta)$ exists, then we can find a further subsequence $\{j_m\}_{m\in\N}$
such that $\vec Y^{(j_m)}$ converges weakly to a limit $\vec Y^\infty = (\vec Y^\infty(-n))_{n\in\N_0}$
as $([0,1]^2)^\N$-valued random variables with the product topology. In particular,
$\vec Y^{(j_m)}(-N)$ converges weakly to $\vec Y^\infty(-N)$ for each $N\in\N_0$.

By Theorem \ref{thm:FCdomain}, the family
\be{seq2}
\left\{\left(\frac{(F_c^kg)_1}{h}\circ \phi^{-1}\right)(\vec y)\right\}_{k\in\N_0}
\ee
is continuous on $[0,1)^2$. In fact, it is also continuous at $\phi(R_\infty)$ with
\be{seq3}
\left(\frac{(F_c^kg)_1}{h}\circ \phi^{-1}\right)(\vec z)=\lambda_{1,\vec z}
\qquad  \forall\, k\in\N_0,\,\vec z \in \phi(R_\infty).
\ee
Indeed, this follows from these observations: (1) $g\in \cH_0^r$, and hence $((g_1/h)\circ\phi^{-1})
(\vec z) = \lambda_{1,\vec z}$ for $\vec z\in \phi(R_\infty)$ and is continuous at $\vec z$;
(2) by (\ref{fcjgeqn2}), $(F_c^kg)_1(\vec\theta)/h(\vec\theta)=\E[(g_1/h)\circ \phi^{-1})(\vec Y^{(k)}(0))]$;
(3) because $Y^{(k)}_i$, $i=1,2$, are martingales while $\phi(R_\infty)= \{(1,0), (0,1), (1,1)\}$
are extremal in $[0,1]^2$, it follows from the Markov inequality that $\hat K^{-k,0}(\phi(\vec\theta), d\vec y)$
converges weakly to the point mass at $\vec z$ as $\phi(\vec\theta)\to \vec z$ for $\vec z \in \phi(R_\infty)$.
By Lemma \ref{lem:supplim}, we can now substitute $j_m$ for $j$ in (\ref{fcjgeqn2}) and take the limit $m\to\infty$, to obtain
\begin{equation}
\lim_{m\to\infty} (F_c^{j_m}g)_1(\vec\theta) = h(\vec\theta)\,
\E\left[\left(\frac{(F_c^N g)_1}{h}\circ\phi^{-1}\right)\left(\vec Y^{\infty}(-N)\right)\right]
\qquad \forall\, N\in\N_0.
\end{equation}
Denote the distribution of $\vec Y^\infty(-N)$ by $\mu_N$. Again by Lemma \ref{lem:supplim}, $\mu_N$
is concentrated on $\phi(R_\infty) \cup [0,1)\times \{0\} \cup \{0\}\times [0,1)$. Consequently,
because $((F_c^Ng)_1/h\circ \phi^{-1})(\vec y)$ vanishes on $\{0\}\times [0,1]$, we have
\begin{equation}
\label{limfcjmg}
\begin{aligned}
&\lim_{m\to\infty} (F_c^{j_m}g)_1(\vec\theta)\\
&= h(\vec\theta)\left(\mu_N\{(1,1)\}\lambda_{1,(\infty,\infty)} + \int_0^1
\left(\frac{(F_c^N g)_1}{h}\circ\phi^{-1}\right)(y_1,0)\, \mu_N(dy_1\times\{0\})\right).
\end{aligned}
\end{equation}
Since $y_1$, $y_2$, $y_1y_2$ are bounded continuous functions on $[0,1]^2$ and since
$\E[Y^{(j_m)}_i(-N)] = \phi_i(\vec\theta)$ and $\E[Y_1^{(j_m)}(-N)Y_2^{(j_m)}(-N)] =
\phi_1(\vec\theta)\phi_2(\vec\theta)$ with $\phi = (\phi_1,\phi_2)$, we must
also have $\int y_i \mu_N(d\vec y)$ $= \phi_i(\vec\theta)$ and $\int y_1 y_2 \mu_N(d\vec y)
= \phi_1(\vec\theta) \phi_2(\vec\theta)$. By our property of the support of $\mu_N$, we
thus find
\begin{eqnarray}
\mu_N\{(1,1)\} &=& \phi_1(\vec\theta)\phi_2(\vec\theta)
= \frac{\theta_1\theta_2}{h(\vec\theta)}, \\
\int y_1\ \mu_N(d y_1\times \{0\}) &=& \int y_1 (1-y_2) \mu_N(d\vec y)
= \phi_1(\vec\theta)(1-\phi_2(\vec\theta))
= \frac{\theta_1}{h(\vec\theta)}.
\end{eqnarray}
Therefore
\begin{eqnarray}
&&
\!\!\!\!\!\!\!\!\!\!\!\!\!\!\!\!\!\!\!\!\!\!\!\!\!\!\!\!\!\!
\Big|\lim_{m\to\infty} (F_c^{j_m}g)_1(\vec\theta) - \lambda_{1,(\infty,\infty)}\theta_1\theta_2
- \lambda_{1,(\infty,0)}\theta_1\Big| \nonumber \\
&\leq& h(\vec\theta)  \sup_{y_1\in [0,1]}
\Big|\big(\frac{(F_c^N g)_1}{h}\circ\phi^{-1}\big)(y_1,0) - \lambda_{1,(\infty,0)}y_1\Big| \nonumber\\
&=& h(\vec\theta) \sup_{x>0} \Big|\frac{(F_c^Ng)_1(x,0)-\lambda_{1,(\infty,0)}x}{1+x}\Big|. \label{subseqlimconv}
\end{eqnarray}
Next, note that $((F_c^Ng)_1)_{N\in\N_0}$ restricted to $(0,\infty)\times \{0\}$ are the iterates
of the renormalization transformation acting on diffusion functions on the halfline with initial
diffusion function $g_1(x,0)$. Since $\lim_{x\to\infty} g_1(x,0)/x = \lambda_{1,(\infty,0)}
\in [0,\infty)$, Theorem 5 of Baillon, Cl\'ement, Greven and den Hollander~\cite{BCGH97} implies
that $\sup_{x>0} |(F_c^Ng)_1(x,0)-\lambda_{1,(\infty,0)}x|/(1+x)\to 0$ as $N\to\infty$. (The
case $\lambda_{1,(\infty,0)}=0$ is not included in Theorem 5 in~\cite{BCGH97}, but an examination
of the proof shows that the same result holds.) Since $N$ can be taken arbitrarily large in
(\ref{subseqlimconv}), we have established the convergence in (\ref{DAFnlim}) along the subsequence
$\{j_m\}_{m\in\N}$. Since $\{(F_c^jg)_1(\vec\theta)\}_{j\in\N_0}$ is uniformly bounded,
(\ref{DAFnlim}) now follows and the proof of Theorem \ref{thm:DA} for $c_n\equiv c$ is complete.
\epr

We now prove  Lemma \ref{lem:supplim}.

\noindent
{\bf Proof of Lemma \ref{lem:supplim}.}
We must prove that the weak limit of $\{\vec Y^{(j_m)}\}_{m\in\N}$, written $\vec Y^\infty$,
satisfies
\be{seq4}
\P\big(\vec Y^\infty(-N) \in \phi(R_\infty) \cup [0,1)\times\{0\}\cup \{0\}
\times [0,1)\big)=1\qquad \forall\ N\in\N_0.
\ee
The proof consists of the following three steps:
\begin{itemize}
\item[(A)] Show that $(Y^\infty_i(-n))_{n\in\N_0}$, $i=1,2$, are backward martingales
on $[0,1]$, i.e.,
\be{seq5}
\E\Big[Y^\infty_i (-k)\, \Big|\, (Y^\infty_i(-n))_{n\geq k+1}\Big]
= Y^\infty_i(-k-1), \qquad i=1, 2,
\ee
implying that $\lim_{n\to\infty} \vec Y^\infty(-n) = \vec Y^\infty(-\infty)$ exists a.s.\ by
the backward martingale convergence theorem (see e.g. Section 4.6 in Durrett~\cite{Du96}).

\item[(B)] Show that $\P\big\{\vec Y^\infty(-\infty) \in
\phi(R_\infty) \cup [0,1)\times\{0\}\cup \{0\}\times [0,1) \big\}=1$.

\item[(C)] Show that $\P\big\{\vec Y^\infty(-N) \in
\phi(R_\infty)\cup [0,1)\times\{0\}\cup \{0\}\times [0,1)\big\}=1$ for all $N\in\N_0$.
\end{itemize}

\noindent
Since $(Y^{(j_m)}_i(-n))_{n\in\N_0}$, $m\in\N, i=1,2$, are bounded backward martingale sequences,
(A) follows from a general result on weak limits of backward martingale sequences, which
we state as Lemma \ref{lem:bmartlim} below. The proof of (B) given below uses Lemma \ref{lem:kerest},
which relies on uniform lower and upper bounds on $\{F_c^ng\}_{n\in\N_0}$, where
assumptions (\ref{DAgcon1}) and $g\in \cH_0$ are crucial. The proof of (C) given below is
achieved after approximating $\vec Y^\infty$ by the Markov chains $\vec Y^{(j_m)}$ and using
the fact that $Y^{(j_m)}_i$, $i=1,2$, are martingales. Note that it is not clear if $\vec Y^\infty$ is
a Markov chain, because $\vec Y^{(j_m)}$ take values in
$[0,1)^2$ while $\vec Y^\infty$ takes values in $[0,1]^2$.
Even though the transition kernels of $\vec Y^{(j_m)}$ are
consistent for $m$ sufficiently large, they may not be (weakly)
continuously extendable to $[0,1]^2\backslash [0,1)^2$.

\bl{lem:bmartlim} {\bf [Weak limits of backward martingales]} \\
For $j\in \N$, let $Z^{(j)}=(Z^{(j)}(-n))_{n\in\N_0}$ be a backward martingale, i.e.,
\be{bamg}
\E\Big[Z^{(j)} (-k)\, \Big|\, (Z^{(j)}(-n))_{n\geq k+1}\Big] = Z^{(j)}(-k-1).
\ee
If $\{Z^{(j)}(0)\}_{j\in \N}$ are uniformly integrable, and $Z^{(j)}$ converges weakly to a
random variable $Z^\infty= (Z^\infty(-n))_{n\in\N_0}$ in the space $\R^\N$ with the product
topology, then $(Z^\infty(-n))_{n\in\N_0}$ is also a backward martingale.
\el

\bpr
Since $\{Z^{(j)}(0)\}_{j\in \N_0}$ are uniformly integrable, we have
\begin{equation}
\forall\,\epsilon>0,\,\exists\ N>0 \mbox{ such that }
\quad \E\Big[|Z^{(j)}(0)|\, 1_{|Z^{(j)}(0)|\geq N}\Big] \leq \epsilon
\qquad \forall\,j\in\N,
\end{equation}
which is easily seen to be equivalent to
\begin{equation}
\forall\,\epsilon>0,\,\exists\ N>0 \mbox{ such that }
\quad \E\Big[\big(|Z^{(j)}(0)|-N\big)^+ \Big] \leq \epsilon
\qquad \forall\,j\in\N.
\end{equation}
Since $f(x)=(|x|-N)^+$ is a convex function, for all $j, k\in \N$ we have, by Jensen's inequality,
\begin{eqnarray}
\E\Big[\big(|Z^{(j)}(-k)|-N\big)^+ \Big] &=& \E\Big[f\left(Z^{(j)}(-k)\right)\Big]
\nonumber \\
&=& \E\Big[f\left( \E\big[Z^{(j)}(0)\,\big|\,(Z^{(j)}(-n))_{n\geq k}\big]\right)\Big]
\nonumber \\
&\leq& \E\Big[\E\big[f\big(Z^{(j)}(0)\big)\,\big|\,(Z^{(j)}(-n))_{n\geq k}\big]\Big]
\nonumber\\
&=&\E\Big[f\big(Z^{(j)}(0)\big)\Big] = \E\Big[\big(|Z^{(j)}(0)|-N\big)^+ \Big].
\end{eqnarray}
Therefore $\{Z^{(j)}(-n)\}_{j\in\N,n\in\N_0}$ is a uniformly integrable family.

For each $k\in\N_0$ and $j\in\N$, and any bounded continuous function $f\colon\,
\R^\N \to \R$, the martingale property of $Z^{(j)}$ implies that
\begin{equation}
\E\Big[ f\left(\big(Z^{(j)}(-n)\big)_{n\geq k+1}\right)\,
\left(Z^{(j)}(-k)-Z^{(j)}(-k-1)\right)\Big] = 0.
\end{equation}
Since $Z^{(j)}$ converges weakly to $Z^{\infty}$, and $\{Z^{j}(-k)\}_{j\in\N}$ and
$\{Z^{j}(-k-1)\}_{j\in\N}$ are uniformly integrable, we may pass to the limit $j\to\infty$
and obtain
\begin{equation}
\label{martlimeqn}
\E\Big[ f\left(\big(Z^{\infty}(-n)\big)_{n\geq k+1}\right)\,
\Big(Z^{\infty}(-k)-Z^{\infty}(-k-1)\Big)\Big] = 0.
\end{equation}
Indeed, the latter is easily verified by applying Skorohod's representation theorem, which allows for
a coupling between $\{Z^{(j)}\}_{j\in\N}$ and $Z^\infty$ such that the convergence is a.s.
From (\ref{martlimeqn}) we have
\begin{equation}
\label{martlimeqn*}
\E\left[ f\left(\big(Z^{\infty}(-n)\big)_{n\geq k+1}\right)\,
\E\Big[Z^{\infty}(-k)-Z^{\infty}(-k-1)\,\Big|\,
\big(Z^{\infty}(-n)\big)_{n\geq k+1}\Big]\, \right] = 0,
\end{equation}
which implies that
\begin{equation}
\E\big[Z^{\infty}(-k)-Z^{\infty}(-k-1)\,\big|\,\big(Z^{\infty}(-n)\big)_{n\geq
k+1}\big]=0 \quad a.s.,
\end{equation}
and thus establishes the martingale property for $Z^\infty$.
\epr

We are now ready to verify (B) and (C).

\medskip\noindent
\underline{(B)}: Note that
\be{phiphi}
\phi(R_\infty)\cup ([0,1)\times\{0\}) \cup (\{0\}\times [0,1))
=([0,1]\times\{0\}) \cup (\{0\}\times[0,1]) \cup (1,1).
\ee
Suppose that (B) fails. Then there exists a $\vec u\in (0,1]^2\backslash (1,1)$ in the support
of the distribution of $\vec Y^\infty(-\infty)$. In particular, for each $\epsilon>0$ there
exist $\delta(\epsilon)>0$ and $N(\epsilon)>0$ such that
\begin{equation}
\P\Big\{\vec Y^\infty(-n) \in \hat B_{\epsilon/2}(\vec u)
\quad \forall\, n \geq N(\epsilon)\Big\} > \delta(\epsilon),
\end{equation}
where $\hat B_{\epsilon/2}(\vec u)=\{\vec y\in [0,1]^2\colon \Vert \vec y-\vec u\Vert\leq \epsilon/2\}$.
Since $\vec Y^{(j_m)}$ converges weakly to $\vec Y^\infty$ as $m\to\infty$, for each $M\in \N$
we can find an $m^* = m^*(M)$ sufficiently large such that
\begin{equation}
\label{longstuckchain}
\P\left\{\vec Y^{(j_{m^*})}(-n) \in \hat B_{\epsilon}(\vec u)\cap(0,1)^2\quad
\forall\,N(\epsilon)\leq n\leq N(\epsilon)+M\right\} \geq \frac12\,\delta(\epsilon).
\end{equation}
We now derive a contradiction with Lemma \ref{lem:kerest} as follows. By assumption (\ref{DAgcon1})
and the fact that $g\in \cH_0^r$, implying $g_1(\vec x)+g_2(\vec x) \leq C(1+x_1)(1+x_2)$ for some
$0<C=C(g)<\infty$, $F_c^ng$ satisfy the same upper and lower bounds for all $n\in\N$. It is then
easy to check that in Lemma \ref{lem:kerest} with ${\cal A} = \{F_c^ng\}_{n\in\N_0}$ condition
(\ref{kbcond2}) is satisfied for all $\vec\theta \in (0,\infty)^2$, and conditions
(\ref{kbcond1}--\ref{kbcond3}) are satisfied for all $\alpha >0$ (and the analogue of
(\ref{kbcond1}) for vertical strips). Since the transition kernel $\hat K^{-n-1,-n}(\phi(\vec\theta),
d\vec y)$ is related to the biased equilibrium measure $\hat\Gamma^{c, F_c^ng}_{\vec\theta, h}(d\vec x)$
through the coordinate change $\phi$, Lemma \ref{lem:kerest}~(i) and (ii) imply that, for $\vec u\in (0,1]^2
\backslash (1,1)$ and $\epsilon>0$ sufficiently small,
\begin{equation}
\inf_{n\in\N_0 \atop \vec v \in \hat B_\epsilon(\vec u)\cap [0,1)^2}
\hat K^{-n-1,-n}\big(\vec v, [0,1)^2\backslash\hat B_\epsilon(\vec u)\big) > 0.
\end{equation}
This uniform rate of escape from $\hat B_\epsilon(\vec u)$ contradicts (\ref{longstuckchain}),
where $M$ can be chosen to be arbitrarily large while $\delta(\epsilon)>0$ remains fixed.

\medskip\noindent
\underline{(C)}: For $\epsilon>0$, let
\begin{equation}
U_\epsilon = \left\{ \vec y\in [0,1]^2 \colon\,
\inf_{\vec z \in \phi(R_\infty)\cup [0,1)\times\{0\}\cup\{0\}\times [0,1)} \Vert
\vec y-\vec z\Vert \leq \epsilon \right\}.
\end{equation}
Since $\lim_{n\to\infty}\vec Y^\infty(-n)=\vec Y^\infty(-\infty)$ a.s., we can choose
$M=M(\epsilon)$ sufficiently large such that $\P(\vec Y^\infty(-M) \in U_\epsilon)
> 1-\epsilon$. Since $\vec Y^{(j_m)}(-M) \to \vec Y^\infty(-M)$ in distribution as $m\to \infty$,
we can choose $m^*(M)$ sufficiently large such that $\P(\vec Y^{(j_m)}(-M) \in U_{2\epsilon})
> 1-2\epsilon$ for all $m\geq m^*$. By the geometry of $\phi(R_\infty)\cup ([0,1)\times\{0\})
\cup (\{0\}\times [0,1))$ and the fact that $Y^{(j_m)}_i(\cdot)$, $i=1,2$, are martingales for the
Markov chain $(\vec Y^{(j_m)}(-n))_{n\in\N_0}$, an elementary application of the Chebychev
inequality shows that, for all $m\geq m^*$ and $L>2$,
\begin{equation}
\P\left(\vec Y^{(j_m)}(0) \in U_{2L\epsilon} \right) \geq (1-2\epsilon)\Big(1-\frac{2}{L}\Big).
\end{equation}
By the weak convergence of $\vec Y^{(j_m)}(0)$ to $\vec Y^\infty(0)$ as $m\to\infty$, the same
holds for $\vec Y^\infty(0)$. Now let $\epsilon\to 0$ and $L\to \infty$ such that $\epsilon L\to 0$.
Then we find that
$$
\P\left(\vec Y^\infty(0)\in \phi(R_\infty)\cup ([0,1)\times\{0\}) \cup (\{0\}\times [0,1))\right) =1.
$$
The same argument works for $\vec Y^\infty(-N)$ for any $N\in\N_0$.
\qed


\section{Proof of Theorem \ref{thm:DA} with varying $c_n$}
\label{S6}

\proof
The proof of Theorem \ref{thm:DA} with varying $c_n$ follows the same line of argument as that for constant $c_n$,
except for a few technical differences, which we now outline. For the
rest of the section, let $(\vec X(-n))_{n\in\N_0}$ denote the backward time-inhomogeneous
Markov chain with transition kernels
\be{eq:vartransker}
\P\big(\vec X(-n) \in d\vec x \,\big|\, \vec X(-n-1)=\vec\theta\big) = \Gamma_{\vec\theta}^{c_n, F^{[n]}g}(d\vec x),
\ee
and let $(\vec X^h(-n))_{n\in\N_0}$ denote $\vec X$ $h$-transformed by $h(\vec x) = (1+x_1)(1+x_2)$, which is still
a harmonic function for $\vec X$. Both $\vec X$ and $\vec X^h$ generalize their counterparts in Section \ref{S5}.
We proceed by first establishing the analogue of Lemma
\ref{lem:supplim}, where $\{\vec Y^{(j)}\}_{j\in\N}$ are now defined
in terms of our current $\vec X$ and $\vec X^h$.

The proof of Lemma \ref{lem:supplim} in Section \ref{S5} is based on
Lemma \ref{lem:kerest}, which no longer applies in our current context, because if $c_n$ can be arbitrarily large, then
we lose the uniformity of the escape probability with respect to $\{\Gamma_{\vec x}^{c_n,F^{[n]}g}\}_{n\in\N_0}$. So, the
first task is to formulate a suitable analogue of Lemma \ref{lem:kerest} for our current $\vec X$ and
$\vec X^h$, which would imply the analogue of Lemma \ref{lem:supplim} for the present context. In the derivation of Theorem \ref{thm:DA}
for constant $c_n$ from Lemma \ref{lem:supplim}, we used the following fact from Baillon, Cl\'ement, Greven and den Hollander~\cite{BCGH97}: for the renormalization transformation $F_c$ acting on one-dimensional diffusion functions $f:[0,\infty)\to [0,\infty)$, where $f$ is positive and continuous on $(0,\infty)$, locally Lipschitz at 0, $f(0)=0$ and $\lim_{x\to\infty} f(x)/x = \lambda \in [0,\infty)$, we have $\sup_{x>0} |(F_c^nf)(x)-\lambda x|/(1+x)\to 0$ as $n\to\infty$. Our second task is therefore to establish the analogous result for $F^{[n]}f$. The two technical points outlined above will be addressed in Lemma \ref{lem:varcnkerest} and Proposition \ref{prop:1dDA} below.

Observe that, by Proposition \ref{prop:momenteqn}, for all $-m\leq -n \leq 0$ and $\vec\theta \in [0,\infty)^2$,
the backward Markov chain $\vec X$ satisfies the moment equations
\begin{eqnarray}
&& \E\big[ \vec X(-n)\, \big|\, \vec X(-m)=\vec\theta\big] = \vec \theta, \label{eq:varcn1stmom}\\
&& \E\big[ X_1(-n) X_2(-n)\, \big|\, \vec X(-m)=\vec\theta \big] = \theta_1 \theta_2, \label{eq:varcnmixmom}\\
&& \E\big[ X_i(-n)^2\, \big|\, \vec X(-m) =\vec\theta \big] = \theta_i^2 + \left(\sum_{j=n}^{m-1} \frac{1}{c_j}\right) (F^{[m]}g)_i(\vec\theta), \quad i=1, 2. \label{eq:varcn2ndmom} \\
&& \E\big[ g_i(\vec X(0)) \, \big|\, \vec X(-m) =\vec\theta \big] = (F^{[m]}g)_i(\vec\theta), \qquad \qquad \qquad \qquad i=1,2,
\label{eq:vargmom}
\end{eqnarray}
From the point of view of variance increment, (\ref{eq:varcn2ndmom}) indicates that the natural time associated with
$(\vec X(-n))_{n\geq 0}$ is not $n$, but rather $\sum_{i=0}^{n-1} c_i^{-1}$. Therefore to obtain a uniform bound on
escape probabilities for the Markov chain $\vec X^h$, we formulate the analogue of Lemma \ref{lem:kerest} as
follows.

\bl{lem:varcnkerest}{\bf [Uniform rate of escape of $(\vec X^h(-n))_{n\geq 0}$ from small balls and thin strips]}\\
Let $(c_n)_{n\in \N_0}$ and $g$ be as in Theorem $\ref{thm:DA}$. Let $(\vec X(-n))_{n\in\N_0}$ denote the inhomogeneous
backward Markov chain with transition kernel $(\ref{eq:vartransker})$, and let $(\vec X^h(-n))_{n\in\N_0}$ denote
$(\vec X(-n))_{n\in\N_0}$ $h$-transformed by $h(\vec x) = (1+x_1)(1+x_2)$. There exists an increasing sequence
$(n_k)_{k\in\N_0} \subset \N_0$ with $n_0=0$ such that $\sum_{i=n_{k}}^{n_{k+1}-1} c_i^{-1} \in [\Lambda^{-1}, \Lambda]$ for some $\Lambda >1$
for all $k\in\N_0$. For $A\subset [0,\infty)^2$, denote $\tau_A^{-m} = \inf\{-j \geq -m : \vec X^h(-j) \notin A\}$. Then
\begin{itemize}
\item[$(i)$] For each $\vec\theta \in (0,\infty)^2$, there exists $\epsilon>0$ such that
\be{eq:varkb1}
\inf_{k\in\N_0 \atop \vec x\in B_\epsilon(\vec\theta)} \P\left( \tau^{-n_{k+1}}_{B_{\epsilon}(\vec \theta)} \leq -n_k\,\Big|\, \vec X^h(-n_{k+1})=\vec x\right) >0.
\ee

\item[$(ii)$] For each $\alpha>0$, there exist $\epsilon, N>0$ such that
\beq
\inf_{k\in\N_0 \atop \vec x\in [N,\infty)\times [\alpha-\epsilon, \alpha+\epsilon]}
\P\left( \tau^{-n_{k+1}}_{[N,\infty)\times [\alpha-\epsilon, \alpha+\epsilon]} \leq -n_k\,\Big|\, \vec X^h(-n_{k+1})=\vec x\right) >0, \label{eq:varkb2}\\
\inf_{k\in\N_0 \atop \vec x\in [\alpha-\epsilon, \alpha+\epsilon]\times [N,\infty)}
\P\left( \tau^{-n_{k+1}}_{[\alpha-\epsilon, \alpha+\epsilon]\times [N,\infty)} \leq -n_k\,\Big|\, \vec X^h(-n_{k+1})=\vec x\right) >0. \label{eq:varkb2'}
\eeq
\end{itemize}
\el
\proof The existence of the increasing sequence $(n_k)_{k\in\N_0}$ with the prescribed property follows immediately from our assumptions that $\inf_{n\in\N_0} c_n>0$ and $\sum_{n\in\N_0} c_n^{-1}=\infty$. The rest of the proof parallels that of Lemma \ref{lem:kerest}. First we prove (\ref{eq:varkb1})--(\ref{eq:varkb2'}) with $\vec X^h$
replaced by $\vec X$. By (\ref{expmomeqn}), for each $m\in\N_0$ and $\vec x\in [0,\infty)^2$, $i=1,2$, conditioned on
$\vec X(-m-1) = \vec x$, we have
\beq
\E\left[\frac{1}{(1+X_i(-m))^2} \right] &=&  \frac{1}{1+x_i} \E\left[\frac{1}{1+X_i(-m)}\right] + \frac{2}{c_m(1+x_i)}\E\left[\frac{(F^{[m]}g)_i(\vec X(-m))}{(1+X_i(-m))^3}\right] \nonumber\\
&\geq& \frac{1}{(1+x_i)^2} + \frac{2}{c_m(1+x_i)}\E\left[\frac{(F^{[m]}g)_i(\vec X(-m))}{(1+X_i(-m))^3}\right], \label{eq:varratmomeqn}
\eeq
where we applied Jensen's inequality. Conditioned on $\vec X(-n_{k+1}) = \vec x$, we can apply
(\ref{eq:varratmomeqn}) iteratively to obtain, for $i=1,2$,
\beq
\label{eq:itervarratmomeqn}
\E\left[\frac{1}{(1+X_i(-n_k))^2}\right]\!\!\!\!\!\!\!\! &&\geq \\
\!\!\!\!\!\!\!\!\!\!\!\!\!\!\!\!&&\!\!\!\!\!\!\!\!\!\!\!\!\!\!\!\! \frac{1}{(1+x_i)^2} + \sum_{m=n_k}^{n_{k+1}-1} \frac{2}{c_m}
\E\left[\frac{1}{1+ X_i(-m-1)} \frac{(F^{[m]}g)_i(\vec X(-m))}{(1+X_i(-m))^3}\right]. \nonumber
\eeq
If (\ref{eq:varkb1}) fails when $\vec X^h$ is replaced by $\vec X$, then there exists $\vec \theta \in (0,\infty)^2$
such that, for all $\epsilon>0$, there exist sequences $k^{(l)}\uparrow \infty$ and $\vec x^{(l)} \in B_\epsilon(\vec\theta)$ (depending on $\epsilon$) such that
\be{eq:escapeassump}
\lim_{l\to\infty} \P\Big(\tau^{-n_{k^{(l)}+1}}_{B_\epsilon(\vec\theta)} \leq -n_{k^{(l)}}\, \Big|\, \vec X(-n_{k^{(l)}+1})=\vec x^{(l)}\Big)=0.
\ee
Now we apply (\ref{eq:itervarratmomeqn}) to $\vec X(-n)$ for $-n_{k^{(l)}+1} \leq -n \leq -n_{k^{(l)}}$ with
$\vec X(-n_{k^{(l)}+1}) = \vec x^{(l)}$. By (\ref{eq:escapeassump}), as $l\to\infty$, the two sides of (\ref{eq:itervarratmomeqn})
satisfy
\beq
\mbox{l.h.s.} &\leq& \frac{1}{(1+\theta_i -\epsilon)^2} + o(1), \label{eq:varlhs}\\
\mbox{r.h.s.} &\geq& \frac{1}{(1+\theta_i+\epsilon)^2} + \sum_{m=n_{k^{(l)}}}^{n_{k^{(l)}+1}-1}\frac{2}{c_m}\,(1-o(1))\, \frac{\delta}{(1+\theta_i+\epsilon)^4}\ ,\label{eq:varrhs}
\eeq
where in (\ref{eq:varrhs}) we have used the assumption that $g_i(\vec x) \geq \alpha_i x_i + \beta_i x_1x_2$ for some
$\alpha_i, \beta_i\geq 0$ and $\alpha_i+\beta_i>0$, $i=1,2$, which implies that $\{F^{[m]}g\}_{m\in \N_0}$ satisfy
the same lower bound and $(F^{[m]}g)_i(\vec x) \geq\delta >0$ uniformly for $\vec x\in B_\epsilon(\vec\theta)$ and $m\in\N_0$. Since $\sum_{m=n_{k}}^{n_{k-1}}\frac{2}{c_m} \geq \Lambda^{-1}>0$ uniformly
for all $k\in\N_0$, (\ref{eq:varlhs}) and (\ref{eq:varrhs}) are incompatible for $\epsilon>0$ sufficiently small
and $l\in\N$ sufficiently large. Therefore (\ref{eq:varkb1}) must hold for $\vec X$ in place of $\vec X^h$. The proof
of (\ref{eq:varkb2})--(\ref{eq:varkb2'}) for $\vec X$ in place of $\vec X^h$ is similar, and we leave the details to
the reader.

To verify that (\ref{eq:varkb1}) also holds for $\vec X^h$, we apply Lemma \ref{lem:densityhtransform} and note
that $h(\vec x)=(1+x_1)(1+x_2)$ is bounded uniformly from above for $\vec x\in B_\epsilon(\vec\theta)$, and bounded uniformly from below by 1 for $\vec x \in [0,\infty)^2$. The proof of (\ref{eq:varkb2})--(\ref{eq:varkb2'}) for $\vec X^h$ is essentially the same as its counterpart in the proof of Lemma \ref{lem:kerest}. Note that by Lemma \ref{lem:densityhtransform}, the law of $\vec X^h\big(\tau^{-n_{k+1}}_{[N,\infty)\times[\alpha-\epsilon, \alpha+\epsilon]} \wedge (-n_k)\big)$ conditioned on $\vec X^h(-n_{k+1})=\vec x\in [N,\infty)\times[\alpha-\epsilon, \alpha+\epsilon]$ is absolutely continuous with respect to the law of $\vec X\big(\tau^{-n_{k+1}}_{[N,\infty)\times[\alpha-\epsilon, \alpha+\epsilon]} \wedge (-n_k)\big)$ conditioned on $\vec X(-n_{k+1})=\vec x$, where the density is $\frac{h(\cdot)}{h(\vec x)}$. As in the proof of Lemma \ref{lem:kerest}, it suffices to show that for any fixed $0<\epsilon<\alpha<\infty$,
\be{eq:halfwayest}
\lim_{x_1\to \infty} \sup_{k\in\N_0 \atop x_2 \in [\alpha-\epsilon, \alpha+\epsilon]}
\P\left(\tau^{-n_{k+1}}_{[\frac{x_1}{2},\infty)\times\R}  \leq -n_k \,\Big|\, \vec X(-n_{k+1})= \vec x \right) = 0.
\ee
Since $(X_1(-n))_{n\leq n_{k+1}}$ is a martingale, by Doob's inequality and (\ref{eq:varcn1stmom}--\ref{eq:varcn2ndmom}), we have
\begin{eqnarray}
&& \!\!\!\!\!\!\!\!\!\!\!\!
\P\left(\tau^{-n_{k+1}}_{[\frac{x_1}{2},\infty)\times\R}  \leq -n_k \,\Big|\, \vec X(-n_{k+1})= \vec x \right) \nonumber\\
&\leq& \P\Big( \sup_{-n_{k+1} \leq -n \leq -n_k} |X_1(-n)-x_1| \geq \frac{x_1}{2}\,\Big|\, \vec X(-n_{k+1})= \vec x \Big) \nonumber\\
&\leq& \frac{16}{x_1^2}\ \E\Big[(X_1(-n_k)-x_1)^2\,\Big|\, \vec X(-n_{k+1})=\vec x\Big] \nonumber\\
&=& \frac{16}{x_1^2} \left(\sum_{n=n_k}^{n_{k+1}-1} \frac{1}{c_n}\right) (F^{[n_{k+1}]}g)_1(\vec x). \label{615}
\end{eqnarray}
Note that $\sum_{n=n_k}^{n_{k+1}-1} \frac{1}{c_n}\leq\Lambda$ uniformly in $k$. Since $g\in \cH_0^r$, we have
$g_1(\vec x)+g_2(\vec x) \leq K(1+x_1)(1+x_2)$ for some $K\in (0, \infty)$, and by Proposition \ref{prop:momenteqn},
$\{F^{[n]}g\}_{n\in\N_0}$ all share the same upper bound. Equation (\ref{eq:halfwayest}) then
follows immediately.
\qed \\
{\bf Remark.} Note that (\ref{eq:varkb1})--(\ref{eq:varkb2'}) with
$\vec X$ in place of $\vec X^h$ are proved using only the assumptions that $\sum_{n\in\N_0}c_n^{-1}=\infty$ and, $\{F^{[n]}g\}_{n\in\N_0}$
have a uniform lower bound which is positive and uniformly bounded away from $0$ on $(a,\infty)^2$ for each $a>0$. Only in
deriving (\ref{eq:varkb2})--(\ref{eq:varkb2'}) from their analogues for $\vec X$, did we use the assumptions that $\inf_{n\in\N_0} c_n>0$ and,
$\{F^{[n]}g\}_{n\in\N_0}$ have a uniform upper bound $\phi=(\phi_1, \phi_2)$, where $\phi_1(x_1, x_2)$ grows sub-quadratically
in $x_1$ and $\phi_2(x_1, x_2)$ grows sub-quadratically in $x_2$.

Using Lemma \ref{lem:varcnkerest} and the fact that $1, x_1, x_2, x_1x_2$ are still harmonic functions for the Markov chain $\{\vec X(-n)\}_{n\in\N_0}$, we
deduce the analogue of Lemma \ref{lem:supplim} in our present context by the same arguments as
in the original proof. To deduce Theorem \ref{thm:DA} with varying $c_n$ from the analogue of Lemma \ref{lem:supplim}, we need to
address the second technical point outlined at the beginning of this section.

\bp{prop:1dDA}{\bf [Convergence to fixed points under $F^{[n]}$: the half line case]} \\
Let $(c_n)_{n\in\N_0}$ satisfy $\sum_{n\in\N_0} c_n^{-1} = \infty$. Let $f(x): [0,\infty) \to [0,\infty)$ be positive and continuous on $(0,\infty)$, locally Lipschitz at 0, $f(0)=0$ and $\lim_{x\to\infty} x^{-1} f(x) = \lambda \in [0,\infty)$. Then we have
\begin{equation}
\label{1dgcon}
\lim_{n\to\infty}\ \sup_{x>0} \left|\frac{(F^{[n]}f)(x) - \lambda x}{1+x}\right| = 0,
\end{equation}
where $F^{[n]}$ are renormalization transformations acting on one-dimensional diffusion functions.
\ep
\proof
Note that we do not require $\inf_{n\in\N_0} c_n >0$ as in Lemma \ref{lem:varcnkerest}.
The case $c_n\equiv c$ is covered by Theorem 5 of Baillon, Cl\'ement, Greven and den Hollander~\cite{BCGH97}. Here we give a proof along the same
line of argument as we have been pursuing so far in this section for the proof of Theorem \ref{thm:DA} with varying
$c_n$, except that we do not need to appeal to the current proposition.

As in Section 2.3 of  Baillon, Cl\'ement, Greven and den
Hollander~\cite{BCGH97}, we make use of the concave upper envelope $f^+$ and the convex lower envelope
$f^-$ of $f$. It is easy to see that $f^+$ and, $f^-$ in the case $\lambda>0$, satisfy the same constraints as specified
for $f$ in the proposition.  Since, for any $c>0$,
$F_c$ is convexity preserving and order preserving by Proposition 3 of~\cite{BCGH97}, together with Jensen's inequality we have, for each $x\in [0,\infty)$, $(F^{[n]}f^+)(x) \downarrow f^+_\infty(x)$ and $(F^{[n]}f^-)(x) \uparrow f^-_\infty(x)$ for some $f^+_\infty$ and $f^-_\infty$ as $n\to\infty$,
and $(F^{[n]}f^-)(x) \leq (F^{[n]}f)(x) \leq (F^{[n]}f^+)(x)$ for all $n\in\N_0$. We claim that it suffices to show that $f^+_\infty(x) = f^-_\infty(x) = \lambda x$. Indeed,
$$
\sup_{x>0} \left|\frac{(F^{[n]}f)(x) - \lambda x}{1+x}\right|
= \sup_{y\in (0,1)} \Big|(1-y)\,(F^{[n]}f)\circ \phi_1^{-1}(y) -\lambda y\Big|,
$$
where $\phi_1(x) = \frac{x}{1+x}$. Since $F_c$ preserves the slope at infinity, we have that
\be{}
\begin{aligned}
\psi_n^+(y) &=(1-y)\, (F^{[n]}f^+)\circ \phi_1^{-1}(y), \\
\psi_n(y)   &=(1-y)\, (F^{[n]}f)\circ \phi_1^{-1}(y), \\
\psi_n^-(u) &= (1-y)\, (F^{[n]}f^-)\circ \phi_1^{-1}(y),
\end{aligned}
\ee
are all continuous functions on $[0,1]$. If $f^+_\infty(x) = f^-_\infty(x) = \lambda x$, then, on $[0,1]$, $\psi_n^+(y)$
decreases monotonically to $\lambda y$ as $n\to\infty$, while $\psi_n^-(y)$ increases monotonically to $\lambda y$ as $n\to\infty$. Since the monotone convergence of a sequence of continuous functions to a continuous limit is necessarily
uniform on compacts, the sup-norm convergence of $\psi_n(y)$ to $\lambda y$ on $[0,1]$ follows since $\psi_n$ is sandwiched between $\psi_n^+$ and $\psi_n^-$.

The proof that $f^+_\infty(x) = \lim_{n\to\infty} (F^{[n]}f^+)(x) = \lambda x$ and $f^-_\infty(x)=\lim_{n\to\infty} (F^{[n]}f^-)(x)=\lambda x$ now
follows the same argument as that used for Theorem~\ref{thm:DA} with varying $c_n$.  First consider the case $\{F^{[n]}f^+\}_{n\in\N_0}$ with
$\lambda>0$. In the proof of Theorem \ref{thm:DA} with varying $c_n$,
we replace $\vec X$ there by the $[0,\infty)$-valued Markov chain
$(X(-n))_{n\in\N_0}$ with transition kernels
$$
\P(X(-n) \in \cdot\,|\, X(-n-1)=x)=\Gamma_{x}^{c_n, F^{[n]}f^+}(\cdot);
$$
$\vec X^h$ is replaced by $X^h$, which is the $h$-transform of $X$ by the harmonic function $h(x) = 1+x$; $\phi(\vec x)$
is replaced by $\phi_1(x) = \frac{x}{1+x}$; in Lemma \ref{lem:supplim}, the relevant boundary points now consist of only
$\{0\} \cup \{1\}$. Lastly, because of the one-dimensional setting, we only need to establish the analogue of (\ref{eq:varkb1}).
By the remark following the proof of Lemma~\ref{lem:varcnkerest}, the only assumptions we need here are
$\sum_{n\in\N_0}c_n^{-1}=\infty$ and, a uniform lower bound on $\{F^{[n]}f^+\}_{n\in\N_0}$ which is positive and bounded
away from 0 on $[a,\infty)$ for each $a>0$. Note that $f^-$ provides such a lower bound. The case $\{F^{[n]}f^-\}_{n\in\N_0}$ with
$\lambda>0$ is identical. For the case $\lambda=0$, we only need to consider $\{F^{[n]}f^+\}_{n\in\N_0}$. Everything remains the same,
except that the uniform lower bound on $\{F^{[n]}f^+\}_{n\in\N_0}$ is now provided by $f^+_\infty$. Indeed, as a limit
of concave functions, $f^+_\infty$ is also concave, hence either $f^+_\infty\equiv 0$, in which case we are done, or
$f^+_\infty$ is positive and non-decreasing on $(0,\infty)$, which is sufficient for the proof of the analogue of
(\ref{eq:varkb1}) to go through.
\qed

With Lemma~\ref{lem:varcnkerest} and Proposition~\ref{prop:1dDA}, we can now proceed as in the proof of Theorem~\ref{thm:DA}
for constant $c_n$ and extend it to varying $c_n$. We leave the details to the reader.
\qed


\appendix

\section{Appendix 1: Moment equations and estimates}
\label{S3.1}

\bp{prop:momenteqn} {\bf [Moment equations and estimates]}\\
Let $g\in \cC$, $\vec\theta\in [0,\infty)^2$, $c>0$ and let $\Gamma_{\vec\theta}^{c,g}$
be any equilibrium distribution of $(\ref{SDEaut})$ with generator $(\ref{generator})$. Let $\vec X=(X_1, X_2)$ be a random variable
with distribution $\Gamma_{\vec\theta}^{c,g}$. Then:
\begin{itemize}
\item[$(i)$]
For any $f(\vec x) \in C_b^2([0,\infty)^2)$ that differs from a function with compact
support by only a constant,
\be{eq:testfunction}
\E^{c,g}_{\vec\theta}\left[(L^{c,g}_{\vec\theta} f)(\vec X)\right]
= \E^{c,g}_{\vec\theta}\left[c \sum_{i=1}^2 (\theta_i-X_i)
\frac{\partial}{\partial x_i} f(\vec X)
+ \sum_{i=1}^2 g_i(\vec X) \frac{\partial^2}{\partial x_i^2} f(\vec X)\right] = 0.
\end{equation}

\item[$(ii)$]
For all $g\in\cH_a$ with $0 \leq a < c$, all $\vec\theta\in [0,\infty)^2$ and $i=1,2$,
\begin{eqnarray}
&& \E^{c,g}_{\vec\theta} [X_i] = \theta_i, \label{eq:1stmoment}\\
&& \E^{c,g}_{\vec\theta} [X_1 X_2] = \theta_1 \theta_2, \label{eq:mixmoment}\\
&& \E^{c,g}_{\vec\theta} [X^2_i] = \theta_i^2
+ \frac{1}{c} \E^{c,g}_{\vec\theta} [g_i(\vec X)]
= \theta^2_i + \frac{1}{c}(F_c g)_i(\vec\theta), \label{eq:2ndmoment}
\end{eqnarray}
where all expectations are finite.

\item[$(iii)$]
Let $g\in\cH_a$ with $0 \leq a < c$, and let $K$ be any compact subset of $[0,\infty)^2$.
Then
\be{eq:2ndlog}
\sup_{c'\geq c,\vec\theta\in K}\
\E^{c'\!,g}_{\vec\theta} \big[(X_1+X_2+2)^2 \log (X_1+X_2+2)\big] < C_{c,K,g}
\ee
for some $C_{c,K,g}<\infty$ depending only on $c$, $K$ and $g$. Consequently,
$g_1$ and $g_2$ are uniformly integrable with respect to
$\{\Gamma_{\vec\theta}^{c',g}\}_{c'\geq c,\, \vec\theta\in K}$.
\end{itemize}
\ep

\bpr
(i) This part follows from the observation that, with our choice of $f$,
\begin{equation}
f(\vec X(t)) - f(\vec X(0)) - \int_0^t(L_{\vec\theta}^{c,g} f)(\vec X(s))ds
\end{equation}
is a martingale. Taking expectation and noting the stationarity of the distribution
of $\vec X(t)$, we obtain (\ref{eq:testfunction}).

\medskip\noindent
(ii) We first prove that the expectations in (\ref{eq:1stmoment}--\ref{eq:2ndmoment})
are all finite. Once this is settled, the equalities will follow easily.

\medskip\noindent
\underline{Finiteness}:
Let $h\in C^2_b([0,\infty))$ be such that $h(r)=r$ for $r\in [0,1]$, $h$ is constant
on $[3,\infty)$, $h'\in [0,1]$ and $h''\in [-1,0]$. Let $h_n(r) = n h(\frac{r}{n})$.
Then $h_n'\in[0,1]$, $h_n''\in[-\frac{1}{n},0]$, and $h_n(r)\uparrow r$, $h_n'(r)
\uparrow 1$, $h_n''(r)\to 0$ as $n\to\infty$.

\medskip\noindent
\underline{(\ref{eq:1stmoment})}:
We apply (\ref{eq:testfunction}) for $f(x_1,x_2) = h_n(\rho_1x_1+\rho_2x_2)$
with fixed $\rho_1,\rho_2>0$. Since (in the formulas below we suppress the
argument)
\be{parth1}
\begin{aligned}
\partial_{x_i} h_n(\rho_1x_1+\rho_2x_2) &= \rho_i h'_n,\\
\partial^2_{x_i} h_n(\rho_1x_1+\rho_2x_2) &= \rho_i^2 h_n'',
\end{aligned}
\ee
and $h_n(\rho_1x_1+\rho_2x_2)$ differs from a function with compact support by
a constant, by substituting the partials into (\ref{eq:testfunction}), we get
\be{Ezero1}
\E^{c,g}_{\vec\theta} \left[c \sum_{i=1}^2 \rho_i (\theta_i-X_i)h_n'
+ \sum_{i=1}^2 \rho_i^2 g_i(\vec X) h_n'' \right]=0,
\ee
which can be rewritten as
\be{eq:1st}
\begin{aligned}
c\,\E^{c,g}_{\vec\theta}\Big[(\rho_1 X_1+\rho_2 X_2)h_n'\Big]
&= c\,\E^{c,g}_{\vec\theta}\Big[(\rho_1\theta_1+\rho_2\theta_2)h_n'\Big]
+ \E^{c,g}_{\vec\theta}\Big[(\rho_1^2g_1+\rho_2^2g_2)h_n''\Big]\\
&\leq c\,\E^{c,g}_{\vec\theta}\Big[\rho_1\theta_1+\rho_2\theta_2\Big]
\end{aligned}
\ee
since $h_n''\leq 0$ and $g_1, g_2\geq 0$. By monotone convergence as $n\to\infty$, we get
\be{mon}
\rho_1 \E^{c,g}_{\vec\theta}[X_1] + \rho_2 \E^{c,g}_{\vec\theta}[X_2]
\leq \rho_1\theta_1 + \rho_2\theta_2.
\ee
Since $\rho_1,\rho_2$ are arbitrary, we obtain $\E^{c,g}_{\vec\theta} [ X_i ] \leq
\theta_i <\infty$, $i=1,2$.

\medskip\noindent
\underline{(\ref{eq:mixmoment})}:
Here we apply (\ref{eq:testfunction}) for $f(x_1,x_2) = h_n\big((1+x_1)(1+x_2)\big)$.
The calculations are similar to that for (\ref{eq:1stmoment}), which we skip.

\medskip\noindent
\underline{(\ref{eq:2ndmoment})}:
Here we apply (\ref{eq:testfunction}) for $f(x_1,x_2) = h_n(\rho_1x_1^2+\rho_2x_2^2)$
with fixed $\rho_1,\rho_2>0$. Since
\be{parth3}
\begin{aligned}
\partial_{x_i} h_n(\rho_1x_1^2+\rho_2x_2^2) &= 2\rho_i x_i h_n',\\
\partial^2_{x_i} h_n(\rho_1x_1^2+\rho_2x_2^2) &= 2\rho_i h_n'
+ 4\rho_i^2 x^2_i h_n'',
\end{aligned}
\ee
by substituting the partials into (\ref{eq:testfunction}), we get
$$
\E^{c,g}_{\vec\theta}\Big[2 c\rho_1X_1 (\theta_1-X_1)h_n'
+ 2c\rho_2 X_2(\theta_2-X_2)h_n'+ 2(\rho_1g_1+\rho_2g_2)h_n'
+ 4(\rho_1^2 X_1^2g_1 + \rho_2^2 X_2^2g_2)h_n''\Big] = 0.
$$
Rearranging terms, we obtain
\be{eq:2nd}
\begin{aligned}
&2c\,\E^{c,g}_{\vec\theta}\Big[(\rho_1 X_1^2+\rho_2 X_2^2)h_n'\Big]\\
&\qquad = 2c\,\E^{c,g}_{\vec\theta}\Big[(\rho_1\theta_1 X_1+\rho_2\theta_2 X_2)h_n'\big]
+ 2\E^{c,g}_{\vec\theta}\Big[(\rho_1g_1+\rho_2g_2)h_n'\Big]\\
&\qquad\qquad + \E_{\vec\theta}^{c,g}\Big[ 4(\rho_1^2 X_1^2g_1
+ \rho_2^2 X_2^2g_2)h_n'' \Big]\\
&\qquad \leq 2c\,(\rho_1\theta_1^2 + \rho_2\theta_2^2)
+ 2\, \E^{c,g}_{\vec\theta}\Big[(\rho_1g_1+\rho_2g_2)h_n'\Big].
\end{aligned}
\ee
Since $g\in\cH_a$ with $0\leq a <c$, we have $g_1(\vec x)+g_2(\vec x) \leq C(1+x_1)(1+x_2)
+ a(x_1^2+x_2^2)$. Substituting this bound into (\ref{eq:2nd}) and setting $\rho_1 = \rho_2 =1$,
using the fact that $\E^{c,g}_{\vec\theta}[X_i]\leq\theta_i$ and $\E^{c,g}_{\vec\theta}
[X_1 X_2]<\infty$, and rearranging terms, we get
\be{Erel}
2(c-a) \,\E^{c,g}_{\vec\theta}\Big[(X_1^2+X_2^2)h_n'\Big] < C'<\infty.
\ee
By monotone convergence as $n\to\infty$, we obtain $\E^{c,g}_{\vec\theta}[ X_i^2]<\infty$.
This also implies $\E^{c,g}_{\vec\theta}[g_i]<\infty$.

\medskip\noindent
\underline{Equality}:
Having thus proved that the expectations in (\ref{eq:1stmoment}--\ref{eq:2ndmoment})
are finite, we are now ready to prove that equality holds. To that end, return to
(\ref{eq:1st}). Since $\E^{c,g}_{\vec\theta}[\rho_1^2 g_1+\rho_2^2 g_2]<\infty$,
$h_n''\in[-\frac{1}{n},0]$ and $h_n''\to 0$ as $n\to\infty$, (\ref{eq:1stmoment})
follows by applying the dominated convergence theorem. By the same argument,
(\ref{eq:2ndmoment}) follows by applying the dominated convergence theorem to
(\ref{eq:2nd}), provided that
\begin{equation}
\label{ubquadr}
(\rho_1^2x_1^2g_1 +\rho_2^2x_2^2g_2)
|h_n''(\rho_1x_1^2+\rho_2x_2^2)| \leq C(x_1^2+x_2^2)
\end{equation}
for some $C<\infty$ independent of $n$. To see the latter, note that $h_n''(r) =
\frac{1}{n} h''(\frac{r}{n}) = \frac{r}{n}h''(\frac{r}{n})
\frac{1}{r} \leq \frac{3}{r}$, since $h''\in[-1,0]$ and $h''(\frac{r}{n})\neq 0$
only when $\frac{r}{n}\leq 3$. The bound in (\ref{ubquadr}) then follows readily.

To verify (\ref{eq:mixmoment}), we apply (\ref{eq:testfunction}) for $f(x_1, x_2) = h_n((x_1+x_2)^2)$
instead of $h_n((1+x_1)(1+x_2))$. This gives
\be{extra}
\E^{c,g}_{\vec\theta}\Big[ 2c(\theta_1+\theta_2-X_1-X_2)(X_1+X_2)h_n'
+ 2(g_1+g_2)h_n' + 4(X_1+X_2)^2(g_1+g_2)h_n''\Big] = 0.
\ee
Since $(x_1+x_2)^2 h_n''((x_1+x_2)^2) \leq 3$ and $\E_{\vec\theta}^{c,g}[g_1+g_2]
<\infty$, we can apply the dominated convergence theorem in (\ref{extra}) as
$n\to\infty$. Then, together with (\ref{eq:1stmoment}) and (\ref{eq:2ndmoment}),
we obtain (\ref{eq:mixmoment}).

\medskip\noindent
(iii) This part follows from similar computations as in part (ii). Let $c'\geq c$
be arbitrary, and abbreviate $\bar X_i = 1+X_i$, $\bar\theta_i = 1+\theta_i$, $\bar x_i
= 1+x_i$ for $i=1,2$. We first show that
\be{bdno1}
\E_{\vec\theta}^{c'\!,g} [ \bar X_1 \bar X_2\log(\bar X_1 + \bar X_2)]<\infty
\ee
by applying (\ref{eq:testfunction}) to $h_n\big(\bar x_1 \bar x_2 \log(\bar x_1
+\bar x_2)\big)$. Then we apply (\ref{eq:testfunction}) to $h_n\big((\bar x_1+\bar x_2)^2\log(\bar x_1
+\bar x_2)\big)$ to prove (\ref{eq:2ndlog}).

\medskip\noindent
\underline{(\ref{bdno1})}:
Let $f(x_1, x_2) = h_n\big(\bar x_1 \bar x_2 \log(\bar x_1 + \bar x_2)\big)$, which
differs from a function with compact support by a constant. Since
$$
\begin{aligned}
\partial_{x_1} h_n\big(\bar x_1\bar x_2 \log(\bar x_1 +\bar x_2)\big)
&= \Big(\bar x_2\log(\bar x_1+\bar x_2)
+ \frac{\bar x_1 \bar x_2}{\bar x_1+\bar x_2}\Big) h_n' ,\\
\partial^2_{x_1} h_n\big(\bar x_1\bar x_2 \log(\bar x_1 + \bar x_2)\big)
&= \Big(\frac{\bar x_2}{\bar x_1 + \bar x_2}
+\frac{\bar x_2^2}{(\bar x_1+\bar x_2)^2}\Big)h_n'
+ \Big(\bar x_2\log(\bar x_1+\bar x_2)
+ \frac{\bar x_1 \bar x_2}{\bar x_1+\bar x_2}\Big)^2 h_n'',
\end{aligned}
$$
and since the same holds if we interchange the indices 1 and 2, by substituting the partials
into (\ref{eq:testfunction}) and noting that $\partial^2_{x_1} h_n\leq 2h_n',
\partial^2_{x_2}h_n \leq 2h_n'$, we get
\be{Ezero4}
\begin{aligned}
\E^{c'\!,g}_{\vec\theta}\Big[\ & c' (\bar\theta_1-\bar X_1)
\Big(\bar X_2\log(\bar X_1+\bar X_2)
+ \frac{\bar X_1\bar X_2}{\bar X_1+\bar X_2}\Big) h_n' \\
& \ + c' (\bar\theta_2-\bar X_2)
\Big(\bar X_1\log(\bar X_1+\bar X_2)
+ \frac{\bar X_1\bar X_2}{\bar X_1+\bar X_2}\Big) h_n' + 2(g_1+g_2)h_n' \Big] \geq 0.
\end{aligned}
\ee
Rearranging terms and noting that $\frac{\bar X_1 \bar X_2}{\bar X_1+\bar X_2}
< \bar X_1 \wedge \bar X_2$, we find that
\be{eq:4th}
2c'\, \E_{\vec\theta}^{c'\!,g} [ \bar X_1\bar X_2\log (\bar X_1+\bar X_2) h_n' ]
\leq \E_{\vec\theta}^{c'\!,g}\Big[ c' (\bar\theta_1\bar X_2 +\bar X_1\bar\theta_2)
(1+\log (\bar X_1+\bar X_2)) + 2(g_1+ g_2)\Big].
\ee
By assumption, $g_1(\vec x)+g_2(\vec x) \leq C(1+x_1)(1+x_2) + a (x_1^2+x_2^2)$.
Substituting this bound into (\ref{eq:4th}), applying monotone convergence as
$n\to\infty$, and noting that (\ref{eq:2ndmoment}) implies that
$\E_{\vec\theta}^{c'\!,g}[ X_1^2+X_2^2] \leq \phi(\theta_1,\theta_2)$
for some quadratic polynomial $\phi$ depending only on $c$ and $g$, we
easily verify that
\be{seclog}
\E_{\vec\theta}^{c'\!,g} [ \bar X_1 \bar X_2 \log (\bar X_1 +\bar X_2) ]
\leq \tilde\phi(\theta_1,\theta_2)
\ee
for some cubic polynomial $\tilde\phi$ depending only on $c$ and $g$.

By applying (\ref{eq:testfunction}) to $h_n\big( (\bar x_1+\bar x_2)^2\log(\bar x_1+\bar x_2) \big)$
and using (\ref{seclog}), it can be shown that
\be{eq:7th}
\E_{\vec\theta}^{c'\!,g} [ (\bar X_1+\bar X_2)^2\log(\bar X_1 +\bar X_2)]
\leq \hat\phi(\theta_1,\theta_2)
\ee
for some cubic polynomial $\hat\phi$ depending only on $c$ and $g$. The uniform bound
in (\ref{eq:2ndlog}) then follows. The calculations, which we omit, are similar as before.

Since $g_1(\vec x)+g_2(\vec x) \leq C(\bar x_1^2+\bar x_2^2)$ for some $C<\infty$,
which by (\ref{eq:7th}) is uniformly integrable with respect to
$\{\Gamma_{\vec\theta}^{c'\!,g}\}_{c'\geq c,\, \vec\theta\in K}$ for any compact
$K\subset [0,\infty)^2$, it follows that $g_1$ and $g_2$ are also uniformly integrable.
\epr
\smallskip

\noindent
{\bf Remark.} By similar computations, it can be shown that (\ref{eq:2ndlog}) is still valid when the
logarithm in the left-hand side of the inequality is raised to an arbitrary power.


\section{Appendix 2: Properties of uniformly elliptic diffusions}
\label{A2}

In this Appendix, we list some facts about uniformly elliptic diffusions that are needed in the proof of
Theorem \ref{thm:SDEeq}. We thank S.R.S.~Varadhan for pointing out some of the relevant results and references
on uniformly elliptic diffusions.

\bt{thm:diffprop}{\bf [Uniformly elliptic diffusions in $\R^d$]} \\
Let $b : \R^d\to \R^d$ be a bounded measurable map, and let $a: \R^d \to S_d$ be a continuous map, where $S_d$ is the space of symmetric non-negative definite
$d\times d$ real matrices. Assume further that $a(\cdot)$ is uniformly elliptic, i.e., there exists $0<\Lambda<\infty$ such that for all
$\vec x, \vec \theta \in\R^d$, $\vec\theta\neq 0$,
$$
\Lambda^{-1} \leq \frac{\langle \vec\theta, a(\vec x)\vec\theta\rangle}{\langle\vec\theta, \vec\theta\rangle} \leq \Lambda.
$$
Then, for each $\vec x\in\R^d$, the martingale problem with generator
\be{eq:appendixL}
Lf = \sum_{i,j=1}^d a_{ij}(\vec x)\frac{\partial^2}{\partial x_i \partial x_j}f(\vec x) + \sum_{i=1}^d b_i(\vec x) \frac{\partial}{\partial x_i} f(\vec x), \qquad f\in C_c^2(\R^d),
\ee
has a unique solution $\P^{\vec x}$ in the space of probability measures on $\Omega=C([0,\infty), \R^d)$ with $\P^{\vec x}(\omega\in \Omega : \omega(0)=\vec x)=1$.
The family of solutions $\{\P^{\vec x}\}_{\vec x\in\R^d}$ defines a strong Feller and strong Markov process
that admits a transition probability density $p_t(\vec x, \vec y)$ with respect to Lebesgue measure for each $t> 0$ and $\vec x\in\R^d$. Furthermore, for each
$t>0$ and $\vec x^* \in\R^d$,
$$
\lim_{\vec x\to \vec x^*} \Vert p_t(\vec x, \cdot)-p_t(\vec x^*, \cdot)\Vert_1 = \lim_{\vec x\to \vec x^*} \int_{\R^d} |p_t(\vec x, \vec y) - p_t(\vec x^*, \vec y)| d\vec y = 0.
$$
\et
{\bf Proof.} All facts follow from results in Stroock and Varadhan~\cite{SV79}. For the well-posedness of the martingale
problem, see Theorem 7.2.1 therein. For the strong Markov property, see Theorem 6.2.2. For the strong Feller property, see Theorem 7.2.4. For the existence of the transition density, see Theorem 9.1.9 and Lemma 9.2.2. Lastly, for the $L_1$-continuity of the transition density, see Theorem 11.4.3.
\qed

\bt{thm:domaindiff}{\bf [Diffusions restricted to bounded domains]}\\
Let $a$ and $b$ satisfy the conditions in Theorem {\rm \ref{thm:diffprop}}, and let $\{\P^{\vec x}\}_{\vec x\in\R^d}$ denote the family of solutions to
the martingale problem with coefficients $(a,b)$ in {\rm (\ref{eq:appendixL})}. If $\bar a : \R^d\to S_d$ and $\bar b: \R^d \to \R^d$ are locally bounded measurable maps
with $\bar a=a$ and $\bar b=b$ on a bounded open set $D$, then for any $\vec x\in D$ and any solution $\bar \P^{\vec x}$ to the martingale problem with
coefficients $(\bar a, \bar b)$, $\bar \P^{\vec x} = \P^{\vec x}$ on ${\cal F}_{\tau_D}$, the sigma-field on $\Omega$ generated by the family of projection maps
$\{\pi_s : \Omega\to\R^d\ |\ \pi_s(\omega)=\omega(s\wedge\tau_D)   \}_{s\geq 0}$, where $\tau_D(\omega) = \inf\{t\geq 0 : \omega(t)\notin D\}$.
\et
{\bf Proof.} See Theorem 10.1.1 in Stroock and Varadhan~\cite{SV79}.

\bc{cor:domaindiffdensity}{\bf [Transition density for diffusions restricted to bounded domains]}\\
Let $b : \R^d\to\R^d$ be a locally bounded measurable map, and let $a : \R^d\to S_d$ be continuous such that the martingale problem with coefficients $a$ and $b$
in {\rm (\ref{eq:appendixL})} is well-posed. Assume further that $a$ is non-degenerate on $\overline D$ for a simply connected bounded open set $D\subset \R^d$ with
smooth boundary. For any $\vec x\in D$, if $\P^{\vec x}$ is the solution of the martingale problem starting from $\vec x$, then, for each $t>0$, the measure
$\mu_t^D(\vec x, \cdot)$ on Borel-measurable sets defined by $\mu_t^D(\vec x, \cdot) = \P^{\vec x}(\omega : t<\tau_D(\omega), \omega(t) \in \cdot)$ admits a
density $p_t^D(\vec x, \vec y)$ with respect to Lebesgue measure. Furthermore, for each $\vec x^*\in D$, there exist $\epsilon, \delta>0$ sufficiently small such
that, for all $\vec x, \vec x' \in B_\epsilon(\vec x^*)$, the ball of radius $\epsilon$ centered at $\vec x^*$, the overlap between $\mu_\delta^D(\vec x, \cdot)$
and $\mu_\delta^D(\vec x', \cdot)$ satisfies
\be{eq:densityoverlap}
\frac{\mu_\delta^D(\vec x, D) + \mu_\delta^D(\vec x', D) - \Vert p_\delta^D(\vec x, \cdot)- p_\delta^D(\vec x', \cdot)\Vert_1}{2} \geq \frac{1}{2}.
\ee
\ec
{\bf Proof.} By our assumptions on $a, b$ and $D$, we can find coefficients $(\bar a, \bar b)$ on $\R^d$ such that $(\bar a, \bar b)=(a,b)$ on $D$, $(\bar a,\bar b)$
are bounded, $a$ is continuous and uniformly elliptic on $\R^d$. For instance, we can define $\bar b=b$ on $D$ and $\bar b\equiv 0$ on $\R^d \backslash D$,
define $\bar a =a$ on $\overline D$ and $\bar a\equiv I$ on $\R^d \backslash B$ where $B$ is a large open ball containing $\overline D$, and on $B\backslash \overline D$ define $\bar a$ to be the
harmonic interpolation between its values on $\partial B$ and $\partial D$. By Theorem \ref{thm:diffprop}, the martingale problem with coefficients $(\bar a, \bar b)$
has a unique family of solutions $\{\bar P^{\vec x}\}_{\vec x\in\R^d}$, which is strong Markov and admits a transition density $\bar p_t(\vec x, \vec y)$ for all $t>0$
and $\vec x\in\R^d$. By Theorem \ref{thm:domaindiff}, for $\vec x\in D$, $\bar \P^{\vec x} = \P^{\vec x}$ on ${\cal F}_{\tau_D}$. In particular,
$\mu_t^D(\vec x, \cdot) = \bar \mu_t^D(\vec x, \cdot) = \bar \P^{\vec x}\{ \omega : t<\tau_D(\omega), \omega(t) \in \cdot\}$. Since $\bar \mu_t^D(\vec x, \cdot)$
is absolutely continuous with respect to $\bar \P^{\vec x}(\omega : \omega(t) \in \cdot)$ with density $\bar p_t(\vec x,\vec y)$,
$\mu_t^D(\vec x, \cdot)=\bar\mu_t^D(\vec x,\cdot)$ also admits a density $p^D_t(\vec x, \vec y)$ with respect to Lebesgue measure for all $\vec x\in D$ and $t>0$.

It is not difficult to see that the left-hand side of (\ref{eq:densityoverlap}) is the mass of the maximal positive
measure that is dominated by both $\mu^D_\delta(\vec x, \cdot)$ and $\mu^D_\delta(\vec x', \cdot)$. To verify (\ref{eq:densityoverlap}),
fix $\vec x^* \in D$ and choose $\epsilon'>0$ such that $B_{2\epsilon'}(\vec x^*) \subset D$. Then we can choose $\delta>0$
sufficiently small such that, for all $\vec x\in B_{\epsilon'}(\vec x^*)$, $\P^{\vec x}(\tau_D\leq \delta) \leq \frac{1}{5}$. To verify this claim, note that,
given $\vec z\in B_{\epsilon'}(\vec x^*)$, if we define $f(\vec x) = \Vert \vec x-\vec z\Vert^2 = \sum_{i=1}^d (x_i - z_i)^2$, then
$$
f(\vec X(t\wedge \tau_D)) - f(\vec X(0)) - \int_0^{t\wedge \tau_D} Lf(\vec X(s)) ds
$$
is a martingale, where $(\vec X(s))_{s\geq 0}$ has law $\P^{\vec z}$. In particular,
\begin{eqnarray}
\label{eq:escapeprob}
(\epsilon')^2 \P^{\vec z}(\tau_D \leq \delta) &\leq& \E[\Vert \vec X(\delta \wedge \tau_D) - \vec z\Vert^2] \\
&=& \E\left[ \int_0^{\delta\wedge \tau_D} 2\sum_{i=1}^d \left(b_i(\vec X(s)) (X_i(s) - z_i) + a_{ii}(\vec X(s))\right) ds \right]
\leq \delta\ C_{D, a, b}, \nonumber
\end{eqnarray}
where $C_{D,a,b}$ depends only on $D$ and $(a,b)$ on $D$. Therefore $\P^{\vec z}(\tau_D \leq \delta) \leq \delta\, C_{D,a,b} (\epsilon')^{-2}$ uniformly for
all $\vec z \in B_{\epsilon'}(\vec x^*)$. Choosing $\delta$ sufficiently small, we then verify the claim.

Applying Theorem \ref{thm:diffprop} to $\{\bar P^{\vec x}\}_{\vec x\in\R^d}$, we can choose $\epsilon\in(0,\epsilon')$ small such that, for all
$\vec x \in B_{\epsilon}(\vec x^*)$, $\Vert \bar p_\delta(\vec x, \cdot) - \bar p_\delta(\vec x^*, \cdot)\Vert_1 \leq \frac{1}{10}$, and hence, for all
$\vec x, \vec x' \in B_{\epsilon}(\vec x^*)$, $\Vert\bar p_\delta(\vec x, \cdot) - \bar p_\delta(\vec x', \cdot)\Vert_1 \leq \frac{1}{5}$. Since
for $\vec z\in B_{\epsilon}(\vec x^*)$, $\Vert \bar p_\delta(\vec z, \cdot) - \bar p_\delta^D(\vec z, \cdot)\Vert_1 = \bar\P^{\vec z}(\tau_D\leq \delta)\leq \frac{1}{5}$,
we have $\Vert \bar p_\delta^D(\vec x, \cdot) - \bar p_\delta^D(\vec x', \cdot)\Vert_1 \leq \frac{3}{5}$ for all $\vec x, \vec x'\in B_\epsilon(\vec x^*)$.
Finally, note that, for $\vec x, \vec x'\in B_\epsilon(\vec x^*)$, $\mu_\delta^D (\vec x, D) = 1 - \bar P^{\vec
  x}(\tau_D \leq \delta) \geq 1- \frac{1}{5}$ and the same holds for $\mu_\delta^D(\vec x', D)$, hence, substitution of all
the estimates into the left-hand side of (\ref{eq:densityoverlap}) yields the desired result. \qed
\smallskip

\noindent{\bf Remark.} Note that the constant on the right-hand side of (\ref{eq:densityoverlap}) can be made arbitrarily close to 1 by choosing $\epsilon, \delta$ sufficiently small.

\bt{thm:supportthm}{\bf [Support theorem for uniformly elliptic diffusions]}\\
Let $a, b, D$ and $\{P^{\vec x}\}_{\vec x\in D}$ be as in Corollary {\rm \ref{cor:domaindiffdensity}}. For any $\vec x\in D$, $\epsilon>0$, and any continuous function
$\psi : [0,t]\to D$ with $\psi(0)=\vec x$,
$$
\P^{\vec x}\Big(\omega : \sup_{0\leq s\leq t} |\omega(s)-\psi(s)| \leq \epsilon\Big)>0.
$$
\et
{\bf Proof.} The support theorem is a classic result of Stroock and Varadhan. The statement above follows Theorem (2.5) in Chapter V of Bass~\cite{Bass98} and
Theorem \ref{thm:domaindiff} above.

\bt{thm:occupation}{\bf [Occupation time measure for uniformly elliptic diffusions]}\\
Let $a$, $b$, $D$ and $\{\P^{\vec x}\}_{\vec x\in D}$ be as in Corollary {\rm \ref{cor:domaindiffdensity}}. If $A\subset D$ has positive Lesbegue measure,
then, for all $\vec x\in D$, $\E^{\vec x}[\int_0^{\tau_D} 1_{\omega(s) \in A} ds] >0$, where $\E^{\vec x}$ denotes expectation with respect to
$\P^{\vec x}$, and $\tau_D = \inf\{t\geq 0 : \omega(t) \notin D\}$.
\et
{\bf Proof.} The statement above follows from Theorem (8.5) in Chapter V of Bass~\cite{Bass98} (which goes back to Krylov) in combination with the support theorem, Theorem \ref{thm:supportthm}, and the Girsanov transformation (see Theorem 7.2.2 in Stroock and Varadhan~\cite{SV79}).
\qed
\bigskip

\noindent
{\bf Acknowledgment:} The work in this paper was supported by
DFG and NWO, as part of the Dutch-German Bilateral Research Group
on ``Mathematics of Random Spatial Models from Physics and Biology''.
AG and JMS were supported by the DFG-grant GR 876/12--1 -- 12--3.
DD was hosted by EURANDOM during two visits and is supported by an NSERC Discovery Grant.
JMS is supported by GA\v CR grant 201/06/1323. JMS and RS received travel support from the ESF scientific program
``Random Dynamics in Spatially Extended Models''. FdH and RS are grateful to
the Pacific Institute for the Mathematical Sciences and the Mathematics
Department of the University of British Columbia, Vancouver, Canada,
for hospitality: FdH from January to August 2006, RS from mid-April
to mid-May 2006 when part of the work in this paper was completed.
RS was a postdoc at EURANDOM from October 2004 to October 2006.
DD, FdH and RS thank Ed Perkins for valuable discussions. The authors
thank the associate editor and the referee for an exceptionally careful reading
of the paper and many helpful suggestions.


\end{document}